\documentclass[reqno,fleqn]{amsart}
\usepackage{amsfonts,amsthm,amssymb}
\usepackage{lipsum}
\usepackage[utf8]{inputenc}
\usepackage[english]{babel}

\usepackage{lmodern,textcomp}
\usepackage{amsmath}
\usepackage{mathtools}
\usepackage{latexsym}
\usepackage{tikz}
\usepackage{fancyvrb}
\usepackage{epsfig}
\usepackage{pstricks,slashbox,multirow}
\usepackage{graphicx}
\usepackage{rotating}
\usetikzlibrary{positioning}
\usepackage{subcaption}
\usepackage{datetime}
\newtheorem{theorem}{Theorem}

\newtheorem{definition}[theorem]{Definition}

\newtheorem{prm}{Problem}
\newtheorem{oprm}{Open Problem}

\newtheorem{rem}[theorem]{Remark}

\title[A study on Type-2 isomorphic $C_n(R)$: Part 5: Type-2 isomorphic $C_{n}(R)$ for $n$ = 48,81,96]{A study on Type-2 isomorphic circulant graphs. \\ Part 5: Type-2 isomorphic circulant graphs of orders 48, 81, 96}

\author{\sc Vilfred Kamalappan} 
\address{Department of Mathematics, Central University of Kerala, Periye, Kasaragod, Kerala, India - 671 316.}
\email{vilfredkamal@gmail.com}



\date{}

\begin{document}

\begin{abstract} This article is the $5^{th}$ part of a detailed study on Type-2 isomorphic circulant graphs having ten parts \cite{v2-1}-\cite{v2-10} and is a continuation of Part 4. Here, we study Type-2 isomorphic circulant graphs of $C_{48}(r_1,r_2,r_3)$,  $C_{81}(r_1,r_2,r_3)$  and $C_{96}(r_1,r_2,r_3,r_4)$. We find that for $n$ = 48,96, total number of pairs of isomorphic circulant graphs of Type-2 w.r.t. $m$ = 2 of the forms  $C_{n}(r_1,r_2,r_3)$ and $C_{n}(s_1,s_2,s_3)$ are 18 and 72, respectively, and total number of triples of isomorphic circulant graphs of Type-2 w.r.t. $m$ = 3 of the form $C_{81}(x_1,x_2,x_3)$, $C_{81}(y_1,y_2,y_3)$ and $C_{81}(z_1,z_2,z_3)$ are 27.
\end{abstract}

\maketitle

	
\section{Introduction}

In \cite{v2-2} to \cite{v2-4}, we studied Type-2 isomorphic circulant graphs of orders 16, 24, 27, 32 and 54 and shown that total number of pairs of Type-2 isomorphic circulant graphs of orders 16, 24 and 32 are 8, 32 and 384, respectively and total number of triples of Type-2 isomorphic circulant graphs of orders 27 and 54 are 12 and 960, respectively.  This paper is a continuation of \cite{v2-4} and using modified definition \ref{d3}, we study Type-2 isomorphic circulant graphs of $C_{48}(r_1,r_2,r_3)$,  $C_{81}(r_1,r_2,r_3)$  and $C_{96}(r_1,r_2,r_3,r_4)$. We find that the total number of pairs of Type-2 isomorphic circulant graphs of the forms  $C_{n}(r_1,r_2,r_3)$ and $C_{n}(s_1,s_2,s_3)$ are 18 and 72 for $n$ = 48, 96, respectively and the total number of triples of Type-2 isomorphic circulant graphs of the form $C_{81}(x_1,x_2,x_3)$, $C_{81}(y_1,y_2,y_3)$ and $C_{81}(z_1,z_2,z_3)$ are 27. For basic definitions and results on isomorphic circulant graphs, refer \cite{v2-1, v2-4}.

\begin{definition}{\rm\cite{ad67}} \quad \label{d1} For $R =$ $\{r_1$, $r_2$, $\ldots$, $r_k\}$ and $S$ = $\{s_1$, $s_2$, $\ldots$, $s_k\}$, circulant graphs $C_n(R)$ and $C_n(S)$ are {\it Adam's isomorphic} if there exists a positive integer $x$ $\ni$ $\gcd(n, x)$ = 1 and $S$ = $\{xr_1$, $xr_2$, $\ldots$, $xr_k\}_n^*$ where $<r_i>_n^*$, the {\it reflexive modular reduction} of a sequence $< r_i >$, is the sequence obtained by reducing each $r_i$ under modulo $n$ to yield $r_i'$ and then replacing all resulting terms $r_i'$ which are larger than $\frac{n}{2}$ by $n-r_i'$. Hereafter, we call Adam's isomorphism as {\em Type-1 isomorphism}.  
\end{definition}

A circulant graph $C_n(R)$ is said to have {\em Cayley Isomorphism (CI) property} if whenever $C_n(S)$ is isomorphic to $C_n(R)$, they are Type-1 isomorphic \cite{v2-1}.

\begin{theorem} \cite{v24} \label{t2} Let $Ad_n(C_n(R))$ = $\{\varphi_{n,x}(C_n(R)) = C_n(xR): x\in\varphi_n \}$. Then, $C_n(S)\in Ad_n(C_n(R))$ if and only if $Ad_n(C_n(R))$ = $Ad_n(C_n(S))$ if and only if $C_n(R)\in Ad_n(C_n(S))$. \hfill $\Box$
\end{theorem}

Definition of Type-2 isomorphism of circulant graph $C_n(R)$ w.r.t. $m$ given in \cite{v2-2-arX} is modified in \cite{v2-1} as follows and hereafter we use the same. 

\begin{definition} \cite{v2-1} \quad  \label{d3} Let $V(K_n) = \{u_0,u_1,u_2,...,u_{n-1}\}$, $V(C_n(R))$ = $\{v_0, v_1, v_2, ... , v_{n-1}\}$, $|R| \geq 3$, $r\in R$ and $m > 1$ and $m^3$ be divisors of $\gcd(n, r)$ and $n$, respectively. Define 1-1 mapping $\theta_{n,m,t} :$ $V(C_n(R)) \rightarrow V(K_n)$ such that $\theta_{n,m,t}(v_x) = u_{x+jtm}$,  $\theta_{n,m,t}((v_x, v_{x+s}))$ = $(\theta_{n,m,t}(v_x),$ $\theta_{n,m,t}(v_{x+s}))$ under subscript arithmetic modulo $n$ and $\theta_{n,m,t}(C_n(R))$ = $C_n(\theta_{n,m,t}(R))$ where $\theta_{n,m,t}(R)$ in $C_n(\theta_{n,m,t}(R))$ is calculated under the reflexive modulo $n$, $\forall$ $x \in \mathbb{Z}_n$, $x = qm+j,$ $0 \leq j \leq m-1$, $s\in R$ and $0 \leq q,t \leq \frac{n}{m} -1$. And for a particular value of $t,$ if  $\theta_{n,m,t}(C_n(R))$ = $C_n(S)$ for some $S$  and  $S \neq yR$ for all $y\in \varphi_n$ under reflexive modulo $n,$ then $C_n(R)$ and $C_n(S)$ are called {\em isomorphic circulant graphs of Type-2 w.r.t. $m$.} 
	
When $C_n(R)$ and $C_n(S)$ are Type-2 isomorphic w.r.t. $m$, then we also say that $C_{kn}(kR)$  and $C_{kn}(kS)$ are Type-2 isomorphic w.r.t. $m$, $k\in\mathbb{N}$. Here, $k.C_n(T)$ = $C_{kn}(kT)$, $k\in\mathbb{N}$. 	 
\end{definition}

\begin{rem} \cite{v2-1}  \label{r4} Following steps are used to establish Type-2 isomorphism w.r.t. $m$ between circulant graphs $C_n(R)$ and $C_n(S)$. (i) $R$ $\neq$ $S$ and $|R| = |S| \geq 3$; (ii) $\exists$ $r\in R,S$ and $m > 1$ $\ni$ $m$ is a divisor of $\gcd(n, r)$, $m^3$ is a divisor of $n$ and for some $t$ $\ni$ $1 \leq t \leq \frac{n}{m} -1$, $\theta_{n,m,t}(C_n(R))$ = $C_n(S)$ and (iii) $S$ $\neq$ $xR$ for all $x\in\varphi_n$ under arithmetic reflexive modulo $n$. 

Thus, if $C_n(R)$ and $C_n(S)$ are Type-2 isomorphic circulant graphs w.r.t. $m$, then there exist $r\in R,S$, $m > 1$ and some $t$ $\ni$ $m$ is a divisor of $\gcd(n, r)$, $m^3$ is a divisor of $n$, $1 \leq t \leq \frac{n}{m} -1$, $\theta_{n,m,t}(C_n(R))$ = $C_n(S)$ and $S$ $\neq$ $xR$ for all $x\in\varphi_n$ under arithmetic reflexive modulo $n$.
\end{rem} 

\begin{rem}  \label{r5} \quad The calculation on jump sizes $r_i$s which are integer multiples of $m$ need not be done under the transformation $\theta_{n,m,t}$, while searching for possible value(s) of $t$ for which the transformed graph $\theta_{n,m,t}(C_n(R))$ is circulant of the form $C_n(S)$ for some $S \subseteq [1, \frac{n}{2}]$, as there is no change in these $r_i$s where $r\in R$ and $m > 1$ and $m^3$ are divisors of $\gcd(n, r)$ and $n$, respectively. 

Thus, if $\theta_{n,m,t}(C_n(R))$ = $C_n(S)$ for some $S$ and thereby $C_n(R)$ $\cong$ $C_n(S)$, then $\theta_{n,m,t}(C_n(R \cup mT))$ = $C_n(S \cup mT)$ for any $T$ and thereby $C_n(R \cup mT)$ $\cong$ $C_n(S \cup mT)$.

Also, for a given $C_n(R)$, w.r.t. different values of $m$, we may get different Type-2 isomorphic circulant graphs.
\end{rem}

\begin{rem} \cite{v2-1} \label{r6} \quad {\rm For given $C_n(R)$ and $C_n(S)$ when either $\theta_{n,m,t}(C_n(R))$ = $C_n(S)$ for some $t$ or $C_n(xR)$ = $C_n(S)$ for some $x$, then $C_n(R)$ and $C_n(S)$ are isomorphic, $0 \leq t \leq \frac{n}{m} -1$ and $x\in\varphi_n$. }
\end{rem}

\begin{theorem}{\rm \cite{v24}}\quad \label{t7} {\rm For $n \geq 2$, $1 \leq 2s-1 \leq 2n-1$, $n \neq 2s-1$, $R$ = $\{2,2s-1, 4n-(2s-1)\}$ and $S$ = $\{ 2,$ $2n-(2s-1)$, $2n+2s-1 \}$, $\theta_{8n,2,n}(C_{8n}(R))$ = $C_{8n}(S)$ = $\theta_{8n,2,3n}(C_{8n}(R)),$ $\theta_{8n,2,n}(C_{8n}(S))$ = $C_{8n}(R)$ = $\theta_{8n,2,3n}(C_{8n}(S))$ and circulant graphs $C_{8n}(R)$ and $C_{8n}(S)$ are Type-2 isomorphic  w.r.t. $m$ = 2. When $n$ = $2s-1$, the two circulant graphs are the same. \hfill $\Box$}
\end{theorem}

\begin{theorem} \cite{v2-6} \label{t8} {\rm Let $p$ be an odd prime number, $1 \leq i \leq p$, $1 \leq x \leq p-1$, $y\in\mathbb{N}_0$, $0 \leq y \leq np-1$, $1 \leq x+yp \leq np^2-1$, $d^{np^3, x+yp}_i = (i-1)xpn+x+yp$,  $R^{np^3, x+yp}_i$ $=$ $\{p$, $d^{np^3, x+yp}_i$, $np^2-d^{np^3, x+yp}_i$, $np^2+d^{np^3, x+yp}_i$, $2np^2-d^{np^3, x+yp}_i$, $2np^2+$ $d^{np^3, x+yp}_i,$ $3np^2-d^{np^3, x+yp}_i$, $3np^2+d^{np^3, x+yp}_i$, . . . , $(p-1)np^2$ - $d^{np^3, x+yp}_i$, $(p-1)np^2+d^{np^3, x+yp}_i$, $np^3-d^{np^3, x+yp}_i$, $np^3-p\}$ and $i,j,n,x\in\mathbb{N}$. Then, for a given set of values of $n$, $p$, $x$ and $y$, $\theta_{np^3,p,jn} (C_{np^3}(R^{np^3, x+yp}_i))$ = $C_{np^3}(R^{np^3, x+yp}_{i+j})$ and the $p$ circulant graphs $C_{np^3}(R^{np^3, x+yp}_i)$ are isomorphic of Type-2 w.r.t.  $p$, $1 \leq i,j \leq p$ where $i+j$ in $R^{np^3, x+yp}_{i+j}$ is calculated under addition modulo $p$ and $C_{np^3}(R^{np^3, x+yp}_0)$ = $C_{np^3}(R^{np^3, x+yp}_p)$. \hfill $\Box$}
\end{theorem}

\section{Main result}

In this section, we find Type-2 isomorphic circulant graphs of $C_{48}(r_1,r_2,r_3)$, $C_{81}(r_1,r_2,r_3)$ and $C_{96}(r_1,r_2,r_3)$ and give a few open problems on these type of isomorphic circulant graphs.

\subsection{On isomorphic circulant graphs $C_{48}(r_1,r_2,r_3)$ of Type-2 w.r.t. $m$ = 2}

In this subsection, we study Type-2 isomorphic circulant graphs of $C_{48}(r_1,r_2,r_3)$. We obtain all pairs of Type-2 isomorphic circulant graphs of the form $C_{48}(r_1,r_2,r_3)$ and $C_{48}(s_1,s_2,s_3)$ and the total number of such pairs are 18 and all are of Type-2 w.r.t. $m$ = 3. We start with presenting some isomorphic circulant graphs of $C_{48}(r_1,r_2,r_3)$ of Type-1 in the next two problems. 

\begin{prm} \label{p1} {\rm The following pairs of circulant graphs are Type-1 isomorphic for $s$ = 4,8,12,16,20,24.  
\begin{enumerate}
     \item [\rm (a)]	$C_{48}(1,s,23)$ and $C_{48}(s,11,13)$; 
 	
     \item [\rm (b)]	$C_{48}(3,s,21)$ and $C_{48}(s,9,15)$; and 
 	
    \item [\rm (c)]  $C_{48}(5,s,19)$ and $C_{48}(s,7,17)$. 
 \end{enumerate} }
\end{prm}
\noindent
{\bf Solution.}\quad Given pairs of circulant graphs are isomorphic of Type-1 by the following. 
\begin{enumerate}
      \item [\rm (a1)] $C_{48}(11(1,4,23)) = C_{48}(4,11,13)$. $\Rightarrow$ $C_{48}(1,4,23)$ and $C_{48}(4,11,13)$ are Type-1 isomorphic;
 	
      \item [\rm (a2)] $C_{48}(11(1,8,23)) = C_{48}(8,11,13)$. $\Rightarrow$ $C_{48}(1,8,23)$ and $C_{48}(8,11,13)$ are Type-1 isomorphic;
 
     \item [\rm (a3)]	 $C_{48}(11(1,12,23)) = C_{48}(11,12,13)$. $\Rightarrow$ $C_{48}(1,12,23)$ and $C_{48}(11,12,13)$ are  Type-1 isomorphic;
 	
     \item [\rm (a4)]	 $C_{48}(11(1,16,23)) = C_{48}(11,13,16)$. $\Rightarrow$ $C_{48}(1,16,23)$ and $C_{48}(11,13,16)$ are  Type-1 isomorphic;
     
     \item [\rm (a5)]	 $C_{48}(11(1,20,23)) = C_{48}(11,13,20)$. $\Rightarrow$ $C_{48}(1,20,23)$ and $C_{48}(11,13,20)$ are  Type-1 isomorphic;
     
      \item [\rm (a6)] $C_{48}(11(1,23,24)) = C_{48}(11,13,24)$. $\Rightarrow$ $C_{48}(1,23,24)$ and $C_{48}(11,13,24)$ are  Type-1 isomorphic;

     \item [\rm (b1)] $C_{48}(11(3,4,21)) = C_{48}(4,9,15)$. $\Rightarrow$ $C_{48}(3,4,21)$ and $C_{48}(4,9,15)$ are Type-1 isomorphic;

     \item [\rm (b2)] $C_{48}(11(3,8,21)) = C_{48}(8,9,15)$. $\Rightarrow$ $C_{48}(3,8,21)$ and $C_{48}(8,9,15)$ are Type-1 isomorphic;

      \item [\rm (b3)] $C_{48}(11(3,12,21)) = C_{48}(9,12,15)$. $\Rightarrow$ $C_{48}(3,12,21)$ and $C_{48}(9,12,15)$ are Type-1 isomorphic;
    
      \item [\rm (b4)] $C_{48}(11(3,16,21)) = C_{48}(9,15,16)$. $\Rightarrow$ $C_{48}(3,16,21)$ and $C_{48}(9,15,16)$ are Type-1 isomorphic;

      \item [\rm (b5)] $C_{48}(11(3,20,21)) = C_{48}(9,15,20)$. $\Rightarrow$ $C_{48}(3,20,21)$ and $C_{48}(9,15,20)$ are Type-1 isomorphic;

      \item [\rm (b6)] $C_{48}(11(3,21,24)) = C_{48}(9,15,24)$. $\Rightarrow$ $C_{48}(3,21,24)$ and $C_{48}(9,15,24)$ are Type-1 isomorphic;

    \item [\rm (c1)] $C_{48}(11(4,5,19)) = C_{48}(4,7,17)$. $\Rightarrow$ $C_{48}(4,5,19)$ and $C_{48}(4,7,17)$ are Type-1 isomorphic;

    \item [\rm (c2)] $C_{48}(11(5,8,19)) = C_{48}(7,8,17)$. $\Rightarrow$ $C_{48}(5,8,19)$ and $C_{48}(7,8,17)$ are Type-1 isomorphic;

    \item [\rm (c3)] $C_{48}(11(5,12,19)) = C_{48}(7,12,17)$. $\Rightarrow$ $C_{48}(5,12,19)$ and $C_{48}(7,12,17)$ are Type-1 isomorphic;

    \item [\rm (c4)] $C_{48}(11(5,16,19)) = C_{48}(7,16,17)$. $\Rightarrow$ $C_{48}(5,16,19)$ and $C_{48}(7,16,17)$ are Type-1 isomorphic;

   \item [\rm (c5)] $C_{48}(11(5,19,20)) = C_{48}(7,17,20)$. $\Rightarrow$ $C_{48}(5,19,20)$ and $C_{48}(7,17,20)$ are Type-1 isomorphic; 

   \item [\rm (c6)] $C_{48}(11(5,19,24)) = C_{48}(7,17,24)$. $\Rightarrow$ $C_{48}(5,19,24)$ and $C_{48}(7,17,24)$ are Type-1 isomorphic.     
   
   \hfill $\Box$
\end{enumerate}

\begin{prm} \label{p2} {\rm The following statements are true.  
\begin{enumerate}
\item [\rm (d)]	$C_{48}(1,3,23)$ $\cong$ $C_{48}(1,21,23)$ $\cong$ $C_{48}(3,7,17)$ and are  Type-1 isomorphic;

\item [\rm (e)]	$C_{48}(1,9,23)$ $\cong$ $C_{48}(1,15,23)$ $\cong$ $C_{48}(7,9,17)$ and are  Type-1 isomorphic;
			
\item [\rm (f)]   $C_{48}(3,5,19)$ $\cong$ $C_{48}(5,19,21)$ $\cong$ $C_{48}(3,11,13)$ and are  Type-1 isomorphic;
			
\item [\rm (g)]	$C_{48}(5,9,19)$ $\cong$ $C_{48}(5,15,19)$ $\cong$ $C_{48}(9,11,13)$ and are  Type-1 isomorphic.	
\end{enumerate} }
\end{prm}
\noindent
{\bf Solution.}\quad Given pairs of circulant graphs are Type-1 isomorphic by the following. 
\begin{enumerate}
\item [\rm (d)] $C_{48}(17(1,3,23))$ = $C_{48}(3,7,17)$ and $C_{48}(23(1,3,23))$ = $C_{48}(1,21,23)$. 

$\Rightarrow$ $C_{48}(1,3,23)$, $C_{48}(1,21,23)$ and $C_{48}(3,7,17)$  are Type-1 isomorphic;
	
\item [\rm (e)] $C_{48}(17(1,9,23))$ = $C_{48}(7,9,17)$ and $C_{48}(23(1,9,23))$ = $C_{48}(1,15,23)$. 

$\Rightarrow$  $C_{48}(1,9,23)$, $C_{48}(1,15,23)$ and $C_{48}(7,9,17)$  are Type-1 isomorphic;
	
\item [\rm (f)] $C_{48}(23(3,5,19))$ = $C_{48}(5,19,21)$ and $C_{48}(23(5,19,21))$ = $C_{48}(3,5,19)$. 

$\Rightarrow$  $C_{48}(3,5,19)$, $C_{48}(5,19,21)$ and $C_{48}(3,11,13)$ are  Type-1 isomorphic;
	
\item [\rm (g)] $C_{48}(17(5,9,19))$ = $C_{48}(9,11,13)$ and $C_{48}(23(5,9,19))$ = $C_{48}(5,15,19)$. 

$\Rightarrow$  $C_{48}(5,9,19)$,  $C_{48}(5,15,19)$ and $C_{48}(9,11,13)$ are  Type-1 isomorphic. 	\hfill $\Box$
\end{enumerate}

In  the next two problems, we present Type-2 isomorphic circulant graphs of $C_{48}(r_1,r_2,r_3)$ w.r.t. $m$ = 2. Here, $n$ = 48 = $2\times 3\times 2^3$ and so $m$ = 2 is the only possible value of $m$ such that $C_{n}(R)$ = $C_{48}(R)$ is having Type-2 isomorphic circulant graph w.r.t. $m$. 
  
 \begin{prm} \label{p3} {\rm Show that the following pairs of circulant graphs are Type-2 isomorphic w.r.t. $m$ = 2  for $s$ = 2, 10, 14, 22.  
 \begin{enumerate}
 	\item [\rm (h)]  $C_{48}(1,s,23)$ and $C_{48}(s,11,13)$;  
 	\item [\rm (i)]  $C_{48}(3,s,21)$ and $C_{48}(s,9,15)$; and 
 	\item [\rm (j)]  $C_{48}(5,s,19)$ and $C_{48}(s,7,17)$ where $\gcd(48, s)$ = 2.  
 \end{enumerate} }	
\end{prm}
 \noindent
 {\bf Solution.}\quad Let $R_1$ = $\{1,s,23\}$, $R_2$ = $\{s,11,13\}$, $S_1$ = $\{3,s,21\}$, $S_2$ = $\{s,9,15\}$, $T_1$ = $\{5,s,19\}$, $T_2$ = $\{s,7,17\}$ and $s$ = 2, 10, 14, 22. Here, $n$ = 48 = $6\times 2^3$ and $r$ = $s$ so that $s\in R_i,S_i,T_i$ and $\gcd(n, s)$ = $\gcd(48, 2)$ = 2 = $\gcd(48,10 )$ = $\gcd(48, 14)$ = $\gcd(48, 22)$ = $\gcd(n, n-s)$ = $\gcd(48, 48-s)$, $i$ = 1,2. This implies, $m$ = 2 is only possible value of $m$ such that given pairs of circulant graphs may be Type-2 isomorphic. Using the definition of $\theta_{n,m,t}$, we get, $\theta_{48,2,6}(s)$ = $s$ and $\theta_{48,2,6}(48-s)$ = $48-s$. Now, consider different cases.
  \begin{enumerate}
 	\item [\rm (h)]  For  $s$ = 2, 10, 14, 22, using the definition of $\theta_{n,m,t}$, we get,  
 \\
 $\theta_{48,2,6}(C_{48}(1,s,23))$ = $\theta_{48,2,6}(C_{48}(1,s,23, 25,48-s,47))$ = $C_{48}(\theta_{48,2,6}(1,s,23, 25,48-s,47))$ 
 
 \hspace{2.7cm} = $C_{48}(13,s,35, 37,48-s,11)$ = $C_{48}(s,11,13)$. 
 
 $\Rightarrow$ $C_{48}(1,s,23)$ $\cong$ $C_{48}(s,11,13)$ for $s$ = 2, 10, 14, 22. Also, for $s$ = 2, 10, 14, 22, we get,
\\ 
$Ad_{48}(C_{48}(1,s,23)$ = $Ad_{48}(C_{48}(1,s,23, 25,48-s,47))$ 
 
 \hfill = $\{C_{48}(x(1,s,23, 25,48-s,47)): x = 1,5,7,11,13,17,19,23,25,29,31,35,37,41,43,47\}$. 
 
 From the above relation, we find $Ad_{48}(C_{48}(1,s,23))$ for $s$ = 2, 10, 14, 22 and will show that $C_{48}(11,s,13) \notin Ad_{48}(C_{48}(1,s,23))$ for $s$ = 2, 10, 14, 22.
 
\item [\rm (h1)]  $Ad_{48}(C_{48}(1,2,23))$ = $\{C_{48}(1,2,23, 25,46,47)$, $C_{48}(5,10,19, 29,38,43)$,
 
 \hfill  $C_{48}(7,14,17, 31,34,41)$, $C_{48}(11,13,22, 26,35,37)\}$  

 \hspace{1.6cm} = $\{C_{48}(1,2,23)$, $C_{48}(5,10,19)\}$, $C_{48}(7,14,17)$, $C_{48}(11,13,22)\}$. 

 $\Rightarrow$ $C_{48}(2,11,13)\notin Ad_{48}(C_{48}(1,2,23))$.
 
 $\Rightarrow$ $C_{48}(1,2,23)$ and $C_{48}(2,11,13)$ are Type-2 isomorphic w.r.t. $m$ = 2.
 
\item [\rm (h2)]  $Ad_{48}(C_{48}(1,10,23))$ = $\{C_{48}(1,10,23, 25,38,47)$, $C_{48}(2,5,19, 29,43,46)$,  

\hfill $C_{48}(7,17,22, 26,31,41)$, $C_{48}(11,13,14, 34,35,37)\}$ 

\hspace{1.6cm} = $\{C_{48}(1,10,23)$, $C_{48}(2,5,19)$, $C_{48}(7,17,22)$, $C_{48}(11,13,14)\}$. 

$\Rightarrow$ $C_{48}(10,11,13)\notin Ad_{48}(C_{48}(1,10,23))$.

$\Rightarrow$ $C_{48}(1,10,23)$ and $C_{48}(10,11,13)$ are Type-2 isomorphic w.r.t. $m$ = 2.

\item [\rm (h3)]  $Ad_{48}(C_{48}(1,14,23))$ = $\{C_{48}(1,14,23, 25,34,47)$, $C_{48}(5,19,22, 26,29,43)$,  

\hfill $C_{48}(2,7,17, 31,41,46)$, $C_{48}(10,11,13, 35,37,38)\}$ 

\hspace{1.6cm} = $\{C_{48}(1,14,23)$, $C_{48}(5,19,22)$, $C_{48}(2,7,17)$, $C_{48}(10,11,13)\}$. 

$\Rightarrow$ $C_{48}(11,13,14)\notin Ad_{48}(C_{48}(1,14,23))$.

$\Rightarrow$ $C_{48}(1,14,23)$ and $C_{48}(11,13,14)$ are Type-2 isomorphic w.r.t. $m$ = 2.

\item [\rm (h4)]  $Ad_{48}(C_{48}(1,22,23))$ = $\{C_{48}(1,22,23, 25,26,47)$, $C_{48}(5,14,19, 29,34,43)$,  

\hfill $C_{48}(7,10,17, 31,38,41)$, $C_{48}(2,11,13, 35,37,46)\}$ 

\hspace{1.6cm} = $\{C_{48}(1,22,23)$, $C_{48}(5,14,19)$, $C_{48}(7,10,17)$, $C_{48}(2,11,13)\}$. 

$\Rightarrow$ $C_{48}(11,13,22)\notin Ad_{48}(C_{48}(1,22,23))$.

$\Rightarrow$ $C_{48}(1,22,23)$ and $C_{48}(11,13,22)$ are Type-2 isomorphic w.r.t. $m$ = 2.

\item [\rm (i)]  For  $s$ = 2, 10, 14, 22, using the definition of $\theta_{n,m,t}$, we get, 
\\ 	
$\theta_{48,2,6}(C_{48}(3,s,21))$ = $\theta_{48,2,6}(C_{48}(3,s,21, 27,48-s,45))$ 

\hspace{2.7cm} = $C_{48}(\theta_{48,2,6}(3,s,21, 27,48-s,45))$ 

\hspace{2.7cm} = $C_{48}(15,s,33, 39,48-s,9)$ = $C_{48}(s,9,15)$. 

$\Rightarrow$ $C_{48}(3,s,21)$ $\cong$ $C_{48}(s,9,15)$ for $s$ = 2, 10, 14, 22. Also, for $s$ = 2, 10, 14, 22, we get,
\\ 
$Ad_{48}(C_{48}(3,s,21))$ = $Ad_{48}(C_{48}(3,s,21, 27,48-s,45))$ 

\hfill = $\{C_{48}(x(3,s,21, 27,48-s,45)): x = 1,5,7,11,13,17,19,23,25,29,31,35,37,41,43,47\}$. 

From the above relation, we find $Ad_{48}(C_{48}(3,s,21))$ for $s$ = 2, 10, 14, 22 and show that $C_{48}(s,9,15) \notin Ad_{48}(C_{48}(3,s,21))$.

\item [\rm (i1)]  $Ad_{48}(C_{48}(2,3,21))$ = $\{C_{48}(2,3,21, 27,45,46)$, $C_{48}(9,10,15, 33,38,39)$,

\hfill  $C_{48}(3,14,21, 27,34,35)$,  $C_{48}(9,15,22, 27,33,39)\}$ 

\hspace{1.6cm} = $\{C_{48}(2,3,21)$, $C_{48}(9,10,15)$, $C_{48}(3,14,21)$, $C_{48}(9,15,22)\}$. 

$\Rightarrow$ $C_{48}(3,9,15)\notin Ad_{48}(C_{48}(2,3,21))$.

$\Rightarrow$ $C_{48}(2,3,21)$ and $C_{48}(3,9,15)$ are Type-2 isomorphic w.r.t. $m$ = 2.

\item [\rm (i2)]  $Ad_{48}(C_{48}(3,10,21))$ = $\{C_{48}(3,10,21, 27,38,45)$, $C_{48}(2,9,15, 33,39,46)$,  

\hfill $C_{48}(3,21,22, 26,27,45)$, $C_{48}(9,14,15, 33,34,39)\}$ 

\hspace{1.6cm} = $\{C_{48}(3,10,21)$, $C_{48}(2,9,15)$, $C_{48}(3,21,22)$, $C_{48}(9,14,15)\}$. 

$\Rightarrow$ $C_{48}(3,9,15)\notin Ad_{48}(C_{48}(3,10,21))$.

$\Rightarrow$ $C_{48}(3,10,21)$ and $C_{48}(3,9,15)$ are Type-2 isomorphic w.r.t. $m$ = 2.

\item [\rm (i3)]  $Ad_{48}(C_{48}(3,14,21))$ = $\{C_{48}(3,14,21, 27,34,45)$, $C_{48}(9,15,22, 26,33,39)$,  

\hfill $C_{48}(2,3,21, 27,45,46)$, $C_{48}(9,10,15, 33,38,39)\}$ 

\hspace{1.6cm} = $\{C_{48}(3,14,21)$, $C_{48}(9,15,22)$, $C_{48}(2,3,21)$, $C_{48}(9,10,15)\}$. 

$\Rightarrow$ $C_{48}(9,14,15)\notin Ad_{48}(C_{48}(3,14,21))$.

$\Rightarrow$ $C_{48}(3,14,21)$ and $C_{48}(9,14,15)$ are Type-2 isomorphic w.r.t. $m$ = 2.

\item [\rm (i4)]  $Ad_{48}(C_{48}(3,21,22))$ = $\{C_{48}(3,21,22, 26,27,45)$, $C_{48}(9,14,15, 33,34,39)$,  

\hfill $C_{48}(3,10,21, 27,38,45)$, $C_{48}(2,9,15, 33,39,46)\}$ 

\hspace{1.6cm} = $\{C_{48}(3,21,22)$, $C_{48}(9,14,15)$, $C_{48}(3,10,21)$, $C_{48}(2,9,15)\}$. 

$\Rightarrow$ $C_{48}(9,15,22)\notin Ad_{48}(C_{48}(3,21,22))$.

$\Rightarrow$ $C_{48}(3,21,22)$ and $C_{48}(9,15,22)$ are Type-2 isomorphic w.r.t. $m$ = 2.

\item [\rm (j)]  For  $s$ = 2, 10, 14, 22, using the definition of $\theta_{n,m,t}$, we get, 
\\ 	
$\theta_{48,2,6}(C_{48}(5,s,19))$ = $\theta_{48,2,6}(C_{48}(5,s,19, 29,48-s,43))$ 

\hspace{2.7cm} = $C_{48}(\theta_{48,2,6}(5,s,19, 29,48-s,43))$ 

\hspace{2.7cm} = $C_{48}(17,s,31, 41,48-s,7)$ = $C_{48}(s,7,17)$. 

$\Rightarrow$ $C_{48}(5,s,19)$ $\cong$ $C_{48}(s,7,17)$ for $s$ = 2, 10, 14, 22. Also, for $s$ = 2, 10, 14, 22, we get,
\\ 
$Ad_{48}(C_{48}(5,s,19))$ = $Ad_{48}(C_{48}(5,s,19, 29,48-s,43))$ 

\hfill = $\{C_{48}(x(5,s,19, 29,48-s,43)): x = 1,5,7,11,13,17,19,23,25,29,31,35,37,41,43,47\}$. 

From the above relation, we find $Ad_{48}(C_{48}(5,s,19))$ for $s$ = 2, 10, 14, 22 and show that $C_{48}(s,7,17) \notin Ad_{48}(C_{48}(5,s,19))$.

\item [\rm (j1)]  $Ad_{48}(C_{48}(2,5,19))$ = $\{C_{48}(2,5,19, 29,43,46)$, $C_{48}(1,10,23, 25,38,47)$,

\hfill  $C_{48}(11,13,14, 34,35,37)$,  $C_{48}(7,17,22, 27,31,41)\}$ 

\hspace{1.6cm} = $\{C_{48}(2,5,19)$, $C_{48}(1,10,23)$, $C_{48}(11,13,14)$, $C_{48}(7,17,22)\}$. 

$\Rightarrow$ $C_{48}(2,7,17)\notin Ad_{48}(C_{48}(2,5,19))$.

$\Rightarrow$ $C_{48}(2,5,19)$ and $C_{48}(2,7,17)$ are Type-2 isomorphic w.r.t. $m$ = 2.

\item [\rm (j2)]  $Ad_{48}(C_{48}(5,10,19))$ = $\{C_{48}(5,10,19, 29,38,43)$, $C_{48}(1,2,23, 25,46,47)$,  

\hfill $C_{48}(11,13,22, 26,35,37)$, $C_{48}(7,14,17, 31,34,41)\}$ 

\hspace{1.6cm} = $\{C_{48}(5,10,19)$, $C_{48}(1,2,23)$, $C_{48}(11,13,22)$, $C_{48}(7,14,17)\}$. 

$\Rightarrow$ $C_{48}(7,10,17)\notin Ad_{48}(C_{48}(5,10,19))$.

$\Rightarrow$ $C_{48}(5,10,19)$ and $C_{48}(7,10,17)$ are Type-2 isomorphic w.r.t. $m$ = 2.

\item [\rm (j3)]  $Ad_{48}(C_{48}(5,14,19))$ = $\{C_{48}(5,14,19, 29,34,43)$, $C_{48}(1,22,23, 25,26,47)$,  

\hfill $C_{48}(2,11,13, 35,37,46)$, $C_{48}(7,10,17, 31,38,41)\}$ 

\hspace{1.6cm} = $\{C_{48}(5,14,19)$, $C_{48}(1,22,23)$, $C_{48}(2,11,13)$, $C_{48}(7,10,17)\}$. 

$\Rightarrow$ $C_{48}(7,14,17)\notin Ad_{48}(C_{48}(5,14,19))$.

$\Rightarrow$ $C_{48}(5,14,19)$ and $C_{48}(7,14,17)$ are Type-2 isomorphic w.r.t. $m$ = 2.

\item [\rm (j4)]  $Ad_{48}(C_{48}(5,19,22))$ = $\{C_{48}(5,19,22, 26,29,43)$, $C_{48}(1,14,23, 25,34,47)$,  

\hfill $C_{48}(10,11,13, 35,37,38)$, $C_{48}(2,7,17, 31,41,46)\}$ 

\hspace{1.6cm} = $\{C_{48}(5,19,22)$, $C_{48}(1,14,23)$, $C_{48}(10,11,13)$, $C_{48}(2,7,17)\}$. 

$\Rightarrow$ $C_{48}(7,17,22)\notin Ad_{48}(C_{48}(5,19,22))$.

$\Rightarrow$ $C_{48}(5,19,22)$ and $C_{48}(7,17,22)$ are Type-2 isomorphic w.r.t. $m$ = 2.  \hfill $\Box$
   \end{enumerate} 
 
 \begin{prm} \label{p4} {\rm Show that the following pairs of circulant graphs are Type-2 isomorphic w.r.t. $m$ = 2 for $s$ = 6, 18.  
\begin{enumerate}
    \item [\rm (x)]  $C_{48}(1,s,23)$ and $C_{48}(s,11,13)$;  
	\item [\rm (y)]  $C_{48}(3,s,21)$ and $C_{48}(s,9,15)$; and 
	\item [\rm (z)]  $C_{48}(5,s,19)$ and $C_{48}(s,7,17)$ where gcd(48, s) = 6.  
	\end{enumerate} }	
\end{prm}
\noindent
{\bf Solution.}\quad Let $R_1$ = $\{1,s,23\}$, $R_2$ = $\{s,11,13\}$, $S_1$ = $\{3,s,21\}$, $S_2$ = $\{s,9,15\}$, $T_1$ = $\{5,s,19\}$, $T_2$ = $\{s,7,17\}$ and $s$ = 6, 18. Here, $n$ = 48 = $6\times 2^3$ and $r$ = $s$ so that $s\in R_i,S_i,T_i$ and $\gcd(n, s)$ = $\gcd(48, 6)$ = 6 = $\gcd(48, 18)$, $i$ = 1,2. This implies, $m$ = 2 is the only possible value of $m$ such that given pairs of circulant graphs may be Type-2 isomorphic w.r.t. $m$. Using the definition of $\theta_{n,m,t}$, we get, $\theta_{48,2,6}(s)$ = $s$ and $\theta_{48,2,6}(48-s)$ = $48-s$ when $s$ = 6 or 18. Now, consider different cases.
\begin{enumerate}
	\item [\rm (x)]  For $s$ = 6,18, using the definition of $\theta_{n,m,t}$, we get,
	\\
	$\theta_{48,2,6}(C_{48}(1,s,23))$ = $\theta_{48,2,6}(C_{48}(1,s,23, 25,48-s,47))$  = $C_{48}(\theta_{48,2,6}(1,s,23, 25,48-s,47))$ 
	
	\hspace{2.7cm} = $C_{48}(13,s,35, 37,48-s,11)$ = $C_{48}(s,11,13)$.
	
	$\Rightarrow$ $C_{48}(1,s,23)$ $\cong$ $C_{48}(s,11,13)$ for $s$ = 6,18. Also, for $s$ = 6,18, we get, 
	\\ 
	$Ad_{48}(C_{48}(1,s,23))$ = $Ad_{48}(C_{48}(1,s,23, 25,48-s,47))$ 
	
	\hfill = $\{C_{48}(x(1,s,23, 25,48-s,47)): x = 1,5,7,11,13,17,19,23,25,29,31,35,37,41,43,47\}$.
	
	From the above relation, we find $Ad_{48}(C_{48}(1,s,23))$ for $s$ = 6, 18, seperately and will show that $C_{48}(11,s,13)\notin Ad_{48}(C_{48}(1,s,23))$ for $s$ = 6, 18.
	
 \item [\rm (x1)] $Ad_{48}(C_{48}(1,6,23))$ = $\{C_{48}(1,6,23)$, $C_{48}(5,18,19)$, $C_{48}(6,7,17)$, $C_{48}(11,13,18)\}$. 
	
	$\Rightarrow$ $C_{48}(6,11,13)\notin Ad_{48}(C_{48}(1,6,23))$.
	
	$\Rightarrow$ $C_{48}(1,6,23)$ and $C_{48}(6,11,13)$ are Type-2 isomorphic w.r.t. $m$ = 2.
	
    \item [\rm (x2)] $Ad_{48}(C_{48}(1,18,23))$ = $\{C_{48}(1,18,23)$, $C_{48}(5,6,19)$, $C_{48}(7,17,18)$, $C_{48}(6,11,13)\}$. 

    $\Rightarrow$ $C_{48}(11,13,18)\notin Ad_{48}(C_{48}(1,18,23))$.

    $\Rightarrow$ $C_{48}(1,18,23)$ and $C_{48}(11,13,18)$ are Type-2 isomorphic w.r.t. $m$ = 2.

\item [\rm (y)]  When $s$ = 6, $C_{48}(3,s,21)$ = $C_{48}(3,6,21)$ = $3.C_{16}(1,2,7)$ and $C_{48}(s,9,15)$ = $C_{48}(6,9,15)$ = $3.C_{16}(2,3,5)$. This implies, $C_{48}(3,6,21)$ and $C_{48}(6,9,15)$ are Type-2 isomorphic w.r.t. $m$ = 2 since $C_{16}(1,2,7)$ and $C_{16}(2,3,5)$ are Type-2 isomorphic w.r.t. $m$ = 2. 

Similarly, when $s$ = 18, $C_{48}(3,s,21)$ = $C_{48}(3,18,21)$ = $3.C_{16}(1,6,7)$ and $C_{48}(s,9,15)$ = $C_{48}(9,15,18)$ = $3.C_{16}(3,5,6)$. This implies that $C_{48}(3,18,21)$ and $C_{48}(9,15,18)$ are Type-2 isomorphic w.r.t. $m$ = 2 since $C_{16}(1,6,7)$ and $C_{16}(2,5,6)$ are Type-2 isomorphic w.r.t. $m$ = 2. 

\item [\rm (z)]  For $s$ = 6,18, using the definition of $\theta_{n,m,t}$, we get,
\\
$\theta_{48,2,6}(C_{48}(5,s,19))$ = $\theta_{48,2,6}(C_{48}(5,s,19, 29,48-s,43))$  = $C_{48}(\theta_{48,2,6}(5,s,19, 29,48-s,43))$ 

\hspace{2.7cm} = $C_{48}(17,s,31, 41,48-s,7)$ = $C_{48}(s,7,17)$.

$\Rightarrow$ $C_{48}(5,s,19)$ $\cong$ $C_{48}(s,7,17)$ for $s$ = 6,18. Also, for $s$ = 6,18, we get, 
\\ 
$Ad_{48}(C_{48}(5,s,19))$ = $Ad_{48}(C_{48}(5,s,19, 29,48-s,43))$ 

\hfill = $\{C_{48}(x(5,s,19, 29,48-s,43)): x = 1,5,7,11,13,17,19,23,25,29,31,35,37,41,43,47\}$.

From the above relation, we find $Ad_{48}(C_{48}(5,s,19))$ for $s$ = 6, 18, seperately and show that $C_{48}(s,7,17)\notin Ad_{48}(C_{48}(5,s,19))$.

\item [\rm (z1)] $Ad_{48}(C_{48}(5,6,19))$ = $\{C_{48}(5,6,19)$, $C_{48}(1,18,23)$ $C_{48}(6,11,13)$, $C_{48}(7,17,18)\}$. 

$\Rightarrow$ $C_{48}(6,7,17)\notin Ad_{48}(C_{48}(5,6,19))$.

$\Rightarrow$ $C_{48}(5,6,19)$ and $C_{48}(6,7,17)$ are Type-2 isomorphic w.r.t. $m$ = 6 = $\gcd(48, 6)$ = $r$. 

\item [\rm (z2)] $Ad_{48}(C_{48}(5,18,19))$ = $\{C_{48}(5,18,19)$, $C_{48}(1,6,23)$, $C_{48}(11,13,18)$, $C_{48}(6,7,17)\}$. 

$\Rightarrow$ $C_{48}(7,17,18)\notin Ad_{48}(C_{48}(5,18,19))$.

$\Rightarrow$ $C_{48}(5,18,19)$ and $C_{48}(7,17,18)$ are Type-2 isomorphic w.r.t. $m$ = 6, $r$ = 18.    \hfill $\Box$
\end{enumerate} 

\vspace{.2cm}
In the above problems, values of $s$ are obtained from the following.
\begin{enumerate}
 	\item [\rm (i)]  $\gcd(48, s)$ = 2 for $s$ = 2,10,14,22, 26,34,38,46; and
 	
 	\item [\rm (ii)]  $\gcd(48, s)$ = 6 for $s$ = 6,18, 30,42. 
 \end{enumerate}

\medskip
\noindent
{\bf Type-1 sets of $C_{48}(r_1,r_2,r_3)$ related to problems \ref{p1} to \ref{p4}}

\medskip
\noindent
We present here the following Type-1 sets which are related to that which occured in problems \ref{p1} to \ref{p4}.
\begin{enumerate}
\item [\rm (a1)]  $Ad_{48}(C_{48}(1,4,23))$ = $\{C_{48}(1,4,23)$, $C_{48}(5,19,20)$, $C_{48}(7,17,20)$, $C_{48}(4,11,13)$

 = $Ad_{48}(C_{48}(5,19,20))$ = $Ad_{48}(C_{48}(7,17,20))$ = $Ad_{48}(C_{48}(4,11,13))$.

\item [\rm (a2)] $Ad_{48}(C_{48}(1,8,23))$ = $\{C_{48}(1,8,23)$, $C_{48}(5,8,19)$, $C_{48}(7,8,17)$, $C_{48}(8,11,13)$

 = $Ad_{48}(C_{48}(5,8,19))$ = $Ad_{48}(C_{48}(7,8,17))$ = $Ad_{48}(C_{48}(8,11,13))$.

\item [\rm (a3)] $Ad_{48}(C_{48}(1,12,23))$ = $\{C_{48}(1,12,23)$, $C_{48}(5,12,19)$, $C_{48}(7,12,17)$, $C_{48}(11,12,13)$

 = $Ad_{48}(C_{48}(5,12,19))$ = $Ad_{48}(C_{48}(7,12,17))$ = $Ad_{48}(C_{48}(11,12,13))$.

\item [\rm (a4)] $Ad_{48}(C_{48}(1,16,23))$ = $\{C_{48}(1,16,23)$, $C_{48}(5,16,19)$, $C_{48}(7,16,17)$, $C_{48}(11,13,16)$

 = $Ad_{48}(C_{48}(5,16,19))$ = $Ad_{48}(C_{48}(7,16,17))$ = $Ad_{48}(C_{48}(11,13,16))$.

\item [\rm (a5)] $Ad_{48}(C_{48}(1,20,23))$ = $\{C_{48}(1,20,23)$, $C_{48}(4,5,19)$, $C_{48}(4,7,17)$, $C_{48}(11,13,20)$

 = $Ad_{48}(C_{48}(4,5,19))$ = $Ad_{48}(C_{48}(4,7,17))$ = $Ad_{48}(C_{48}(11,13,20))$.

\item [\rm (a6)] $Ad_{48}(C_{48}(1,23,24))$ = $\{C_{48}(1,23,24)$, $C_{48}(5,19,24)$, $C_{48}(7,17,24)$, $C_{48}(11,13,24)$

 = $Ad_{48}(C_{48}(5,19,24))$ = $Ad_{48}(C_{48}(7,17,24))$ = $Ad_{48}(C_{48}(11,13,24))$.

\item [\rm (b1)] $Ad_{48}(C_{48}(3,4,21))$ = $\{C_{48}(3,4,21)$, $C_{48}(9,15,20)$, $C_{48}(3,20,21)$, $C_{48}(4,9,15)$

 = $Ad_{48}(C_{48}(9,15,20))$ = $Ad_{48}(C_{48}(3,20,21))$ = $Ad_{48}(C_{48}(4,9,15))$.

\item [\rm (b2)] $Ad_{48}(C_{48}(3,8,21))$ = $\{C_{48}(3,8,21)$, $C_{48}(8,9,15)$, $C_{48}(3,8,21)$, $C_{48}(8,9,15)$

 = $Ad_{48}(C_{48}(8,9,15))$ = $Ad_{48}(C_{48}(3,8,21))$ = $Ad_{48}(C_{48}(8,9,15))$.

\item [\rm (b3)] $Ad_{48}(C_{48}(3,12,21))$ = $\{C_{48}(3,12,21)$, $C_{48}(9,12,15)$, $C_{48}(3,12,21)$, $C_{48}(9,12,15)$

 = $Ad_{48}(C_{48}(9,12,15))$ = $Ad_{48}(C_{48}(3,12,21))$ = $Ad_{48}(C_{48}(9,12,15))$.

\item [\rm (b4)] $Ad_{48}(C_{48}(3,16,21))$ = $\{C_{48}(3,16,21)$, $C_{48}(9,15,16)$, $C_{48}(3,16,21)$, $C_{48}(9,15,16)$

 = $Ad_{48}(C_{48}(9,15,16))$ = $Ad_{48}(C_{48}(3,16,21))$ = $Ad_{48}(C_{48}(9,15,16))$.

\item [\rm (b5)] $Ad_{48}(C_{48}(3,20,21))$ = $\{C_{48}(3,20,21)$, $C_{48}(4,9,15)$, $C_{48}(3,4,21)$, $C_{48}(9,15,20)$

 = $Ad_{48}(C_{48}(4,9,15))$ = $Ad_{48}(C_{48}(3,4,21))$ = $Ad_{48}(C_{48}(9,15,20))$.

\item [\rm (b6)] $Ad_{48}(C_{48}(3,21,24))$ = $\{C_{48}(3,21,24)$, $C_{48}(9,15,24)$, $C_{48}(3,21,24)$, $C_{48}(9,15,24)$

 = $Ad_{48}(C_{48}(9,15,24))$ = $Ad_{48}(C_{48}(3,21,24))$ = $Ad_{48}(C_{48}(9,15,24))$.

\item [\rm (c1)] $Ad_{48}(C_{48}(4,5,19))$ = $\{C_{48}(4,5,19)$, $C_{48}(1,20,23)$, $C_{48}(11,13,20)$, $C_{48}(4,7,17)$

 = $Ad_{48}(C_{48}(1,20,23))$ = $Ad_{48}(C_{48}(11,13,20))$ = $Ad_{48}(C_{48}(4,7,17))$.

\item [\rm (c2)] $Ad_{48}(C_{48}(5,8,19))$ = $\{C_{48}(5,8,19)$, $C_{48}(1,8,23)$, $C_{48}(8,11,13)$, $C_{48}(7,8,17)$

 = $Ad_{48}(C_{48}(1,8,23))$ = $Ad_{48}(C_{48}(8,11,13))$ = $Ad_{48}(C_{48}(7,8,17))$.

\item [\rm (c3)] $Ad_{48}(C_{48}(5,12,19))$ = $\{C_{48}(5,12,19)$, $C_{48}(1,12,23)$, $C_{48}(11,12,13)$, $C_{48}(7,12,17)$

 = $Ad_{48}(C_{48}(1,12,23))$ = $Ad_{48}(C_{48}(11,12,13))$ = $Ad_{48}(C_{48}(7,12,17))$.

\item [\rm (c4)] $Ad_{48}(C_{48}(5,16,19))$ = $\{C_{48}(5,16,19)$, $C_{48}(1,16,23)$, $C_{48}(11,13,16)$, $C_{48}(7,16,17)$

 = $Ad_{48}(C_{48}(1,16,23))$ = $Ad_{48}(C_{48}(11,13,16))$ = $Ad_{48}(C_{48}(7,16,17))$.

\item [\rm (c5)] $Ad_{48}(C_{48}(5,19,20))$ = $\{C_{48}(5,19,20)$, $C_{48}(1,4,23)$, $C_{48}(4,11,13)$, $C_{48}(7,17,20)$

 = $Ad_{48}(C_{48}(1,4,23))$ = $Ad_{48}(C_{48}(4,11,13))$ = $Ad_{48}(C_{48}(7,17,20))$.

\item [\rm (c6)] $Ad_{48}(C_{48}(5,19,24))$ = $\{C_{48}(5,19,24)$, $C_{48}(1,23,24)$, $C_{48}(11,13,24)$, $C_{48}(7,17,24)$

 = $Ad_{48}(C_{48}(1,23,24))$ = $Ad_{48}(C_{48}(11,13,24))$ = $Ad_{48}(C_{48}(7,17,24))$.

\item [\rm (d)]  $Ad_{48}(C_{48}(1,3,23))$ = $\{C_{48}(1,3,23)$, $C_{48}(5,15,19)$, $C_{48}(7,17,21)$, $C_{48}(11,13,15)$, 

\hfill $C_{48}(9,11,13)$, $C_{48}(3,7,17)$, $C_{48}(5,9,19)$, $C_{48}(1,21,23)\}$

 = $Ad_{48}(C_{48}(5,15,19))$ = $Ad_{48}(C_{48}(7,17,21))$ = $Ad_{48}(C_{48}(11,13,15))$ = $Ad_{48}(C_{48}(9,11,13))$ 

= $Ad_{48}(C_{48}(3,7,17))$ = $Ad_{48}(C_{48}(5,9,19))$ = = $Ad_{48}(C_{48}(1,21,23))$.

\item [\rm (e)]  $Ad_{48}(C_{48}(1,9,23))$ = $\{C_{48}(1,9,23)$, $C_{48}(3,5,19)$, $C_{48}(7,15,17)$, $C_{48}(3,11,13)$, 

\hfill $C_{48}(11,13,21)$, $C_{48}(7,9,17)$, $C_{48}(5,19,21)$, $C_{48}(1,15,23)\}$

 = $Ad_{48}(C_{48}(3,5,19))$ = $Ad_{48}(C_{48}(7,15,17))$ = $Ad_{48}(C_{48}(3,11,13))$ = $Ad_{48}(C_{48}(11,13,21))$ 

= $Ad_{48}(C_{48}(7,9,17))$ = $Ad_{48}(C_{48}(5,19,21))$ = = $Ad_{48}(C_{48}(1,15,23))$.

\item [\rm (f)]  $Ad_{48}(C_{48}(3,5,19))$ = $\{C_{48}(3,5,19)$, $C_{48}(1,15,23)$, $C_{48}(11,13,21)$, $C_{48}(7,15,17)$, 

\hfill $C_{48}(7,9,17)$, $C_{48}(3,11,13)$, $C_{48}(1,9,23)$, $C_{48}(5,19,21)\}$

 = $Ad_{48}(C_{48}(1,15,23))$ = $Ad_{48}(C_{48}(11,13,21))$ = $Ad_{48}(C_{48}(7,15,17))$ = $Ad_{48}(C_{48}(7,9,17))$ 

= $Ad_{48}(C_{48}(3,11,13))$ = $Ad_{48}(C_{48}(1,9,23))$ = = $Ad_{48}(C_{48}(5,19,21))$.

\item [\rm (g)]  $Ad_{48}(C_{48}(5,9,19))$ = $\{C_{48}(5,9,19)$, $C_{48}(1,3,23)$, $C_{48}(11,13,15)$, $C_{48}(3,7,17)$, 

\hfill $C_{48}(7,17,21)$, $C_{48}(9,11,13)$, $C_{48}(1,21,23)$, $C_{48}(5,15,19)\}$

 = $Ad_{48}(C_{48}(1,3,23))$ = $Ad_{48}(C_{48}(11,13,15))$ = $Ad_{48}(C_{48}(3,7,17))$ = $Ad_{48}(C_{48}(7,17,21))$ 

= $Ad_{48}(C_{48}(9,11,13))$ = $Ad_{48}(C_{48}(1,21,23))$ = = $Ad_{48}(C_{48}(5,15,19))$.

\item [\rm (h1)]  $Ad_{48}(C_{48}(1,2,23))$ = $\{C_{48}(1,2,23)$, $C_{48}(5,10,19)$, $C_{48}(7,14,17)$, $C_{48}(11,13,22)\}$

\hfill = $Ad_{48}(C_{48}(5,10,19))$ = $Ad_{48}(C_{48}(7,14,17))$ = $Ad_{48}(C_{48}(11,13,22))$.

\item [\rm (h2)]  $Ad_{48}(C_{48}(1,10,23))$ = $\{C_{48}(1,10,23)$, $C_{48}(2,5,19)$, $C_{48}(7,17,22)$, $C_{48}(11,13,14)\}$

\hfill = $Ad_{48}(C_{48}(2,5,19))$ = $Ad_{48}(C_{48}(7,17,22))$ = $Ad_{48}(C_{48}(11,13,14))$.

\item [\rm (h3)]  $Ad_{48}(C_{48}(1,14,23))$ = $\{C_{48}(1,14,23)$, $C_{48}(5,19,22)$, $C_{48}(2,7,17)$, $C_{48}(10,11,13)\}$

\hfill = $Ad_{48}(C_{48}(5,19,22))$ = $Ad_{48}(C_{48}(2,7,17))$ = $Ad_{48}(C_{48}(10,11,13))$. 

\item [\rm (h4)]  $Ad_{48}(C_{48}(1,22,23))$ = $\{C_{48}(1,22,23)$, $C_{48}(5,14,19)$, $C_{48}(7,10,17)$, $C_{48}(2,11,13)\}$

\hfill = $Ad_{48}(C_{48}(5,14,19))$ = $Ad_{48}(C_{48}(7,10,17))$ = $Ad_{48}(C_{48}(2,11,13))$.

\item [\rm (i1)]  $Ad_{48}(C_{48}(2,3,21))$ = $\{C_{48}(2,3,21)$, $C_{48}(9,10,15)$, $C_{48}(3,14,21)$, $C_{48}(9,15,22)\}$

\hfill = $Ad_{48}(C_{48}(9,10,15))$ = $Ad_{48}(C_{48}(3,14,21))$ = $Ad_{48}(C_{48}(9,15,22))$. 

\item [\rm (i2)]  $Ad_{48}(C_{48}(3,10,21))$ = $\{C_{48}(3,10,21)$, $C_{48}(2,9,15)$, $C_{48}(3,21,22)$, $C_{48}(9,14,15)\}$

\hfill = $Ad_{48}(C_{48}(2,9,15))$ = $Ad_{48}(C_{48}(3,21,22))$ = $Ad_{48}(C_{48}(9,14,15))$. 

\item [\rm (i3)]  $Ad_{48}(C_{48}(3,14,21))$ = $\{C_{48}(3,14,21)$, $C_{48}(9,15,22)$, $C_{48}(2,3,21)$, $C_{48}(9,10,15)\}$

\hfill = $Ad_{48}(C_{48}(9,15,22))$ = $Ad_{48}(C_{48}(2,3,21))$ = $Ad_{48}(C_{48}(9,10,15))$.

\item [\rm (i4)]  $Ad_{48}(C_{48}(3,21,22))$ = $\{C_{48}(3,21,22)$, $C_{48}(9,14,15)$, $C_{48}(3,10,21)$, $C_{48}(2,9,15)\}$

\hfill = $Ad_{48}(C_{48}(9,14,15))$ = $Ad_{48}(C_{48}(3,10,21))$ = $Ad_{48}(C_{48}(2,9,15))$. 

\item [\rm (j1)]  $Ad_{48}(C_{48}(2,5,19))$ = $\{C_{48}(2,5,19)$, $C_{48}(1,10,23)$, $C_{48}(11,13,14)$, $C_{48}(7,17,22)\}$

\hfill = $Ad_{48}(C_{48}(1,10,23))$ = $Ad_{48}(C_{48}(11,13,14))$ = $Ad_{48}(C_{48}(7,17,22))$. 

\item [\rm (j2)]  $Ad_{48}(C_{48}(5,10,19))$ = $\{C_{48}(5,10,19)$, $C_{48}(1,2,23)$, $C_{48}(11,13,22)$, $C_{48}(7,14,17)\}$

\hfill = $Ad_{48}(C_{48}(1,2,23))$ = $Ad_{48}(C_{48}(11,13,22))$ = $Ad_{48}(C_{48}(7,14,17))$. 

\item [\rm (j3)]  $Ad_{48}(C_{48}(5,14,19))$ = $\{C_{48}(5,14,19)$, $C_{48}(1,22,23)$, $C_{48}(2,11,13)$, $C_{48}(7,10,17)\}$

\hfill = $Ad_{48}(C_{48}(1,22,23))$ = $Ad_{48}(C_{48}(2,11,13))$ = $Ad_{48}(C_{48}(7,10,17))$.

\item [\rm (j4)]  $Ad_{48}(C_{48}(5,19,22))$ = $\{C_{48}(5,19,22)$, $C_{48}(1,14,23)$, $C_{48}(10,11,13)$, $C_{48}(2,7,17)\}$

\hfill = $Ad_{48}(C_{48}(1,14,23))$ = $Ad_{48}(C_{48}(10,11,13))$ = $Ad_{48}(C_{48}(2,7,17))$.  

 \item [\rm (x1)] $Ad_{48}(C_{48}(1,6,23))$ = $\{C_{48}(1,6,23)$, $C_{48}(5,18,19)$, $C_{48}(6,7,17)$, $C_{48}(11,13,18)\}$

\hfill = $Ad_{48}(C_{48}(5,18,19))$ = $Ad_{48}(C_{48}(6,7,17))$ = $Ad_{48}(C_{48}(11,13,18))$.

\item [\rm (x2)] $Ad_{48}(C_{48}(1,18,23))$ = $\{C_{48}(1,18,23)$, $C_{48}(5,6,19)$, $C_{48}(7,17,18)$, $C_{48}(6,11,13)\}$

\hfill = $Ad_{48}(C_{48}(5,6,19))$ = $Ad_{48}(C_{48}(7,17,18))$ = $Ad_{48}(C_{48}(6,11,13))$.

 \item [\rm (y1)] $Ad_{48}(C_{48}(1,6,21))$ = $\{C_{48}(1,6,21)$, $C_{48}(5,9,18)$, $C_{48}(3,6,7)$, $C_{48}(9,11,18)$, $C_{48}(13,15,18)$, 

\hfill $C_{48}(6,17,21)$, $C_{48}(15,18,19)$, $C_{48}(3,6,23)\}$  

= $Ad_{48}(C_{48}(5,9,18))$ = $Ad_{48}(C_{48}(3,6,7))$ = $Ad_{48}(C_{48}(9,11,18))$ = $Ad_{48}(C_{48}(13,15,18))$ 

= $Ad_{48}(C_{48}(6,17,21))$ = $Ad_{48}(C_{48}(15,18,19))$ = $Ad_{48}(C_{48}(3,6,23))$.

\item [\rm (y2)] $Ad_{48}(C_{48}(1,18,21))$ = $\{C_{48}(1,18,21)$, $C_{48}(5,6,9)$, $C_{48}(3,7,18)$, $C_{48}(6,9,11)$, $C_{48}(6,13,15)$, 

\hfill $C_{48}(17,18,21)$, $C_{48}(6,15,19)$, $C_{48}(3,18,23)\}$  

= $Ad_{48}(C_{48}(5,6,9))$ = $Ad_{48}(C_{48}(3,7,18))$ = $Ad_{48}(C_{48}(6,9,11))$ = $Ad_{48}(C_{48}(6,13,15))$

 = $Ad_{48}(C_{48}(17,18,21))$ = $Ad_{48}(C_{48}(6,15,19))$ = $Ad_{48}(C_{48}(3,18,23))$.

\item [\rm (z1)] $Ad_{48}(C_{48}(5,6,19))$ = $\{C_{48}(5,6,19)$, $C_{48}(1,18,23)$ $C_{48}(6,11,13)$, $C_{48}(7,17,18)\}$

\hfill = $Ad_{48}(C_{48}(1,18,23))$ = $Ad_{48}(C_{48}(6,11,13))$ = $Ad_{48}(C_{48}(7,17,18))$.

\item [\rm (z2)] $Ad_{48}(C_{48}(5,18,19))$ = $\{C_{48}(5,18,19)$, $C_{48}(1,6,23)$, $C_{48}(11,13,18)$, $C_{48}(6,7,17)\}$

\hfill = $Ad_{48}(C_{48}(1,6,23))$ = $Ad_{48}(C_{48}(11,13,18))$ = $Ad_{48}(C_{48}(6,7,17))$.
\end{enumerate} 
 
\subsection{On isomorphic circulant graphs $C_{81}(r_1,r_2,r_3)$ of Type-2 w.r.t. $m$ = 3}

In this subsection, we study Type-2 isomorphic circulant graphs of $C_{81}(r_1,r_2,r_3)$. We obtain all triples of Type-2 isomorphic circulant graphs of the form $C_{81}(x_1,x_2,x_3)$, $C_{81}(y_1,y_2,y_3)$ and $C_{81}(z_1,z_2,z_3)$ and the total number of such triples are 27 and all are of Type-2 w.r.t. $m$ = 3. Here, $n$ = 81 = $3\times 3^3$ and so $m$ = 3 is the only possible value of $m$ such that $C_{n}(R)$ = $C_{81}(R)$ is Type-2 isomorphic w.r.t. $m$. We start with a problem retaled to Type-1 isomorphic circulant graphs of $C_{81}(r_1,r_2,r_3)$.

\begin{prm} \label{p5} {\rm For $s$ = 9,18,27, 36 and $\gcd(81, s)$ = $s$, show that the following statements are true.   
\begin{enumerate}
 \item [\rm (a)] $C_{81}(1,s,26,28)$, $C_{81}(8,s,19,35)$, $C_{81}(s,10,17,37)$ are Type-1 isomorphic; 
 \item [\rm (b)] $C_{81}(2,s,25,29)$, $C_{81}(7,s,20,34)$, $C_{81}(s,11,16,38)$ are Type-1 isomorphic.  
 \item [\rm (c)] $C_{81}(4,s,23,31)$, $C_{81}(s,13,14,40)$, $C_{81}(5,s,22,32)$ are Type-1 isomorphic.  
\end{enumerate} }	
\end{prm}
\noindent
{\bf Solution.}\quad Given triples of circulant graphs are Type-1 isomorphic by the following. 
\begin{enumerate}
 \item [\rm (a1)] $C_{81}(8(1,9,26,28))$ = $C_{81}(8,9,19,35)$ and $C_{81}(10(1,9,26,28))$ = $C_{81}(9,10,17,37)$.

$\Rightarrow$ $C_{81}(1,9,26,28)$, $C_{81}(8,9,19,35)$, $C_{81}(9,10,17,37)$ are Type-1 isomorphic.

\item [\rm (a2)] $C_{81}(8(1,18,26,28))$ = $C_{81}(8,18,19,35)$ and $C_{81}(10(1,18,26,28))$ = $C_{81}(10,17,18,37)$.

$\Rightarrow$ $C_{81}(1,18,26,28)$, $C_{81}(8,18,19,35)$, $C_{81}(10,17,18,37)$ are Type-1 isomorphic.

\item [\rm (a3)] $C_{81}(8(1,26,27,28))$ = $C_{81}(8,19,27,35)$ and $C_{81}(10(1,26,27,28))$ = $C_{81}(10,17,27,37)$.

$\Rightarrow$ $C_{81}(1,26,27,28)$, $C_{81}(8,19,27,35)$, $C_{81}(10,17,27,37)$ are Type-1 isomorphic.

\item [\rm (a4)] $C_{81}(8(1,26,28,36))$ = $C_{81}(8,19,35,36)$ and $C_{81}(10(1,26,28,36))$ = $C_{81}(10,17,36,37)$.

$\Rightarrow$ $C_{81}(1,26,28,36)$, $C_{81}(8,19,35,36)$, $C_{81}(10,17,36,37)$ are Type-1 isomorphic.

\item [\rm (b1)] $C_{81}(8(2,9,25,29))$ = $C_{81}(9,11,16,38)$ and $C_{81}(10(2,9,25,29))$ = $C_{81}(7,9,20,34)$.  
 
$\Rightarrow$ $C_{81}(2,9,25,29)$, $C_{81}(9,11,16,38)$, $C_{81}(7,9,20,34)$ are Type-1 isomorphic.

\item [\rm (b2)] $C_{81}(8(2,18,25,29))$ = $C_{81}(11,16,18,38)$ and $C_{81}(10(2,18,25,29))$ = $C_{81}(7,18,20,34)$.  
 
$\Rightarrow$ $C_{81}(2,18,25,29)$, $C_{81}(11,16,18,38)$, $C_{81}(7,18,20,34)$ are Type-1 isomorphic.

\item [\rm (b3)] $C_{81}(8(2,25,27,29))$ = $C_{81}(11,16,27,38)$ and $C_{81}(10(2,25,27,29))$ = $C_{81}(7,20,27,34)$  
 
$\Rightarrow$ $C_{81}(2,25,27,29)$, $C_{81}(11,16,27,38)$, $C_{81}(7,20,27,34)$ are Type-1 isomorphic.

\item [\rm (b4)] $C_{81}(8(2,25,29,36))$ = $C_{81}(11,16,36,38)$ and $C_{81}(10(2,25,29,36))$ = $C_{81}(7,20,34,36)$.  
 
$\Rightarrow$ $C_{81}(2,25,29,36)$, $C_{81}(11,16,36,38)$, $C_{81}(7,20,34,36)$ are Type-1 isomorphic.

 \item [\rm (c1)] $C_{81}(8(4,9,23,31))$ = $C_{81}(5,9,22,32)$ and $C_{81}(10(4,9,23,31))$ = $C_{81}(9,13,14,40)$.  
 
$\Rightarrow$ $C_{81}(4,9,23,31)$, $C_{81}(5,9,22,32)$, $C_{81}(9,13,14,40)$ are Type-1 isomorphic.

 \item [\rm (c2)] $C_{81}(8(4,18,23,31))$ = $C_{81}(5,18,22,32)$ and $C_{81}(10(4,18,23,31))$ = $C_{81}(13,14,18,40)$.  
 
$\Rightarrow$ $C_{81}(4,18,23,31)$, $C_{81}(5,18,22,32)$, $C_{81}(13,14,18,40)$ are Type-1 isomorphic.

  \item [\rm (c3)] $C_{81}(8(4,23,27,31))$ = $C_{81}(5,22,27,32)$ and $C_{81}(10(4,23,27,31))$ = $C_{81}(13,14,27,40)$.  
 
$\Rightarrow$ $C_{81}(4,23,27,31)$, $C_{81}(5,22,27,32)$, $C_{81}(13,14,27,40)$ are Type-1 isomorphic.

  \item [\rm (c4)] $C_{81}(8(4,23,31,36))$ = $C_{81}(5,22,32,36)$ and $C_{81}(10(4,23,31,36))$ = $C_{81}(13,14,36,40)$.  
 
$\Rightarrow$ $C_{81}(4,23,31,36)$, $C_{81}(5,22,32,36)$, $C_{81}(13,14,36,40)$ are Type-1 isomorphic.   \hfill $\Box$
\end{enumerate} 

In the next two problems present Type-2 isomorphic circulant graphs of $C_{81}(r_1,r_2,r_3)$ w.r.t. $m$ = 3.

\begin{prm} \label{p6} {\rm For $s$ = 3,6,12,15,21,24,30,33,39, show that the following triples of circulant graphs are Type-2 isomorphic w.r.t. $m$ = 3.  
\begin{enumerate}
 \item [\rm (d)] $C_{81}(1,s,26,28)$, $C_{81}(s,8,19,35)$, $C_{81}(s,10,17,37)$;  
 \item [\rm (e)] $C_{81}(2,s,25,29)$, $C_{81}(s,7,20,34)$, $C_{81}(s,11,16,38)$; and   
 \item [\rm (f)] $C_{81}(s,4,23,31)$, $C_{81}(s,13,14,40)$, $C_{81}(s,5,22,32)$ where $\gcd(81, s)$ = 3.  
\end{enumerate} }	
\end{prm}
\noindent
{\bf Solution.}\quad Let $R_1$ = $\{1,s,26,28\}$, $S_1$ = $\{s,10,17,37\}$, $T_1$ = $\{s,8,19,35\}$, 

\hspace{2cm} $R_2$ = $\{2,s,25,29\}$, $S_2$ = $\{s,11,16,38\}$, $T_2$ = $\{s,7,20,34\}$, 

\hspace{2cm} $R_3$ = $\{s,4,23,31\}$, $S_3$ = $\{s,5,22,32\}$, $T_3$ =$\{s,13,14,40\}$ where $s$ = 3,6,12,15,21,24,30,33,39 and $\gcd(81, s)$ = 3. 

Here, $n$ = 81 = $3\times 3^3$ and $r$ = s so that $s\in R_i,S_i,T_i$ and $m$ = $\gcd(81, s)$ = 3 where $s$ = 3,6,12,15,21,24,30,33,39 and $i$ = 1,2,3. 

We have $\theta_{81,3,t}(s)$ = s and $\theta_{81,3,t}(81-s)$ = $81-3$, $0 \leq t \leq 26$ and $s$ = 3,6,12,15,21,24,30,33,39. 

This implies, $m$ = 3 is only possible value of $m$ such that given triples of circulant graphs may be Type-2 isomorphic w.r.t. $m$. Using the definition of $\theta_{n,m,t}$, we get, $\theta_{81,3,t}(s)$ = $s$ and $\theta_{81,3,t}(81-s)$ = $81-s$ for $s$ = $3i$, $i$ = 1 to 13 and $0 \leq t \leq 26$. Now, consider different cases.
  \begin{enumerate}
\item [\rm (d)]  For  $s$ = 3,6,12,15,21,24,30,33,39, using the definition of $\theta_{n,m,t}$, we get,  
 \\
 $\theta_{81,3,3}(C_{81}(1,s,26,28))$ = $\theta_{81,3,3}(C_{48}(1,s,26,28, 53,55,81-s,80))$ 

 \hspace{2.7cm} = $\theta_{81,3,3}(C_{81}(1,s,26,28, 53,55,81-s,80))$ 

 \hspace{2.7cm} = $C_{81}(10,s,44,37, 71,64,81-s,17)$ = $C_{81}(s,10,17,37)$ and 
 \\
 $\theta_{81,3,6}(C_{81}(1,s,26,28))$ = $\theta_{81,3,6}(C_{48}(1,s,26,28, 53,55,81-s,80))$ 

 \hspace{2.7cm} = $\theta_{81,3,6}(C_{81}(1,s,26,28, 53,55,81-s,80))$  
 
 \hspace{2.7cm} = $C_{81}(19,s,62,46, 8,73,81-s,35)$ = $C_{81}(s,8,19,35)$.

 $\Rightarrow$ $C_{81}(1,s,26,28)$ $\cong$ $C_{81}(s,10,17,37)$ $\cong$ $C_{81}(s,8,19,35)$ for $s$ = 3,6,12,15,21,24,30,33,39. Also, for $s$ = 3,6,12,15,21,24,30,33,39, we get,
\\ 
$Ad_{81}(C_{81}(1,s,26,28))$ = $Ad_{81}(C_{81}(1,s,26,28, 53,55,81-s,80))$ 
 
 \hspace{2.9cm}  = $\{C_{81}(x(1,s,26,28, 53,55,81-s,80)): x = 1,2,4,5,7,8,10,11,13,14$, 

\hfill $16,17,19,20,22,23,25,26,28,29,31,32,34,35,37,38,40\}$. 
 
 From the above relation, we find, for $s$ = 3,6,12,15,21,24,30,33,39, $Ad_{81}(C_{81}(1,s,26,28))$ and show that $C_{81}(s,10,17,37),C_{81}(s,8,19,35) \notin Ad_{81}(C_{81}(1,s,26,28))$ as follows.
 
\item [\rm (d1)]  $Ad_{81}(C_{81}(1,3,26,28))$ = $\{C_{81}(1,3,26,28)$, $C_{81}(2,6,25,29)$, $C_{81}(4,12,23,31)$, 

$C_{81}(5,15,22,32)$, $C_{81}(7,20,21,34)$, $C_{81}(8,19,24,35)$, $C_{81}(10,17,30,37)$, 

\hfill $C_{81}(11,16,33,38)$, $C_{81}(13,14,39,40)\}$.  

 $\Rightarrow$ $C_{81}(3,10,17,37),C_{81}(3,8,19,35)\notin Ad_{81}(C_{81}(1,3,26,28))$.
 
 $\Rightarrow$ $C_{81}(1,3,26,28)$, $C_{81}(3,10,17,37)$ $\&$ $C_{81}(3,8,19,35)$ are Type-2 isomorphic w.r.t. $m$ = 3.
 
\item [\rm (d2)]  $Ad_{81}(C_{81}(1,6,26,28))$ = $\{C_{81}(1,6,26,28)$, $C_{81}(2,12,25,29)$, $C_{81}(4,23,24,31)$, 

$C_{81}(5,22,30,32)$, $C_{81}(7,20,34,39)$, $C_{81}(8,19,33,35)$, $C_{81}(10,17,21,37)$, 

\hfill $C_{81}(11,15,16,38)$, $C_{81}(3,13,14,40)\}$.  

 $\Rightarrow$ $C_{81}(6,10,17,37),C_{81}(6,8,19,35)\notin Ad_{81}(C_{81}(1,6,26,28))$.
 
 $\Rightarrow$ $C_{81}(1,6,26,28)$, $C_{81}(6,10,17,37)$ $\&$ $C_{81}(6,8,19,35)$ are Type-2 isomorphic w.r.t. $m$ = 3.
 
\item [\rm (d3)]  $Ad_{81}(C_{81}(1,12,26,28))$ = $\{C_{81}(1,12,26,28)$, $C_{81}(2,24,25,29)$, $C_{81}(4,23,31,33)$, 

$C_{81}(5,21,22,32)$, $C_{81}(3,7,20,34)$, $C_{81}(8,15,19,35)$, $C_{81}(10,17,37,39)$, 

\hfill $C_{81}(11,16,30,38)$, $C_{81}(6,13,14,40)\}$.  

 $\Rightarrow$ $C_{81}(10,12,17,37),C_{81}(8,12,19,35)\notin Ad_{81}(C_{81}(1,12,26,28))$.
 \\
 $\Rightarrow$ $C_{81}(1,12,26,28)$, $C_{81}(10,12,17,37)$ $\&$ $C_{81}(8,12,19,35)$ are Type-2 isomorphic w.r.t. $m$ = 3.
 
\item [\rm (d4)]  $Ad_{81}(C_{81}(1,15,26,28))$ = $\{C_{81}(1,15,26,28)$, $C_{81}(2,25,29,30)$, $C_{81}(4,21,23,31)$, 

$C_{81}(5,6,22,32)$, $C_{81}(7,20,24,34)$, $C_{81}(8,19,35,39)$, $C_{81}(10,12,17,37)$, 

\hfill $C_{81}(3,11,16,38)$, $C_{81}(13,14,33,40)\}$.  

 $\Rightarrow$ $C_{81}(10,15,17,37),C_{81}(8,15,19,35)\notin Ad_{81}(C_{81}(1,15,26,28))$.
 \\
 $\Rightarrow$ $C_{81}(1,15,26,28)$, $C_{81}(10,15,17,37)$ $\&$ $C_{81}(8,15,19,35)$ are Type-2 isomorphic w.r.t. $m$ = 3.
 
\item [\rm (d5)]  $Ad_{81}(C_{81}(1,21,26,28))$ = $\{C_{81}(1,21,26,28)$, $C_{81}(2,25,29,39)$, $C_{81}(3,4,23,31)$, 

$C_{81}(5,22,24,32)$, $C_{81}(7,15,20,34)$, $C_{81}(6,8,19,35)$, $C_{81}(10,17,33,37)$, 

\hfill $C_{81}(11,12,16,38)$, $C_{81}(13,14,30,40)\}$.  

 $\Rightarrow$ $C_{81}(10,17,21,37),C_{81}(8,19,21,35)\notin Ad_{81}(C_{81}(1,21,26,28))$.
 \\
 $\Rightarrow$ $C_{81}(1,21,26,28)$, $C_{81}(10,17,21,37)$ $\&$ $C_{81}(8,19,21,35)$ are Type-2 isomorphic w.r.t. $m$ = 3.
 
\item [\rm (d6)]  $Ad_{81}(C_{81}(1,24,26,28))$ = $\{C_{81}(1,24,26,28)$, $C_{81}(2,25,29,33)$, $C_{81}(4,15,23,31)$, 

$C_{81}(5,22,32,39)$, $C_{81}(6,7,20,34)$, $C_{81}(8,19,30,35)$, $C_{81}(3,10,17,37)$, 

\hfill $C_{81}(11,16,21,38)$, $C_{81}(12,13,14,40)\}$.  

 $\Rightarrow$ $C_{81}(10,17,24,37),C_{81}(8,19,24,35)\notin Ad_{81}(C_{81}(1,24,26,28))$.
 \\
 $\Rightarrow$ $C_{81}(1,24,26,28)$, $C_{81}(10,17,24,37)$ $\&$ $C_{81}(8,19,24,35)$ are Type-2 isomorphic w.r.t. $m$ = 3.
 
\item [\rm (d7)]  $Ad_{81}(C_{81}(1,26,28,30))$ = $\{C_{81}(1,26,28,30)$, $C_{81}(2,21,25,29)$, $C_{81}(4,23,31,39)$, 

$C_{81}(5,12,22,32)$, $C_{81}(7,20,33,34)$, $C_{81}(3,8,19,35)$, $C_{81}(10,17,24,37)$, 

\hfill $C_{81}(6,11,16,38)$, $C_{81}(13,14,15,40)\}$.  

 $\Rightarrow$ $C_{81}(10,17,30,37),C_{81}(8,19,30,35)\notin Ad_{81}(C_{81}(1,26,28,30))$.
 \\
 $\Rightarrow$ $C_{81}(1,26,28,30)$, $C_{81}(10,17,30,37)$ $\&$ $C_{81}(8,19,30,35)$ are Type-2 isomorphic w.r.t. $m$ = 3.
 
\item [\rm (d8)]  $Ad_{81}(C_{81}(1,26,28,33))$ = $\{C_{81}(1,26,28,33)$, $C_{81}(2,21,25,29)$, $C_{81}(4,23,31,39)$, 

$C_{81}(5,12,22,32)$, $C_{81}(7,20,33,34)$, $C_{81}(3,8,19,35)$, $C_{81}(10,17,24,37)$, 

\hfill $C_{81}(6,11,16,38)$, $C_{81}(13,14,15,40)\}$.  

 $\Rightarrow$ $C_{81}(10,17,33,37),C_{81}(8,19,33,35)\notin Ad_{81}(C_{81}(1,26,28,33))$.
 \\
 $\Rightarrow$ $C_{81}(1,26,28,33)$, $C_{81}(10,17,33,37)$ $\&$ $C_{81}(8,19,33,35)$ are Type-2 isomorphic w.r.t. $m$ = 3.
 
\item [\rm (d9)]  $Ad_{81}(C_{81}(1,26,28,39))$ = $\{C_{81}(1,26,28,39)$, $C_{81}(2,3,25,29)$, $C_{81}(4,6,23,31)$, 

$C_{81}(5,22,32,33)$, $C_{81}(7,20,30,34)$, $C_{81}(8,12,19,35)$, $C_{81}(10,15,17,37)$, 

\hfill $C_{81}(11,16,24,38)$, $C_{81}(13,14,21,40)\}$.  

 $\Rightarrow$ $C_{81}(10,17,37,39),C_{81}(8,19,35,39)\notin Ad_{81}(C_{81}(1,26,28,39))$.
 \\
 $\Rightarrow$ $C_{81}(1,26,28,39)$, $C_{81}(10,17,37,39)$ $\&$ $C_{81}(8,19,35,39)$ are Type-2 isomorphic w.r.t. $m$ = 3.
 
\item [\rm (e)]  For  $s$ = 3,6,12,15,21,24,30,33,39, using the definition of $\theta_{n,m,t}$, we get,  
 \\
$\theta_{81,3,3}(C_{81}(2,s,25,29))$ = $\theta_{81,3,3}(C_{81}(2,s,25,29, 52,56,81-s,79))$ 
	
	\hspace{3.1cm} = $C_{81}(20,s,34,47, 61,74,81-s,7)$ = $C_{81}(s,7,20,34)$ and 
	\\
	$\theta_{81,3,6}(C_{81}(2,s,25,29))$ = $\theta_{81,3,6}(C_{81}(2,s,25,29, 52,56,81-s,79))$ 
	
	\hspace{3.1cm} = $C_{81}(38,s,43,65, 70,11,81-s,16)$ = $C_{81}(s,11,16,38)$.

 $\Rightarrow$ $C_{81}(2,s,25,29)$ $\cong$ $C_{81}(s,7,20,34)$ $\cong$ $C_{81}(s,11,16,38)$ for $s$ = 3,6,12,15,21,24,30,33,39. Also, for $s$ = 3,6,12,15,21,24,30,33,39, we get,
\\ 
$Ad_{81}(C_{81}(2,s,25,29))$ = $Ad_{81}(C_{81}(2,s,25,29, 52,56,81-s,79))$ 
 
 \hspace{2.9cm}  = $\{C_{81}(x(2,s,25,29, 52,56,81-s,79)): x = 1,2,4,5,7,8,10,11,13,14$, 

\hfill $16,17,19,20,22,23,25,26,28,29,31,32,34,35,37,38,40\}$. 
 
 From the above relation, we find, for $s$ = 3,6,12,15,21,24,30,33,39, $Ad_{81}(C_{81}(2,s,25,29))$ and show that $C_{81}(s,7,20,34),C_{81}(s,11,16,38) \notin Ad_{81}(C_{81}(2,s,25,29))$ as follows.
 
\item [\rm (e1)]  $Ad_{81}(C_{81}(2,3,25,29))$ =  $\{C_{81}(2,3,25,29)$, $C_{81}(4,6,23,31)$, $C_{81}(8,12,19,35)$,   

$C_{81}(10,15,17,37)$, $C_{81}(13,14,21,40)$, $C_{81}(11,16,24,38)$, $C_{81}(7,20,30,34)$,  

\hfill $C_{81}(5,22,32,33)$, $C_{81}(1,26,28,39)\}$. 

 $\Rightarrow$ $C_{81}(3,7,20,34),C_{81}(3,11,16,38)\notin Ad_{81}(C_{81}(2,3,25,29))$.
 
 $\Rightarrow$ $C_{81}(2,3,25,29)$, $C_{81}(3,7,20,34)$ $\&$ $C_{81}(3,11,16,38)$ are Type-2 isomorphic w.r.t. $m$ = 3.

\item [\rm (e2)]  $Ad_{81}(C_{81}(2,6,25,29))$ =  $\{C_{81}(2,6,25,29)$, $C_{81}(4,12,23,31)$, $C_{81}(8,19,24,35)$,   

$C_{81}(10,17,30,37)$, $C_{81}(13,14,39,40)$, $C_{81}(11,16,33,38)$, $C_{81}(7,20,21,34)$,  

\hfill $C_{81}(5,15,22,32)$, $C_{81}(1,3,26,28)\}$. 

 $\Rightarrow$ $C_{81}(6,7,20,34),C_{81}(6,11,16,38)\notin Ad_{81}(C_{81}(2,6,25,29))$.
 
 $\Rightarrow$ $C_{81}(2,6,25,29)$, $C_{81}(6,7,20,34)$ $\&$ $C_{81}(6,11,16,38)$ are Type-2 isomorphic w.r.t. $m$ = 3.

\item [\rm (e3)]  $Ad_{81}(C_{81}(2,12,25,29))$ =  $\{C_{81}(2,12,25,29)$, $C_{81}(4,23,24,31)$, $C_{81}(8,19,33,35)$,   

$C_{81}(10,17,21,37)$, $C_{81}(3,13,14,40)$, $C_{81}(11,15,16,38)$, $C_{81}(7,20,34,39)$,  

\hfill $C_{81}(5,22,30,32)$, $C_{81}(1,6,26,28)\}$. 

 $\Rightarrow$ $C_{81}(7,12,20,34),C_{81}(11,12,16,38)\notin Ad_{81}(C_{81}(2,12,25,29))$.
 \\
 $\Rightarrow$ $C_{81}(2,12,25,29)$, $C_{81}(7,12,20,34)$ $\&$ $C_{81}(11,12,16,38)$ are Type-2 isomorphic w.r.t. $m$ = 3.

\item [\rm (e4)]  $Ad_{81}(C_{81}(2,15,25,29))$ =  $\{C_{81}(2,15,25,29)$, $C_{81}(4,23,30,31)$, $C_{81}(8,19,21,35)$,   

$C_{81}(6,10,17,37)$, $C_{81}(13,14,24,40)$, $C_{81}(11,16,38,39)$, $C_{81}(7,12,20,34)$,  

\hfill $C_{81}(3,5,22,32)$, $C_{81}(1,26,28,33)\}$. 

 $\Rightarrow$ $C_{81}(7,15,20,34),C_{81}(11,15,16,38)\notin Ad_{81}(C_{81}(2,15,25,29))$.
 \\
 $\Rightarrow$ $C_{81}(2,15,25,29)$, $C_{81}(7,15,20,34)$ $\&$ $C_{81}(11,15,16,38)$ are Type-2 isomorphic w.r.t. $m$ = 3.

\item [\rm (e5)]  $Ad_{81}(C_{81}(2,21,25,29))$ =  $\{C_{81}(2,21,25,29)$, $C_{81}(4,23,31,39)$, $C_{81}(3,8,19,35)$,   

$C_{81}(10,17,24,37)$, $C_{81}(13,14,15,40)$, $C_{81}(6,11,16,38)$, $C_{81}(7,20,33,34)$,  

\hfill $C_{81}(5,12,22,32)$, $C_{81}(1,26,28,30)\}$. 

 $\Rightarrow$ $C_{81}(7,20,21,34),C_{81}(11,16,21,38)\notin Ad_{81}(C_{81}(2,21,25,29))$.
 \\
 $\Rightarrow$ $C_{81}(2,21,25,29)$, $C_{81}(7,20,21,34)$ $\&$ $C_{81}(11,16,21,38)$ are Type-2 isomorphic w.r.t. $m$ = 3.

\item [\rm (e6)]  $Ad_{81}(C_{81}(2,24,25,29))$ =  $\{C_{81}(2,24,25,29)$, $C_{81}(4,23,31,33)$, $C_{81}(8,15,19,35)$,   

$C_{81}(10,17,37,39)$, $C_{81}(6,13,14,40)$, $C_{81}(11,16,30,38)$, $C_{81}(3,7,20,34)$,  

\hfill $C_{81}(5,21,22,32)$, $C_{81}(1,12,26,28)\}$. 

 $\Rightarrow$ $C_{81}(7,20,24,34),C_{81}(11,16,24,38)\notin Ad_{81}(C_{81}(2,24,25,29))$.
 \\
 $\Rightarrow$ $C_{81}(2,24,25,29)$, $C_{81}(7,20,24,34)$ $\&$ $C_{81}(11,16,24,38)$ are Type-2 isomorphic w.r.t. $m$ = 3.

\item [\rm (e7)]  $Ad_{81}(C_{81}(2,25,29,30))$ =  $\{C_{81}(2,25,29,30)$, $C_{81}(4,21,23,31)$, $C_{81}(8,19,35,39)$,   

$C_{81}(10,12,17,37)$, $C_{81}(13,14,33,40)$, $C_{81}(3,11,16,38)$, $C_{81}(7,20,24,34)$,  

\hfill $C_{81}(5,6,22,32)$, $C_{81}(1,15,26,28)\}$. 

 $\Rightarrow$ $C_{81}(7,20,30,34),C_{81}(11,16,30,38)\notin Ad_{81}(C_{81}(2,25,29,30))$.
 \\
 $\Rightarrow$ $C_{81}(2,25,29,30)$, $C_{81}(7,20,30,34)$ $\&$ $C_{81}(11,16,30,38)$ are Type-2 isomorphic w.r.t. $m$ = 3.

\item [\rm (e8)]  $Ad_{81}(C_{81}(2,25,29,33))$ =  $\{C_{81}(2,25,29,33)$, $C_{81}(4,15,23,31)$, $C_{81}(8,19,30,35)$,   

$C_{81}(3,10,17,37)$, $C_{81}(12,13,14,40)$, $C_{81}(11,16,21,38)$, $C_{81}(6,7,20,34)$,  

\hfill $C_{81}(5,22,32,39)$, $C_{81}(1,24,26,28)\}$. 

 $\Rightarrow$ $C_{81}(7,20,33,34),C_{81}(11,16,33,38)\notin Ad_{81}(C_{81}(2,25,29,33))$.
 \\
 $\Rightarrow$ $C_{81}(2,25,29,33)$, $C_{81}(7,20,33,34)$ $\&$ $C_{81}(11,16,33,38)$ are Type-2 isomorphic w.r.t. $m$ = 3.

\item [\rm (e9)]  $Ad_{81}(C_{81}(2,25,29,39))$ =  $\{C_{81}(2,25,29,39)$, $C_{81}(3,4,23,31)$, $C_{81}(6,8,19,35)$,   

$C_{81}(10,17,33,37)$, $C_{81}(13,14,30,40)$, $C_{81}(11,12,16,38)$, $C_{81}(7,15,20,34)$,  

\hfill $C_{81}(5,22,24,32)$, $C_{81}(1,21,26,28)\}$. 

 $\Rightarrow$ $C_{81}(7,20,34,39),C_{81}(11,16,38,39)\notin Ad_{81}(C_{81}(2,25,29,39))$.
 \\
 $\Rightarrow$ $C_{81}(2,25,29,39)$, $C_{81}(7,20,34,39)$ $\&$ $C_{81}(11,16,38,39)$ are Type-2 isomorphic w.r.t. $m$ = 3.

 \item [\rm (f)]  For  $s$ = 3,6,12,15,21,24,30,33,39, using the definition of $\theta_{n,m,t}$, we get,  
 \\
$\theta_{81,3,3}(C_{81}(s,4,23,31))$ = $\theta_{81,3,3}(C_{81}(s,4,23,31, 50,58,77,81-s))$ 
	
	\hspace{3.1cm} = $C_{81}(s,13,41,40, 68,67,14,81-s)$ = $C_{81}(s,13,14,40)$ and 
	\\
	$\theta_{81,3,6}(C_{81}(s,4,23,31))$ = $\theta_{81,3,6}(C_{81}(s,4,23,31, 50,58,77,81-s))$ 
	
	\hspace{3.1cm} = $C_{81}(s,22,59,49, 5,76,32,81-s)$ = $C_{81}(s,5,22,32)$.
		
	$\Rightarrow$ $C_{81}(s,4,23,31)$ $\cong$ $C_{81}(s,13,14,40)$ $\cong$ $C_{81}(s,5,22,32)$ for $s$ = 3,6,12,15,21,24,30,33,39. Also, for $s$ = 3,6,12,15,21,24,30,33,39, we get,
\\ 
$Ad_{81}(C_{81}(s,4,23,31))$ = $Ad_{81}(C_{81}(s,4,23,31, 50,58,77,81-s))$ 
 
 \hspace{2.9cm}  = $\{C_{81}(x(s,4,23,31, 50,58,77,81-s)): x = 1,2,4,5,7,8,10,11,13,14$, 

\hfill $16,17,19,20,22,23,25,26,28,29,31,32,34,35,37,38,40\}$. 
 
 From the above relation, we find, for $s$ = 3,6,12,15,21,24,30,33,39, $Ad_{81}(C_{81}(s,4,23,31))$ and show that $C_{81}(s,13,14,40),C_{81}(s,5,22,32) \notin Ad_{81}(C_{81}(s,4,23,31))$ as follows.
 
\item [\rm (f1)]  $Ad_{81}(C_{81}(3,4,23,31))$ =  $\{C_{81}(3,4,23,31)$, $C_{81}(6,8,19,35)$, $C_{81}(11,12,16,38)$,  

  $C_{81}(7,15,20,34)$, $C_{81}(1,21,26,28)$, $C_{81}(5,22,24,32)$, $C_{81}(13,14,30,40)$, 

\hfill $C_{81}(10,17,33,37)$, $C_{81}(2,25,29,39)\}$. 

 $\Rightarrow$ $C_{81}(3,13,14,40),C_{81}(3,5,22,32)\notin Ad_{81}(C_{81}(3,4,23,31))$.
 
 $\Rightarrow$ $C_{81}(3,13,14,40)$, $C_{81}(3,5,22,32)$ $\&$ $C_{81}(3,4,23,31)$ are Type-2 isomorphic w.r.t. $m$ = 3.

\item [\rm (f2)]  $Ad_{81}(C_{81}(4,6,23,31))$ =  $\{C_{81}(4,6,23,31)$, $C_{81}(8,12,19,35)$, $C_{81}(11,16,24,38)$,  

  $C_{81}(7,20,30,34)$, $C_{81}(1,26,28,39)$, $C_{81}(5,22,32,33)$, $C_{81}(13,14,21,40)$, 

\hfill $C_{81}(10,15,17,37)$, $C_{81}(2,3,25,29)\}$. 

 $\Rightarrow$ $C_{81}(6,13,14,40),C_{81}(5,6,22,32)\notin Ad_{81}(C_{81}(4,6,23,31))$.
 
 $\Rightarrow$ $C_{81}(6,13,14,40)$, $C_{81}(5,6,22,32)$ $\&$ $C_{81}(4,6,23,31)$ are Type-2 isomorphic w.r.t. $m$ = 3.

\item [\rm (f3)]  $Ad_{81}(C_{81}(4,12,23,31))$ =  $\{C_{81}(4,12,23,31)$, $C_{81}(8,19,24,35)$, $C_{81}(11,16,33,38)$,  

  $C_{81}(7,20,21,34)$, $C_{81}(1,3,26,28)$, $C_{81}(5,15,22,32)$, $C_{81}(13,14,39,40)$, 

\hfill $C_{81}(10,17,30,37)$, $C_{81}(2,6,25,29)\}$. 

 $\Rightarrow$ $C_{81}(12,13,14,40),C_{81}(5,12,22,32)\notin Ad_{81}(C_{81}(4,12,23,31))$.
 \\
 $\Rightarrow$ $C_{81}(12,13,14,40)$, $C_{81}(5,12,22,32)$ $\&$ $C_{81}(4,12,23,31)$ are Type-2 isomorphic w.r.t. $m$ = 3.

\item [\rm (f4)]  $Ad_{81}(C_{81}(4,15,23,31))$ =  $\{C_{81}(4,15,23,31)$, $C_{81}(8,19,30,35)$, $C_{81}(11,16,21,38)$,  

  $C_{81}(6,7,20,34)$, $C_{81}(1,24,26,28)$, $C_{81}(5,22,32,39)$, $C_{81}(12,13,14,40)$, 

\hfill $C_{81}(3,10,17,37)$, $C_{81}(2,25,29,33)\}$. 

 $\Rightarrow$ $C_{81}(13,14,15,40),C_{81}(5,15,22,32)\notin Ad_{81}(C_{81}(4,15,23,31))$.
 \\
 $\Rightarrow$ $C_{81}(13,14,15,40)$, $C_{81}(5,15,22,32)$ $\&$ $C_{81}(4,15,23,31)$ are Type-2 isomorphic w.r.t. $m$ = 3.

\item [\rm (f5)]  $Ad_{81}(C_{81}(4,21,23,31))$ =  $\{C_{81}(4,21,23,31)$, $C_{81}(8,19,30,35)$, $C_{81}(11,16,21,38)$,  

  $C_{81}(6,7,20,34)$, $C_{81}(1,24,26,28)$, $C_{81}(5,22,32,39)$, $C_{81}(12,13,14,40)$, 

\hfill $C_{81}(3,10,17,37)$, $C_{81}(2,25,29,33)\}$. 

 $\Rightarrow$ $C_{81}(13,14,21,40),C_{81}(5,21,22,32)\notin Ad_{81}(C_{81}(4,21,23,31))$.
 \\
 $\Rightarrow$ $C_{81}(13,14,21,40)$, $C_{81}(5,21,22,32)$ $\&$ $C_{81}(4,21,23,31)$ are Type-2 isomorphic w.r.t. $m$ = 3.

\item [\rm (f6)]  $Ad_{81}(C_{81}(4,23,24,31))$ =  $\{C_{81}(4,23,24,31)$, $C_{81}(8,19,33,35)$, $C_{81}(11,15,16,38)$,  

  $C_{81}(7,20,34,39)$, $C_{81}(1,6,26,28)$, $C_{81}(5,22,30,32)$, $C_{81}(3,13,14,40)$, 

\hfill $C_{81}(10,17,21,37)$, $C_{81}(2,12,25,29)\}$. 

 $\Rightarrow$ $C_{81}(13,14,24,40),C_{81}(5,22,24,32)\notin Ad_{81}(C_{81}(4,23,24,31))$.
 \\
 $\Rightarrow$ $C_{81}(13,14,24,40)$, $C_{81}(5,22,24,32)$ $\&$ $C_{81}(4,23,24,31)$ are Type-2 isomorphic w.r.t. $m$ = 3.

\item [\rm (f7)]  $Ad_{81}(C_{81}(4,23,30,31))$ =  $\{C_{81}(4,23,30,31)$, $C_{81}(8,19,21,35)$, $C_{81}(11,16,38,39)$,  

  $C_{81}(7,12,20,34)$, $C_{81}(1,26,28,33)$, $C_{81}(3,5,22,32)$, $C_{81}(13,14,24,40)$, 

\hfill $C_{81}(6,10,17,37)$, $C_{81}(2,15,25,29)\}$. 

 $\Rightarrow$ $C_{81}(13,14,30,40),C_{81}(5,22,30,32)\notin Ad_{81}(C_{81}(4,23,30,31))$.
 \\
 $\Rightarrow$ $C_{81}(13,14,30,40)$, $C_{81}(5,22,30,32)$ $\&$ $C_{81}(4,23,30,31)$ are Type-2 isomorphic w.r.t. $m$ = 3.

\item [\rm (f8)]  $Ad_{81}(C_{81}(4,23,31,33))$ =  $\{C_{81}(4,23,31,33)$, $C_{81}(8,15,19,35)$, $C_{81}(11,16,30,38)$,  

  $C_{81}(3,7,20,34)$, $C_{81}(1,12,26,28)$, $C_{81}(5,21,22,32)$, $C_{81}(6,13,14,40)$, 

\hfill $C_{81}(10,17,37,39)$, $C_{81}(2,24,25,29)\}$. 

 $\Rightarrow$ $C_{81}(13,14,33,40),C_{81}(5,22,32,33)\notin Ad_{81}(C_{81}(4,23,31,33))$.
 \\
 $\Rightarrow$ $C_{81}(13,14,33,40)$, $C_{81}(5,22,32,33)$ $\&$ $C_{81}(4,23,31,33)$ are Type-2 isomorphic w.r.t. $m$ = 3.

\item [\rm (f9)]  $Ad_{81}(C_{81}(4,23,31,39))$ =  $\{C_{81}(4,23,31,39)$, $C_{81}(3,8,19,35)$, $C_{81}(6,11,16,38)$,  

  $C_{81}(7,20,33,34)$, $C_{81}(1,26,28,30)$, $C_{81}(5,12,22,32)$, $C_{81}(13,14,15,40)$, 

\hfill $C_{81}(10,17,24,37)$, $C_{81}(2,21,25,29)\}$. 

 $\Rightarrow$ $C_{81}(13,14,39,40),C_{81}(5,22,32,39)\notin Ad_{81}(C_{81}(4,23,31,39))$.
 
This implies, $C_{81}(13,14,39,40)$, $C_{81}(5,22,32,39)$ $\&$ $C_{81}(4,23,31,39)$ are Type-2 isomorphic w.r.t. $m$ = 3.     \hfill $\Box$
\end{enumerate} 

\medskip
\noindent
{\bf Type-1 sets of $C_{81}(r_1,r_2,r_3)$ related to problems \ref{p5} and \ref{p6}}

\medskip
We present here the following  Type-1 sets which are related to that occured in problems \ref{p5} and \ref{p6}.

\begin{enumerate}
 \item [\rm (a1)] $Ad_{81}(C_{81}(1,9,26,28))$ = $\{C_{81}(1,9,26,28)$, $C_{81}(2,18,25,29)$, $C_{81}(4,23,31,36)$, $C_{81}(5,22,32,36)$, 

\hfill $C_{81}(7,18,20,34)$, $C_{81}(8,9,19,35)$, $C_{81}(9,10,17,37)$, $C_{81}(11,16,18,38)$, $C_{81}(13,14,36,40)\}$;

\item [\rm (a2)] $Ad_{81}(C_{81}(1,18,26,28))$ = $\{C_{81}(1,18,26,28)$, $C_{81}(2,25,29,36)$, $C_{81}(4,9,23,31)$, $C_{81}(5,9,22,32)$, 

\hfill $C_{81}(7,20,34,36)$, $C_{81}(8,18,19,35)$, $C_{81}(10,17,18,37)$, $C_{81}(11,16,36,38)$, $C_{81}(9,13,14,40)\}$;

\item [\rm (a3)] $Ad_{81}(C_{81}(1,26,27,28))$ = $\{C_{81}(1,26,27,28)$, $C_{81}(2,25,27,29)$, $C_{81}(4,23,27,31)$,

\hfill $C_{81}(5,22,27,33)$,  $C_{81}(7,20,27,34)$, $C_{81}(8,19,27,35)$, $C_{81}(10,17,27,37)$, 

\hfill $C_{81}(11,16,27,38)$, $C_{81}(13,14,27,40)\}$;

\item [\rm (a4)] $Ad_{81}(C_{81}(1,26,28,36))$ = $\{C_{81}(1,26,28,36)$, $C_{81}(2,9,25,29)$, $C_{81}(4,18,23,31)$, $C_{81}(5,18,22,33)$, 

\hfill $C_{81}(7,9,20,34)$, $C_{81}(8,19,35,36)$, $C_{81}(10,17,36,37)$, $C_{81}(9,11,16,38)$, $C_{81}(13,14,18,40)\}$;

\item [\rm (b1)] $Ad_{81}(C_{81}(2,9,25,29))$ = $\{C_{81}(2,9,25,29)$, $C_{81}(4,18,23,31)$, $C_{81}(8,19,35,36)$, 

\hfill $C_{81}(10,17,36,37)$, $C_{81}(13,14,18,40)$, $C_{81}(9,11,16,38)$, $C_{81}(7,9,20,34)\}$;  
 
\item [\rm (b2)] $Ad_{81}(C_{81}(2,18,25,29))$ = $\{C_{81}(2,18,25,29)$, $C_{81}(4,23,31,36)$, $C_{81}(8,9,19,35)$, 

\hfill $C_{81}(9,10,17,37)$, $C_{81}(13,14,36,40)$, $C_{81}(11,16,18,38)$, $C_{81}(7,18,20,34)\}$;  
 
\item [\rm (b3)] $Ad_{81}(C_{81}(2,25,27,29))$ = $\{C_{81}(2,25,27,29)$, $C_{81}(4,23,27,31)$, $C_{81}(8,19,27,35)$, 

\hfill $C_{81}(10,17,27,37)$, $C_{81}(13,14,27,40)$, $C_{81}(11,16,27,38)$, $C_{81}(7,20,27,34)$  
 
\item [\rm (b4)] $Ad_{81}(C_{81}(2,25,29,36))$ = $\{C_{81}(2,25,29,36)$, $C_{81}(4,9,23,31)$, $C_{81}(8,18,19,35)$, 

\hfill $C_{81}(10,17,18,37)$, $C_{81}(9,13,14,40)$, $C_{81}(11,16,36,38)$, $C_{81}(7,20,34,36)\}$;  
 
 \item [\rm (c1)] $Ad_{81}(C_{81}(4,9,23,31))$ = $\{C_{81}(4,9,23,31)$, $C_{81}(8,18,19,35)$, $C_{81}(11,16,36,38)$, 

\hfill $C_{81}(7,20,34,36)$, $C_{81}(1,18,26,28)$, $C_{81}(5,9,22,32)$, $C_{81}(9,13,14,40)$, 

\hfill $C_{81}(10,17,18,37)$, $C_{81}(2,25,29,36)\}$;  
 
 \item [\rm (c2)] $Ad_{81}(C_{81}(4,18,23,31))$ = $\{C_{81}(5,18,22,32)$,  $C_{81}(8,19,35,36)$, $C_{81}(9,11,16,38)$, 

\hfill $C_{81}(7,9,20,34)$, $C_{81}(1,26,28,36)$, $C_{81}(5,18,22,32)$, $C_{81}(13,14,18,40)$, 

\hfill $C_{81}(10,17,36,37)$, $C_{81}(2,9,25,29)\}$;  
 
 \item [\rm (c3)] $Ad_{81}(C_{81}(4,23,27,31))$ = $\{C_{81}(4,23,27,31)$,  $C_{81}(8,19,27,35)$, $C_{81}(11,16,27,38)$, 

\hfill $C_{81}(7,20,27,34)$, $C_{81}(1,26,27,28)$, $C_{81}(5,22,27,32)$, $C_{81}(13,14,27,40)$, 

\hfill $C_{81}(10,17,27,37)$, $C_{81}(2,25,27,29)\}$;  
 
 \item [\rm (c4)] $Ad_{81}(C_{81}(4,23,31,36))$ = $\{C_{81}(4,23,31,36)$,  $C_{81}(8,9,19,35)$, $C_{81}(11,16,18,38)$, 

\hfill $C_{81}(7,18,20,34)$, $C_{81}(1,9,26,28)$, $C_{81}(5,22,32,36)$, $C_{81}(13,14,36,40)$, 

\hfill $C_{81}(9,10,17,37)$, $C_{81}(2,18,25,29)\}$;   

\item [\rm (d1)]  $Ad_{81}(C_{81}(1,3,26,28))$ = $\{C_{81}(1,3,26,28)$, $C_{81}(2,6,25,29)$, $C_{81}(4,12,23,31)$, 

\hfill $C_{81}(5,15,22,32)$, $C_{81}(7,20,21,34)$, $C_{81}(8,19,24,35)$, $C_{81}(10,17,30,37)$, 

\hfill $C_{81}(11,16,33,38)$, $C_{81}(13,14,39,40)\}$.  

\item [\rm (d2)]  $Ad_{81}(C_{81}(1,6,26,28))$ = $\{C_{81}(1,6,26,28)$, $C_{81}(2,12,25,29)$, $C_{81}(4,23,24,31)$, 

\hfill $C_{81}(5,22,30,32)$, $C_{81}(7,20,34,39)$, $C_{81}(8,19,33,35)$, $C_{81}(10,17,21,37)$, 

\hfill $C_{81}(11,15,16,38)$, $C_{81}(3,13,14,40)\}$.  

\item [\rm (d3)]  $Ad_{81}(C_{81}(1,12,26,28))$ = $\{C_{81}(1,12,26,28)$, $C_{81}(2,24,25,29)$, $C_{81}(4,23,31,33)$, 

\hfill $C_{81}(5,21,22,32)$, $C_{81}(3,7,20,34)$, $C_{81}(8,15,19,35)$, $C_{81}(10,17,37,39)$, 

\hfill $C_{81}(11,16,30,38)$, $C_{81}(6,13,14,40)\}$.  

\item [\rm (d4)]  $Ad_{81}(C_{81}(1,15,26,28))$ = $\{C_{81}(1,15,26,28)$, $C_{81}(2,25,29,30)$, $C_{81}(4,21,23,31)$, 

\hfill $C_{81}(5,6,22,32)$, $C_{81}(7,20,24,34)$, $C_{81}(8,19,35,39)$, $C_{81}(10,12,17,37)$, 

\hfill $C_{81}(3,11,16,38)$, $C_{81}(13,14,33,40)\}$.  

\item [\rm (d5)]  $Ad_{81}(C_{81}(1,21,26,28))$ = $\{C_{81}(1,21,26,28)$, $C_{81}(2,25,29,39)$, $C_{81}(3,4,23,31)$, 

\hfill $C_{81}(5,22,24,32)$, $C_{81}(7,15,20,34)$, $C_{81}(6,8,19,35)$, $C_{81}(10,17,33,37)$, 

\hfill $C_{81}(11,12,16,38)$, $C_{81}(13,14,30,40)\}$.  

\item [\rm (d6)]  $Ad_{81}(C_{81}(1,24,26,28))$ = $\{C_{81}(1,24,26,28)$, $C_{81}(2,25,29,33)$, $C_{81}(4,15,23,31)$, 

\hfill $C_{81}(5,22,32,39)$, $C_{81}(6,7,20,34)$, $C_{81}(8,19,30,35)$, $C_{81}(3,10,17,37)$, 

\hfill $C_{81}(11,16,21,38)$, $C_{81}(12,13,14,40)\}$.  

\item [\rm (d7)]  $Ad_{81}(C_{81}(1,26,28,30))$ = $\{C_{81}(1,26,28,30)$, $C_{81}(2,21,25,29)$, $C_{81}(4,23,31,39)$, 

\hfill $C_{81}(5,12,22,32)$, $C_{81}(7,20,33,34)$, $C_{81}(3,8,19,35)$, $C_{81}(10,17,24,37)$, 

\hfill $C_{81}(6,11,16,38)$, $C_{81}(13,14,15,40)\}$.  

\item [\rm (d8)]  $Ad_{81}(C_{81}(1,26,28,33))$ = $\{C_{81}(1,26,28,33)$, $C_{81}(2,21,25,29)$, $C_{81}(4,23,31,39)$, 

\hfill $C_{81}(5,12,22,32)$, $C_{81}(7,20,33,34)$, $C_{81}(3,8,19,35)$, $C_{81}(10,17,24,37)$, 

\hfill $C_{81}(6,11,16,38)$, $C_{81}(13,14,15,40)\}$.  

\item [\rm (d9)]  $Ad_{81}(C_{81}(1,26,28,39))$ = $\{C_{81}(1,26,28,39)$, $C_{81}(2,3,25,29)$, $C_{81}(4,6,23,31)$, 

\hfill $C_{81}(5,22,32,33)$, $C_{81}(7,20,30,34)$, $C_{81}(8,12,19,35)$, $C_{81}(10,15,17,37)$, 

\hfill $C_{81}(11,16,24,38)$, $C_{81}(13,14,21,40)\}$.  

\item [\rm (e1)]  $Ad_{81}(C_{81}(2,3,25,29))$ =  $\{C_{81}(2,3,25,29)$, $C_{81}(4,6,23,31)$, $C_{81}(8,12,19,35)$,   

\hfill $C_{81}(10,15,17,37)$, $C_{81}(13,14,21,40)$, $C_{81}(11,16,24,38)$, $C_{81}(7,20,30,34)$,  

\hfill $C_{81}(5,22,32,33)$, $C_{81}(1,26,28,39)\}$. 

\item [\rm (e2)]  $Ad_{81}(C_{81}(2,6,25,29))$ =  $\{C_{81}(2,6,25,29)$, $C_{81}(4,12,23,31)$, $C_{81}(8,19,24,35)$,   

\hfill $C_{81}(10,17,30,37)$, $C_{81}(13,14,39,40)$, $C_{81}(11,16,33,38)$, $C_{81}(7,20,21,34)$,  

\hfill $C_{81}(5,15,22,32)$, $C_{81}(1,3,26,28)\}$. 

\item [\rm (e3)]  $Ad_{81}(C_{81}(2,12,25,29))$ =  $\{C_{81}(2,12,25,29)$, $C_{81}(4,23,24,31)$, $C_{81}(8,19,33,35)$,   

\hfill $C_{81}(10,17,21,37)$, $C_{81}(3,13,14,40)$, $C_{81}(11,15,16,38)$, $C_{81}(7,20,34,39)$,  

\hfill $C_{81}(5,22,30,32)$, $C_{81}(1,6,26,28)\}$. 

\item [\rm (e4)]  $Ad_{81}(C_{81}(2,15,25,29))$ =  $\{C_{81}(2,15,25,29)$, $C_{81}(4,23,30,31)$, $C_{81}(8,19,21,35)$,   

\hfill $C_{81}(6,10,17,37)$, $C_{81}(13,14,24,40)$, $C_{81}(11,16,38,39)$, $C_{81}(7,12,20,34)$,  

\hfill $C_{81}(3,5,22,32)$, $C_{81}(1,26,28,33)\}$. 

\item [\rm (e5)]  $Ad_{81}(C_{81}(2,21,25,29))$ =  $\{C_{81}(2,21,25,29)$, $C_{81}(4,23,31,39)$, $C_{81}(3,8,19,35)$,   

\hfill $C_{81}(10,17,24,37)$, $C_{81}(13,14,15,40)$, $C_{81}(6,11,16,38)$, $C_{81}(7,20,33,34)$,  

\hfill $C_{81}(5,12,22,32)$, $C_{81}(1,26,28,30)\}$. 

\item [\rm (e6)]  $Ad_{81}(C_{81}(2,24,25,29))$ =  $\{C_{81}(2,24,25,29)$, $C_{81}(4,23,31,33)$, $C_{81}(8,15,19,35)$,   

\hfill $C_{81}(10,17,37,39)$, $C_{81}(6,13,14,40)$, $C_{81}(11,16,30,38)$, $C_{81}(3,7,20,34)$,  

\hfill $C_{81}(5,21,22,32)$, $C_{81}(1,12,26,28)\}$. 

\item [\rm (e7)]  $Ad_{81}(C_{81}(2,25,29,30))$ =  $\{C_{81}(2,25,29,30)$, $C_{81}(4,21,23,31)$, $C_{81}(8,19,35,39)$,   

\hfill $C_{81}(10,12,17,37)$, $C_{81}(13,14,33,40)$, $C_{81}(3,11,16,38)$, $C_{81}(7,20,24,34)$,  

\hfill $C_{81}(5,6,22,32)$, $C_{81}(1,15,26,28)\}$. 

\item [\rm (e8)]  $Ad_{81}(C_{81}(2,25,29,33))$ =  $\{C_{81}(2,25,29,33)$, $C_{81}(4,15,23,31)$, $C_{81}(8,19,30,35)$,   

\hfill $C_{81}(3,10,17,37)$, $C_{81}(12,13,14,40)$, $C_{81}(11,16,21,38)$, $C_{81}(6,7,20,34)$,  

\hfill $C_{81}(5,22,32,39)$, $C_{81}(1,24,26,28)\}$. 

\item [\rm (e9)]  $Ad_{81}(C_{81}(2,25,29,39))$ =  $\{C_{81}(2,25,29,39)$, $C_{81}(3,4,23,31)$, $C_{81}(6,8,19,35)$,   

\hfill $C_{81}(10,17,33,37)$, $C_{81}(13,14,30,40)$, $C_{81}(11,12,16,38)$, $C_{81}(7,15,20,34)$,  

\hfill $C_{81}(5,22,24,32)$, $C_{81}(1,21,26,28)\}$. 

\item [\rm (f1)]  $Ad_{81}(C_{81}(3,4,23,31))$ =  $\{C_{81}(3,4,23,31)$, $C_{81}(6,8,19,35)$, $C_{81}(11,12,16,38)$,  

\hfill   $C_{81}(7,15,20,34)$, $C_{81}(1,21,26,28)$, $C_{81}(5,22,24,32)$, $C_{81}(13,14,30,40)$, 

\hfill $C_{81}(10,17,33,37)$, $C_{81}(2,25,29,39)\}$. 

\item [\rm (f2)]  $Ad_{81}(C_{81}(4,6,23,31))$ =  $\{C_{81}(4,6,23,31)$, $C_{81}(8,12,19,35)$, $C_{81}(11,16,24,38)$,  

\hfill   $C_{81}(7,20,30,34)$, $C_{81}(1,26,28,39)$, $C_{81}(5,22,32,33)$, $C_{81}(13,14,21,40)$, 

\hfill $C_{81}(10,15,17,37)$, $C_{81}(2,3,25,29)\}$. 

\item [\rm (f3)]  $Ad_{81}(C_{81}(4,12,23,31))$ =  $\{C_{81}(4,12,23,31)$, $C_{81}(8,19,24,35)$, $C_{81}(11,16,33,38)$,  

\hfill   $C_{81}(7,20,21,34)$, $C_{81}(1,3,26,28)$, $C_{81}(5,15,22,32)$, $C_{81}(13,14,39,40)$, 

\item [\rm (f4)]  $Ad_{81}(C_{81}(4,15,23,31))$ =  $\{C_{81}(4,15,23,31)$, $C_{81}(8,19,30,35)$, $C_{81}(11,16,21,38)$,  

\hfill   $C_{81}(6,7,20,34)$, $C_{81}(1,24,26,28)$, $C_{81}(5,22,32,39)$, $C_{81}(12,13,14,40)$, 

\hfill $C_{81}(3,10,17,37)$, $C_{81}(2,25,29,33)\}$. 

\item [\rm (f5)]  $Ad_{81}(C_{81}(4,21,23,31))$ =  $\{C_{81}(4,21,23,31)$, $C_{81}(8,19,30,35)$, $C_{81}(11,16,21,38)$,  

\hfill   $C_{81}(6,7,20,34)$, $C_{81}(1,24,26,28)$, $C_{81}(5,22,32,39)$, $C_{81}(12,13,14,40)$, 

\hfill $C_{81}(3,10,17,37)$, $C_{81}(2,25,29,33)\}$. 

\item [\rm (f6)]  $Ad_{81}(C_{81}(4,23,24,31))$ =  $\{C_{81}(4,23,24,31)$, $C_{81}(8,19,33,35)$, $C_{81}(11,15,16,38)$,  

\hfill   $C_{81}(7,20,34,39)$, $C_{81}(1,6,26,28)$, $C_{81}(5,22,30,32)$, $C_{81}(3,13,14,40)$, 

\hfill $C_{81}(10,17,21,37)$, $C_{81}(2,12,25,29)\}$. 

\item [\rm (f7)]  $Ad_{81}(C_{81}(4,23,30,31))$ =  $\{C_{81}(4,23,30,31)$, $C_{81}(8,19,21,35)$, $C_{81}(11,16,38,39)$,  

\hfill   $C_{81}(7,12,20,34)$, $C_{81}(1,26,28,33)$, $C_{81}(3,5,22,32)$, $C_{81}(13,14,24,40)$, 

\hfill $C_{81}(6,10,17,37)$, $C_{81}(2,15,25,29)\}$. 

\item [\rm (f8)]  $Ad_{81}(C_{81}(4,23,31,33))$ =  $\{C_{81}(4,23,31,33)$, $C_{81}(8,15,19,35)$, $C_{81}(11,16,30,38)$,  

\hfill   $C_{81}(3,7,20,34)$, $C_{81}(1,12,26,28)$, $C_{81}(5,21,22,32)$, $C_{81}(6,13,14,40)$, 

\hfill $C_{81}(10,17,37,39)$, $C_{81}(2,24,25,29)\}$. 

\item [\rm (f9)]  $Ad_{81}(C_{81}(4,23,31,39))$ =  $\{C_{81}(4,23,31,39)$, $C_{81}(3,8,19,35)$, $C_{81}(6,11,16,38)$,  

\hfill   $C_{81}(7,20,33,34)$, $C_{81}(1,26,28,30)$, $C_{81}(5,12,22,32)$, $C_{81}(13,14,15,40)$, 

\hfill $C_{81}(10,17,24,37)$, $C_{81}(2,21,25,29)\}$.
\end{enumerate}

\subsection{On isomorphic circulant graphs $C_{96}(r_1,r_2,r_3)$ of Type-2 w.r.t. $m$ = 2}

In this section, we find all pairs of Type-2 isomorphic circulant graphs of the form $C_{96}(r_1,r_2,r_3)$ and $C_{96}(s_1,s_2,s_3)$ and the total number of such pairs are 72 and all are of Type-2 w.r.t. $m$ = 2. Here, 96 = $2^2 \times 3\times 2^3$ and so $m$ = 2 is the only possible value of $m$ such that $C_{96}(r_1,r_2,r_3)$ is having Type-2 isomorphic circulant graph(s) w.r.t. $m$. We start with presenting some  pairs of isomorphic circulant graphs $C_{96}(r_1,r_2,r_3)$ and $C_{96}(s_1,s_2,s_3)$ of Type-1.

 \begin{prm} \label{p7} {\rm For $s$ = 4,8,12,16,20,24,28,32,36,40,44,48, show that the following pairs of circulant graphs are Type-1.  
 \begin{enumerate}
 	\item [\rm (a)]  $C_{96}(1,s,47)$ and $C_{48}(s,23,25)$;  
 	\item [\rm (b)]  $C_{96}(3,s,45)$ and $C_{48}(s,21,27)$;  
 	\item [\rm (c)]  $C_{96}(5,s,43)$ and $C_{48}(s,19,29)$;  
 	\item [\rm (d)]  $C_{96}(7,s,41)$ and $C_{48}(s,17,31)$;  
 	\item [\rm (e)]  $C_{96}(9,s,39)$ and $C_{48}(s,15,33)$; and  
 	\item [\rm (f)]  $C_{96}(11,s,37)$ and $C_{48}(s,13,35)$ where $\gcd(48, s) > 2$.  
 \end{enumerate} }	
\end{prm}
 \noindent
 {\bf Solution.}\quad For $s$ = 4,8,12,16,20,24,28,32,36,40,44,48, given pairs of circulant graphs are isomorphic of Type-1 by the following.
 \begin{enumerate}
\item [\rm (a1)]  $C_{96}(23(1,4,47))$ = $C_{96}(23(1,4,47, 49,92,95))$ = $C_{48}(23,92,25, 71,4,73)$ = $C_{96}((4,23,25))$;  

\item [\rm (a2)]  $C_{96}(23(1,8,47))$ = $C_{96}(23(1,8,47, 49,88,95))$ = $C_{48}(23,88,25, 71,8,73)$ = $C_{96}((8,23,25))$;  

\item [\rm (a3)]  $C_{96}(23(1,12,47))$ = $C_{96}(23(1,12,47, 49,84,95))$ = $C_{48}(23,84,25, 71,12,73)$ = $C_{96}((12,23,25))$;  

\item [\rm (a4)]  $C_{96}(23(1,16,47))$ = $C_{96}(23(1,16,47, 49,80,95))$ = $C_{48}(23,80,25, 71,16,73)$ = $C_{96}((16,23,25))$;  

\item [\rm (a5)]  $C_{96}(23(1,20,47))$ = $C_{96}(23(1,20,47, 49,76,95))$ = $C_{48}(23,76,25, 71,20,73)$ = $C_{96}((20,23,25))$;  

\item [\rm (a6)]  $C_{96}(23(1,24,47))$ = $C_{96}(23(1,24,47, 49,72,95))$ = $C_{48}(23,72,25, 71,24,73)$ = $C_{96}((23,24,25))$;  

\item [\rm (a7)]  $C_{96}(23(1,28,47))$ = $C_{96}(23(1,28,47, 49,68,95))$ = $C_{48}(23,68,25, 71,28,73)$ = $C_{96}((23,25,28))$;  

\item [\rm (a8)]  $C_{96}(23(1,32,47))$ = $C_{96}(23(1,32,47, 49,64,95))$ = $C_{48}(23,64,25, 71,32,73)$ = $C_{96}((23,25,32))$;  

\item [\rm (a9)]  $C_{96}(23(1,36,47))$ = $C_{96}(23(1,36,47, 49,60,95))$ = $C_{48}(23,60,25, 71,36,73)$ = $C_{96}((23,25,36))$;  

\item [\rm (a10)]  $C_{96}(23(1,40,47))$ = $C_{96}(23(1,40,47, 49,56,95))$ = $C_{48}(23,56,25, 71,40,73)$ = $C_{96}((23,25,40))$;  

\item [\rm (a11)]  $C_{96}(23(1,44,47))$ = $C_{96}(23(1,44,47, 49,52,95))$ = $C_{48}(23,52,25, 71,44,73)$ = $C_{96}((23,25,44))$;  

\item [\rm (a12)]  $C_{96}(23(1,47,48))$ = $C_{96}(23(1,47, 48, 49,95))$ = $C_{48}(23,25, 48, 71,73)$ = $C_{96}((23,25,48))$;  

\item [\rm (b1)]  $C_{96}(23(3,4,45))$ = $C_{96}(23(3,4,45, 51,92,93))$ = $C_{96}(69,92,75, 21,4,27)$ = $C_{48}(4,21,27)$;  

\item [\rm (b2)]  $C_{96}(23(3,8,45))$ = $C_{96}(23(3,8,45, 51,88,93))$ = $C_{96}(69,88,75, 21,8,27)$ = $C_{48}(8,21,27)$;  

\item [\rm (b3)]  $C_{96}(23(3,12,45))$ = $C_{96}(23(3,12,45, 51,84,93))$ = $C_{96}(69,84,75, 21,12,27)$ = $C_{48}(12,21,27)$;  

\item [\rm (b4)]  $C_{96}(23(3,16,45))$ = $C_{96}(23(3,16,45, 51,80,93))$ = $C_{96}(69,80,75, 21,16,27)$ = $C_{48}(16,21,27)$;  

\item [\rm (b5)]  $C_{96}(23(3,20,45))$ = $C_{96}(23(3,20,45, 51,76,93))$ = $C_{96}(69,76,75, 21,20,27)$ = $C_{48}(20,21,27)$;  

\item [\rm (b6)]  $C_{96}(23(3,24,45))$ = $C_{96}(23(3,24,45, 51,72,93))$ = $C_{96}(69,72,75, 21,24,27)$ = $C_{48}(21,24,27)$;  

\item [\rm (b7)]  $C_{96}(23(3,28,45))$ = $C_{96}(23(3,28,45, 51,68,93))$ = $C_{96}(69,68,75, 21,28,27)$ = $C_{48}(21,27,28)$;  

\item [\rm (b8)]  $C_{96}(23(3,32,45))$ = $C_{96}(23(3,32,45, 51,64,93))$ = $C_{96}(69,64,75, 21,32,27)$ = $C_{48}(21,27,32)$;  

\item [\rm (b9)]  $C_{96}(23(3,36,45))$ = $C_{96}(23(3,36,45, 51,60,93))$ = $C_{96}(69,60,75, 21,36,27)$ = $C_{48}(21,27,36)$;  

\item [\rm (b10)]  $C_{96}(23(3,40,45))$ = $C_{96}(23(3,40,45, 51,56,93))$ = $C_{96}(69,56,75, 21,40,27)$ = $C_{48}(21,27,40)$;  

\item [\rm (b11)]  $C_{96}(23(3,44,45))$ = $C_{96}(23(3,44,45, 51,52,93))$ = $C_{96}(69,52,75, 21,44,27)$ = $C_{48}(21,27,44)$;  

\item [\rm (b12)]  $C_{96}(23(3,48,45))$ = $C_{96}(23(3,48,45, 51,48,93))$ = $C_{96}(69,48,75, 21,48,27)$ = $C_{48}(21,27,48)$;  

\item [\rm (c1)]  $C_{96}(23(4,5,43))$ = $C_{96}(23(4,5,43, 53,91,92))$ = $C_{96}(92,19,29, 67,77,4)$ = $C_{48}(4,19,29)$;  

\item [\rm (c2)]  $C_{96}(23(5,8,43))$ = $C_{96}(23(5,8,43, 53,88,91))$ = $C_{96}(19,88,29, 67,8,77)$ = $C_{48}(8,19,29)$;  

\item [\rm (c3)]  $C_{96}(23(5,12,43))$ = $C_{96}(23(5,12,43, 53,84,91))$ = $C_{96}(19,84,29, 67,12,77)$ = $C_{48}(12,19,29)$;  

\item [\rm (c4)]  $C_{96}(23(5,16,43))$ = $C_{96}(23(5,16,43, 53,80,91))$ = $C_{96}(19,80,29, 67,16,77)$ = $C_{48}(16,19,29)$;  

\item [\rm (c5)]  $C_{96}(23(5,20,43))$ = $C_{96}(23(5,20,43, 53,76,91))$ = $C_{96}(19,76,29, 67,20,77)$ = $C_{48}(19,20,29)$;  

\item [\rm (c6)]  $C_{96}(23(5,24,43))$ = $C_{96}(23(5,24,43, 53,72,91))$ = $C_{96}(19,72,29, 67,24,77)$ = $C_{48}(19,24,29)$;  

\item [\rm (c7)]  $C_{96}(23(5,28,43))$ = $C_{96}(23(5,28,43, 53,68,91))$ = $C_{96}(19,68,29, 67,28,77)$ = $C_{48}(19,28,29)$;  

\item [\rm (c8)]  $C_{96}(23(5,32,43))$ = $C_{96}(23(5,32,43, 53,64,91))$ = $C_{96}(19,64,29, 67,32,77)$ = $C_{48}(19,29,32)$;  

\item [\rm (c9)]  $C_{96}(23(5,36,43))$ = $C_{96}(23(5,36,43, 53,60,91))$ = $C_{96}(19,60,29, 67,36,77)$ = $C_{48}(19,29,36)$;  

\item [\rm (c10)]  $C_{96}(23(5,40,43))$ = $C_{96}(23(5,40,43, 53,56,91))$ = $C_{96}(19,56,29, 67,40,77)$ = $C_{48}(19,29,40)$;  

\item [\rm (c11)]  $C_{96}(23(5,43,44))$ = $C_{96}(23(5,43,44, 53,52,91))$ = $C_{96}(19,52,29, 67,44,77)$ = $C_{48}(19,29,44)$;  

\item [\rm (c12)]  $C_{96}(23(5,43,48))$ = $C_{96}(23(5,43, 48, 53,91))$ = $C_{96}(19,29, 48, 67,77)$ = $C_{48}(19,29,48)$;  

\item [\rm (d1)]  $C_{96}(23(4,7,41))$ = $C_{96}(23(4,7,41, 55,89,92))$ = $C_{96}(92,65,79, 17,31,4)$ = $C_{48}(4,17,31)$;  

\item [\rm (d2)]  $C_{96}(23(7,8,41))$ = $C_{96}(23(7,8,41, 55,88,89))$ = $C_{96}(65,88,79, 17,8,31)$ = $C_{48}(8,17,31)$;  

\item [\rm (d3)]  $C_{96}(23(7,12,41))$ = $C_{96}(23(7,12,41, 55,84,89))$ = $C_{96}(65,84,79, 17,12,31)$ = $C_{48}(12,17,31)$;  

\item [\rm (d4)]  $C_{96}(23(7,16,41))$ = $C_{96}(23(7,16,41, 55,80,89))$ = $C_{96}(65,80,79, 17,16,31)$ = $C_{48}(16,17,31)$;  

\item [\rm (d5)]  $C_{96}(23(7,20,41))$ = $C_{96}(23(7,20,41, 55,76,89))$ = $C_{96}(65,76,79, 17,20,31)$ = $C_{48}(17,20,31)$;  

\item [\rm (d6)]  $C_{96}(23(7,24,41))$ = $C_{96}(23(7,24,41, 55,72,89))$ = $C_{96}(65,72,79, 17,24,31)$ = $C_{48}(17,24,31)$;  

\item [\rm (d7)]  $C_{96}(23(7,28,41))$ = $C_{96}(23(7,28,41, 55,68,89))$ = $C_{96}(65,68,79, 17,28,31)$ = $C_{48}(17,28,31)$;  

\item [\rm (d8)]  $C_{96}(23(7,32,41))$ = $C_{96}(23(7,32,41, 55,64,89))$ = $C_{96}(65,64,79, 17,32,31)$ = $C_{48}(17,31,32)$;  

\item [\rm (d9)]  $C_{96}(23(7,36,41))$ = $C_{96}(23(7,36,41, 55,60,89))$ = $C_{96}(65,60,79, 17,36,31)$ = $C_{48}(17,31,36)$;  

\item [\rm (d10)]  $C_{96}(23(7,40,41))$ = $C_{96}(23(7,40,41, 55,56,89))$ = $C_{96}(65,56,79, 17,40,31)$ = $C_{48}(17,31,40)$;  

\item [\rm (d11)]  $C_{96}(23(7,41,44))$ = $C_{96}(23(7,41,44, 52,55,89))$ = $C_{96}(65,79,52, 44,17,31)$ = $C_{48}(17,31,44)$;  

\item [\rm (d12)]  $C_{96}(23(7,41,48))$ = $C_{96}(23(7,41, 48, 55,89))$ = $C_{96}(65,79, 48, 17,31)$ = $C_{48}(17,31,48)$;  

\item [\rm (e1)]  $C_{96}(23(4,9,39))$ = $C_{96}(23(4,9,39, 57,87,92))$ = $C_{96}(92,15,33, 63,81,4)$ = $C_{48}(4,15,33)$;  

\item [\rm (e2)]  $C_{96}(23(8,9,39))$ = $C_{96}(23(8,9,39, 57,87,88))$ = $C_{96}(88,15,33, 63,81,8)$ = $C_{48}(8,15,33)$;  

\item [\rm (e3)]  $C_{96}(23(9,12,39))$ = $C_{96}(23(9,12,39, 57,84,87))$ = $C_{96}(15,84,33, 63,12,81)$ = $C_{48}(12,15,33)$;  

\item [\rm (e4)]  $C_{96}(23(9,16,39))$ = $C_{96}(23(9,16,39, 57,80,87))$ = $C_{96}(15,80,33, 63,16,81)$ = $C_{48}(15,16,33)$;  

\item [\rm (e5)]  $C_{96}(23(9,20,39))$ = $C_{96}(23(9,20,39, 57,76,87))$ = $C_{96}(15,76,33, 63,20,81)$ = $C_{48}(15,20,33)$;  

\item [\rm (e6)]  $C_{96}(23(9,24,39))$ = $C_{96}(23(9,24,39, 57,72,87))$ = $C_{96}(15,72,33, 63,24,81)$ = $C_{48}(15,24,33)$;  

\item [\rm (e7)]  $C_{96}(23(9,28,39))$ = $C_{96}(23(9,28,39, 57,68,87))$ = $C_{96}(15,68,33, 63,28,81)$ = $C_{48}(15,28,33)$;  

\item [\rm (e8)]  $C_{96}(23(9,32,39))$ = $C_{96}(23(9,32,39, 57,64,87))$ = $C_{96}(15,64,33, 63,32,81)$ = $C_{48}(15,32,33)$;  

\item [\rm (e9)]  $C_{96}(23(9,36,39))$ = $C_{96}(23(9,36,39, 57,60,87))$ = $C_{96}(15,60,33, 63,36,81)$ = $C_{48}(15,33,36)$;  

\item [\rm (e10)]  $C_{96}(23(9,39,40))$ = $C_{96}(23(9,39,40, 56,57,87))$ = $C_{96}(15,33,56, 40,63,81)$ = $C_{48}(15,33,40)$;  

\item [\rm (e11)]  $C_{96}(23(9,39,44))$ = $C_{96}(23(9,39,44, 52,57,87))$ = $C_{96}(15,33,52, 44,63,81)$ = $C_{48}(15,33,44)$;  

\item [\rm (e12)]  $C_{96}(23(9,39,48))$ = $C_{96}(23(9,39, 48, 57,87))$ = $C_{96}(15,33, 48, 63,81)$ = $C_{48}(15,33,48)$;  

\item [\rm (f1)]  $C_{96}(23(4,11,37))$ = $C_{96}(23(4,11,37, 59,85,92))$ = $C_{96}(92,61,83, 13,35,4)$ = $C_{48}(4,13,35)$;  

\item [\rm (f2)]  $C_{96}(23(8,11,37))$ = $C_{96}(23(8,11,37, 59,85,88))$ = $C_{96}(88,61,83, 13,35,8)$ = $C_{48}(8,13,35)$;  

\item [\rm (f3)]  $C_{96}(23(11,12,37))$ = $C_{96}(23(11,12,37, 59,84,85))$ = $C_{96}(61,84,83, 13,12,35)$ = $C_{48}(12,13,35)$;  

\item [\rm (f4)]  $C_{96}(23(11,16,37))$ = $C_{96}(23(11,16,37, 59,80,85))$ = $C_{96}(61,80,83, 13,16,35)$ = $C_{48}(13,16,35)$;  

\item [\rm (f5)]  $C_{96}(23(11,20,37))$ = $C_{96}(23(11,20,37, 59,76,85))$ = $C_{96}(61,76,83, 13,20,35)$ = $C_{48}(13,20,35)$;  

\item [\rm (f6)]  $C_{96}(23(11,24,37))$ = $C_{96}(23(11,24,37, 59,72,85))$ = $C_{96}(61,72,83, 13,24,35)$ = $C_{48}(13,24,35)$;  

\item [\rm (f7)]  $C_{96}(23(11,28,37))$ = $C_{96}(23(11,28,37, 59,68,85))$ = $C_{96}(61,68,83, 13,28,35)$ = $C_{48}(13,28,35)$;  

\item [\rm (f8)]  $C_{96}(23(11,32,37))$ = $C_{96}(23(11,32,37, 59,64,85))$ = $C_{96}(61,64,83, 13,32,35)$ = $C_{48}(13,32,35)$;  

\item [\rm (f9)]  $C_{96}(23(11,36,37))$ = $C_{96}(23(11,36,37, 59,60,85))$ = $C_{96}(61,60,83, 13,36,35)$ = $C_{48}(13,35,36)$;  

 \item [\rm (f10)]  $C_{96}(23(11,37,40))$ = $C_{96}(23(11,37,40, 56,59,85))$ = $C_{96}(61,83,56, 40,13,35)$ = $C_{48}(13,35,40)$;  

 \item [\rm (f11)]  $C_{96}(23(11,37,44))$ = $C_{96}(23(11,37,44, 52,59,85))$ = $C_{96}(61,52,83, 13,44,35)$ = $C_{48}(13,35,44)$;  

\item [\rm (f12)]  $C_{96}(23(11,37,48))$ = $C_{96}(23(11,37, 48, 59,85))$ = $C_{96}(61,83, 48, 13,35)$ = $C_{48}(13,35,48)$.  
\end{enumerate}

In the next two problems, we obtain all circulant graphs Type-2 isomorphic to $C_{96}(r_1,r_2,r_3)$ and all are Type-2 isomorphic w.r.t. $m$ = 2.
  
 \begin{prm} \label{p8} {\rm For $s$ = 2,10,14,22,26,34,38,46, show that the following pairs of circulant graphs are Type-2 isomorphic w.r.t. $m$ = 2.  
 \begin{enumerate}
\item [\rm (g)]  $C_{96}(1,s,47)$ and $C_{48}(s,23,25)$;  
\item [\rm (h)]  $C_{96}(3,s,45)$ and $C_{48}(s,21,27)$;  
\item [\rm (i)]  $C_{96}(5,s,43)$ and $C_{48}(s,19,29)$;  
\item [\rm (j)]  $C_{96}(7,s,41)$ and $C_{48}(s,17,31)$;  
\item [\rm (k)]  $C_{96}(9,s,39)$ and $C_{48}(s,15,33)$; and  
\item [\rm (l)]  $C_{96}(11,s,37)$ and $C_{48}(s,13,35)$ 
where $\gcd(48, s)$ = 2.  
 \end{enumerate} }	
\end{prm}
 \noindent
 {\bf Solution.}\quad For $s$ = 2, 10, 14, 22, 26, 34, 38, 46, let

$R_1$ = $\{1,s,47\}$, $R_2$ = $\{3,s,45\}$, $R_3$ = $\{5,s,43\}$, $R_4$ = $\{7,s,41\}$, $R_5$ = $\{9,s,39\}$, $R_6$ = $\{11,s,37\}$,  

$S_1$ = $\{s,23,25\}$, $S_2$ = $\{s,21,27\}$, $S_3$ = $\{s,19,29\}$, $S_4$ = $\{s,17,31\}$, $S_5$ = $\{s,15,33\}$, $S_6$ = $\{s,13,35\}$. Here, $n$ = 96 and $r$ = $s$ so that $s\in R_i,S_i$ and $\gcd(n, s)$ = $\gcd(96, 2)$ = 2 = $\gcd(96,s )$ for $s$ = 2, 10, 14, 22, 26, 34, 38, 46 and $i$ = 1 to 6. Using the definition of $\theta_{n,m,t}$, we get, $\theta_{96,2,12}(s)$ = $s$ and $\theta_{96,2,12}(96-s)$ = $96-s$. Now, consider different cases.
\begin{enumerate}
\item [\rm (g)]  For  $s$ = 2, 10, 14, 22, 26, 34, 38, 46, using the definition of $\theta_{n,m,t}$, we get,  
 \\
 $\theta_{96,2,12}(C_{96}(1,s,47))$ = $\theta_{96,2,12}(C_{96}(1,s,47, 49,96-s,95))$ = $C_{96}(\theta_{96,2,12}(1,s,47, 49,96-s,95))$ 
 
 \hspace{2.8cm} = $C_{96}(25,s,71, 73,96-s,23)$ = $C_{48}(s,23,25)$. 
 
 $\Rightarrow$ $C_{96}(1,s,47)$ $\cong$ $C_{96}(s,23,25)$ for $s$ = 2, 10, 14, 22, 26, 34, 38, 46. Also, for $s$ = 2, 10, 14, 22, 26, 34, 38, 46, we get,
\\ 
$Ad_{96}(C_{96}(1,s,47)$ = $Ad_{96}(C_{96}(1,s,47, 49,96-s,95))$ 
  = $\{C_{96}(x(1,s,47, 49,96-s,95)): x = 1$, 
$5,7,11,13,17,19,23,25,29,31,35,37,41,43,47,49,53,55,59,61,65,67,71,77,79,83,85,89,91\}$. 
 
 From the above relation, we find $Ad_{96}(C_{96}(1,s,47))$ for $s$ = 2, 10, 14, 22, 26, 34, 38, 46 and will show that $C_{96}(s,23,25) \notin Ad_{96}(C_{96}(1,s,47))$.
 
\item [\rm (g1)]  $Ad_{96}(C_{96}(1,2,47))$ = $\{C_{96}(1,2,47)$, $C_{96}(5,10,43)$, $C_{96}(7,14,41)$, $C_{96}(11,22,37)$, 
 
 \hfill  $C_{96}(13,26,35)$, $C_{96}(17,31,34)$, $C_{96}(19,29,38)$, $C_{96}(23,25,46)\}$.  

 $\Rightarrow$ $C_{96}(2,23,25) \notin Ad_{96}(C_{96}(1,2,47))$.
 
 $\Rightarrow$ $C_{96}(1,2,47)$ and $C_{96}(2,23,25)$ are Type-2 isomorphic w.r.t. $m$ = 2.
 
\item [\rm (g2)]  $Ad_{96}(C_{96}(1,10,47))$ = $\{C_{96}(1,10,47)$, $C_{96}(5,43,46)$, $C_{96}(7,26,41)$, $C_{96}(11,14,37)$, 
 
 \hfill  $C_{96}(13,34,35)$, $C_{96}(17,22,31)$, $C_{96}(2,19,29)$, $C_{96}(23,25,38)\}$.  
 
 $\Rightarrow$ $C_{96}(10,23,25) \notin Ad_{96}(C_{96}(1,10,47))$.
 
 $\Rightarrow$ $C_{96}(1,10,47)$ and $C_{96}(10,23,25)$ are Type-2 isomorphic w.r.t. $m$ = 2.
\item [\rm (g3)]  $Ad_{96}(C_{96}(1,14,47))$ = $\{C_{96}(1,14,47)$, $C_{96}(5,26,43,)$, $C_{96}(2,7,41)$, $C_{96}(11,37,38)$, 
 
 \hfill  $C_{96}(10,13,35)$, $C_{96}(17,31,46)$, $C_{96}(19,22,29)$, $C_{96}(23,25,34)\}$.  
 
 $\Rightarrow$ $C_{96}(14,23,25) \notin Ad_{96}(C_{96}(1,14,47))$.
 
 $\Rightarrow$ $C_{96}(1,14,47)$ and $C_{96}(14,23,25)$ are Type-2 isomorphic w.r.t. $m$ = 2.

\item [\rm (g4)]  $Ad_{96}(C_{96}(1,22,47))$ = $\{C_{96}(1,22,47)$, $C_{96}(5,14,43,)$, $C_{96}(7,38,41)$, $C_{96}(11,37,46)$, 
 
 \hfill  $C_{96}(2,13,35)$, $C_{96}(10,17,31)$, $C_{96}(19,29,34)$, $C_{96}(23,25,26)\}$.  
 
 $\Rightarrow$ $C_{96}(22,23,25) \notin Ad_{96}(C_{96}(1,22,47))$.
 
 $\Rightarrow$ $C_{96}(1,22,47)$ and $C_{96}(22,23,25)$ are Type-2 isomorphic w.r.t. $m$ = 2.

\item [\rm (g5)]  $Ad_{96}(C_{96}(1,26,47))$ = $\{C_{96}(1,26,47)$, $C_{96}(5,34,43)$, $C_{96}(7,10,41)$, $C_{96}(2,11,37)$, 
 
 \hfill  $C_{96}(13,35,46)$, $C_{96}(17,31,38)$, $C_{96}(14,19,29)$, $C_{96}(22,23,25)\}$.  
 
 $\Rightarrow$ $C_{96}(23,25,26) \notin Ad_{96}(C_{96}(1,26,47))$.
 
 $\Rightarrow$ $C_{96}(1,26,47)$ and $C_{96}(23,25,26)$ are Type-2 isomorphic w.r.t. $m$ = 2.

\item [\rm (g6)]  $Ad_{96}(C_{96}(1,34,47))$ = $\{C_{96}(1,34,47)$, $C_{96}(5,22,43)$, $C_{96}(7,41,46)$, $C_{96}(10,11,37)$, 
 
 \hfill  $C_{96}(13,35,38)$, $C_{96}(2,17,31)$, $C_{96}(19,26,29)$, $C_{96}(14,23,25)\}$.  
 
 $\Rightarrow$ $C_{96}(23,25,34) \notin Ad_{96}(C_{96}(1,34,47))$.
 
 $\Rightarrow$ $C_{96}(1,34,47)$ and $C_{96}(23,25,34)$ are Type-2 isomorphic w.r.t. $m$ = 2.

\item [\rm (g7)]  $Ad_{96}(C_{96}(1,38,47))$ = $\{C_{96}(1,38,47)$, $C_{96}(2,5,43)$, $C_{96}(7,22,41)$, $C_{96}(11,34,37)$, 
 
 \hfill  $C_{96}(13,14,35)$, $C_{96}(17,26,31)$, $C_{96}(19,29,46)$, $C_{96}(10,23,25)\}$.  
 
 $\Rightarrow$ $C_{96}(23,25,38) \notin Ad_{96}(C_{96}(1,38,47))$.
 
 $\Rightarrow$ $C_{96}(1,38,47)$ and $C_{96}(23,25,38)$ are Type-2 isomorphic w.r.t. $m$ = 2.

\item [\rm (g8)]  $Ad_{96}(C_{96}(1,46,47))$ = $\{C_{96}(1,46,47)$, $C_{96}(5,38,43)$, $C_{96}(7,34,41)$, $C_{96}(11,26,37)$, 
 
 \hfill  $C_{96}(13,22,35)$, $C_{96}(14,17,31)$, $C_{96}(10,19,29)$, $C_{96}(2,23,25)\}$.  
 
 $\Rightarrow$ $C_{96}(23,25,46) \notin Ad_{96}(C_{96}(1,46,47))$.
 
 $\Rightarrow$ $C_{96}(1,46,47)$ and $C_{96}(23,25,46)$ are Type-2 isomorphic w.r.t. $m$ = 2.

\item [\rm (h)]  For  $s$ = 2, 10, 14, 22, 26, 34, 38, 46, using the definition of $\theta_{n,m,t}$, we get,  
 \\
 $\theta_{96,2,12}(C_{96}(3,s,45))$ = $\theta_{96,2,12}(C_{96}(3,s,45, 51,96-s,93))$ = $C_{96}(\theta_{96,2,12}(3,s,45, 51,96-s,93))$ 
 
 \hspace{2.7cm} = $C_{96}(27,s,69, 75,96-s,21)$ = $C_{48}(s,21,27)$. 
 
 $\Rightarrow$ $C_{96}(3,s,45)$ $\cong$ $C_{96}(s,21,27)$ for $s$ = 2, 10, 14, 22, 26, 34, 38, 46. Also, for $s$ = 2, 10, 14, 22, 26, 34, 38, 46, we get,
\\ 
$Ad_{96}(C_{96}(3,s,45)$ = $Ad_{96}(C_{96}(3,s,45, 51,96-s,93))$ 
  = $\{C_{96}(x(3,s,45, 51,96-s,93)): x = 1$, 
$5,7,11,13,17,19,23,25,29,31,35,37,41,43,47,49,53,55,59,61,65,67,71,77,79,83,85,89,91\}$. 
 
 From the above relation, we find $Ad_{96}(C_{96}(3,s,45))$ for $s$ = 2, 10, 14, 22, 26, 34, 38, 46 and will show that $C_{96}(s,21,27) \notin Ad_{96}(C_{96}(3,s,45))$.
 
\item [\rm (h1)]  $Ad_{96}(C_{96}(2,3,45))$ = $\{C_{96}(2,3,45)$, $C_{96}(10,15,33)$, $C_{96}(14,21,27)$, $C_{96}(15,22,33)$,  
 
 \hfill  $C_{96}(9,26,39)$, $C_{96}(3,34,45)$, $C_{96}(9,38,39)$, $C_{96}(21,27,46)\}$.  

 $\Rightarrow$ $C_{96}(2,21,27) \notin Ad_{96}(C_{96}(2,3,45))$.
 
 $\Rightarrow$ $C_{96}(2,3,45)$ and $C_{96}(2,21,27)$ are Type-2 isomorphic w.r.t. $m$ = 2.
 
\item [\rm (h2)]  $Ad_{96}(C_{96}(3,10,45))$ = $\{C_{96}(3,10,45)$, $C_{96}(15,33,46)$, $C_{96}(21,26,27)$, $C_{96}(14,15,33)$,  
 
 \hfill  $C_{96}(9,34,39)$, $C_{96}(3,22,45)$, $C_{96}(2,9,39)$, $C_{96}(21,27,38)\}$.  

 $\Rightarrow$ $C_{96}(10,21,27) \notin Ad_{96}(C_{96}(3,10,45))$.
 
 $\Rightarrow$ $C_{96}(3,10,45)$ and $C_{96}(10,21,27)$ are Type-2 isomorphic w.r.t. $m$ = 2.

\item [\rm (h3)]  $Ad_{96}(C_{96}(3,14,45))$ = $\{C_{96}(3,14,45)$, $C_{96}(15,26,33,)$, $C_{96}(2,21,27)$, $C_{96}(15,33,38)$,  
 
 \hfill  $C_{96}(9,10,39)$, $C_{96}(3,45,46)$, $C_{96}(9,22,39)$, $C_{96}(21,27,34)\}$.  

 $\Rightarrow$ $C_{96}(14,21,27) \notin Ad_{96}(C_{96}(3,14,45))$.
 
 $\Rightarrow$ $C_{96}(3,14,45)$ and $C_{96}(14,21,27)$ are Type-2 isomorphic w.r.t. $m$ = 2.

\item [\rm (h4)]  $Ad_{96}(C_{96}(3,22,45))$ = $\{C_{96}(3,22,45)$, $C_{96}(14,15,33)$, $C_{96}(21,27,38)$, $C_{96}(15,33,46)$,  
 
 \hfill  $C_{96}(2,9,39)$, $C_{96}(3,10,45)$, $C_{96}(9,34,39)$, $C_{96}(21,26,27)\}$.  

 $\Rightarrow$ $C_{96}(21,22,27) \notin Ad_{96}(C_{96}(3,22,45))$.
 
 $\Rightarrow$ $C_{96}(3,22,45)$ and $C_{96}(21,22,27)$ are Type-2 isomorphic w.r.t. $m$ = 2.

\item [\rm (h5)]  $Ad_{96}(C_{96}(3,26,45))$ = $\{C_{96}(3,26,45)$, $C_{96}(15,33,34)$, $C_{96}(10,21,27)$, $C_{96}(2,15,33)$,  
 
 \hfill  $C_{96}(9,39,46)$, $C_{96}(3,38,45)$, $C_{96}(9,14,39)$, $C_{96}(21,22,27)\}$.  

 $\Rightarrow$ $C_{96}(21,26,27) \notin Ad_{96}(C_{96}(3,26,45))$.
 
 $\Rightarrow$ $C_{96}(3,26,45)$ and $C_{96}(21,26,27)$ are Type-2 isomorphic w.r.t. $m$ = 2.

\item [\rm (h6)]  $Ad_{96}(C_{96}(3,34,45))$ = $\{C_{96}(3,34,45)$, $C_{96}(15,22,33)$, $C_{96}(21,27,46)$, $C_{96}(10,15,33)$,  
 
 \hfill  $C_{96}(9,38,39,)$, $C_{96}(2,3,45)$, $C_{96}(9,26,39)$, $C_{96}(14,21,27)\}$.  

 $\Rightarrow$ $C_{96}(21,27,34) \notin Ad_{96}(C_{96}(3,34,45))$.
 
 $\Rightarrow$ $C_{96}(3,34,45)$ and $C_{96}(21,27,34)$ are Type-2 isomorphic w.r.t. $m$ = 2.

\item [\rm (h7)]  $Ad_{96}(C_{96}(3,38,45))$ = $\{C_{96}(3,38,45)$, $C_{96}(2,15,33)$, $C_{96}(21,22,27)$, $C_{96}(15,33,34)$,  
 
 \hfill  $C_{96}(9,14,39,)$, $C_{96}(3,26,45)$, $C_{96}(9,39,46)$, $C_{96}(10,21,27)\}$.  

 $\Rightarrow$ $C_{96}(21,27,38) \notin Ad_{96}(C_{96}(3,38,45))$.
 
 $\Rightarrow$ $C_{96}(3,38,45)$ and $C_{96}(21,27,38)$ are Type-2 isomorphic w.r.t. $m$ = 2.

\item [\rm (h8)]  $Ad_{96}(C_{96}(3,45,46))$ = $\{C_{96}(3,45,46)$, $C_{96}(15,33,38)$, $C_{96}(21,27,34)$, $C_{96}(15,26,33)$,  
 
 \hfill  $C_{96}(9,22,39,)$, $C_{96}(3,14,45)$, $C_{96}(9,10,39)$, $C_{96}(2,21,27)\}$.  

 $\Rightarrow$ $C_{96}(21,27,46) \notin Ad_{96}(C_{96}(3,45,46))$.
 
 $\Rightarrow$ $C_{96}(3,45,46)$ and $C_{96}(21,27,46)$ are Type-2 isomorphic w.r.t. $m$ = 2.

\item [\rm (i)]  For  $s$ = 2, 10, 14, 22, 26, 34, 38, 46, using the definition of $\theta_{n,m,t}$, we get,  
 \\
 $\theta_{96,2,12}(C_{96}(5,s,43))$ = $\theta_{96,2,12}(C_{96}(5,s,43, 53,96-s,91))$ = $C_{96}(\theta_{96,2,12}(5,s,43, 53,96-s,91))$ 
 
 \hspace{2.7cm} = $C_{96}(29,s,67, 77,96-s,19)$ = $C_{48}(s,19,29)$. 
 
 $\Rightarrow$ $C_{96}(5,s,43)$ $\cong$ $C_{96}(s,19,29)$ for $s$ = 2, 10, 14, 22, 26, 34, 38, 46. Also, for $s$ = 2, 10, 14, 22, 26, 34, 38, 46, we get,
\\ 
$Ad_{96}(C_{96}(5,s,43)$ = $Ad_{96}(C_{96}(5,s,43, 53,96-s,91))$ 
  = $\{C_{96}(x(5,s,43, 53,96-s,91)): x = 1$, 
$5,7,11,13,17,19,23,25,29,31,35,37,41,43,47,49,53,55,59,61,65,67,71,77,79,83,85,89,91\}$. 
 
 From the above relation, we find $Ad_{96}(C_{96}(5,s,43))$ for $s$ = 2, 10, 14, 22, 26, 34, 38, 46 and will show that $C_{96}(s,19,29) \notin Ad_{96}(C_{96}(5,s,43))$.
 
\item [\rm (i1)]  $Ad_{96}(C_{96}(2,5,43))$ = $\{C_{96}(2,5,43)$, $C_{96}(10,23,25)$, $C_{96}(13,14,35)$, $C_{96}(7,22,41)$,  
 
 \hfill  $C_{96}(17,26,31)$, $C_{96}(11,34,37)$, $C_{96}(1,38,47)$, $C_{96}(19,29,46)\}$.  

 $\Rightarrow$ $C_{96}(2,19,29) \notin Ad_{96}(C_{96}(2,5,43))$.
 
 $\Rightarrow$ $C_{96}(2,5,43)$ and $C_{96}(2,19,29)$ are Type-2 isomorphic w.r.t. $m$ = 2.
 
\item [\rm (i2)]  $Ad_{96}(C_{96}(5,10,43))$ = $\{C_{96}(5,10,43)$, $C_{96}(23,25,46)$, $C_{96}(13,26,35)$, $C_{96}(7,14,41)$,  
 
 \hfill  $C_{96}(17,31,34)$, $C_{96}(11,22,37)$, $C_{96}(1,2,47)$, $C_{96}(19,29,38)\}$.  

 $\Rightarrow$ $C_{96}(10,19,29) \notin Ad_{96}(C_{96}(5,10,43))$.
 
 $\Rightarrow$ $C_{96}(5,10,43)$ and $C_{96}(10,19,29)$ are Type-2 isomorphic w.r.t. $m$ = 2.
 
\item [\rm (i3)]  $Ad_{96}(C_{96}(5,14,43))$ = $\{C_{96}(5,14,43)$, $C_{96}(23,25,26)$, $C_{96}(2,13,35)$, $C_{96}(7,38,41)$,  
 
 \hfill  $C_{96}(10,17,31)$, $C_{96}(11,37,46)$, $C_{96}(1,22,47)$, $C_{96}(19,29,34)\}$.  

 $\Rightarrow$ $C_{96}(14,19,29) \notin Ad_{96}(C_{96}(5,14,43))$.
 
 $\Rightarrow$ $C_{96}(5,14,43)$ and $C_{96}(14,19,29)$ are Type-2 isomorphic w.r.t. $m$ = 2.
 
\item [\rm (i4)]  $Ad_{96}(C_{96}(5,22,43))$ = $\{C_{96}(5,22,43)$, $C_{96}(14,23,25)$, $C_{96}(13,35,38)$, $C_{96}(7,41,46)$,  
 
 \hfill  $C_{96}(2,17,31)$, $C_{96}(10,11,37)$, $C_{96}(1,34,47)$, $C_{96}(19,26,29)\}$.  

 $\Rightarrow$ $C_{96}(19,22,29) \notin Ad_{96}(C_{96}(5,22,43))$.
 
 $\Rightarrow$ $C_{96}(5,22,43)$ and $C_{96}(19,22,29)$ are Type-2 isomorphic w.r.t. $m$ = 2.
 
\item [\rm (i5)]  $Ad_{96}(C_{96}(5,26,43))$ = $\{C_{96}(5,26,43)$, $C_{96}(23,25,34)$, $C_{96}(10,13,35)$, $C_{96}(2,7,41)$,  
 
 \hfill  $C_{96}(17,31,46)$, $C_{96}(11,37,38)$, $C_{96}(1,14,47)$, $C_{96}(19,22,29)\}$.  

 $\Rightarrow$ $C_{96}(19,26,29) \notin Ad_{96}(C_{96}(5,26,43))$.
 
 $\Rightarrow$ $C_{96}(5,26,43)$ and $C_{96}(19,26,29)$ are Type-2 isomorphic w.r.t. $m$ = 2.
 
\item [\rm (i6)]  $Ad_{96}(C_{96}(5,34,43))$ = $\{C_{96}(5,34,43)$, $C_{96}(22,23,25)$, $C_{96}(13,35,46)$, $C_{96}(7,10,41)$,  
 
 \hfill  $C_{96}(17,31,38)$, $C_{96}(2,11,37)$, $C_{96}(1,26,47)$, $C_{96}(14,19,29)\}$.  

 $\Rightarrow$ $C_{96}(19,29,34) \notin Ad_{96}(C_{96}(5,34,43))$.
 
 $\Rightarrow$ $C_{96}(5,34,43)$ and $C_{96}(19,29,34)$ are Type-2 isomorphic w.r.t. $m$ = 2.
 
\item [\rm (i7)]  $Ad_{96}(C_{96}(5,38,43))$ = $\{C_{96}(5,38,43)$, $C_{96}(2,23,25)$, $C_{96}(13,22,35)$, $C_{96}(7,34,41)$,  
 
 \hfill  $C_{96}(14,17,31)$, $C_{96}(11,26,37)$, $C_{96}(1,46,47)$, $C_{96}(10,19,29)\}$.  

 $\Rightarrow$ $C_{96}(19,29,38) \notin Ad_{96}(C_{96}(5,38,43))$.
 
 $\Rightarrow$ $C_{96}(5,38,43)$ and $C_{96}(19,29,38)$ are Type-2 isomorphic w.r.t. $m$ = 2.
 
\item [\rm (i8)]  $Ad_{96}(C_{96}(5,43,46))$ = $\{C_{96}(5,43,46)$, $C_{96}(,23,25)$, $C_{96}(13,,35)$, $C_{96}(7,,41)$,  
 
 \hfill  $C_{96}(,17,31)$, $C_{96}(11,,37)$, $C_{96}(1,,47)$, $C_{96}(,19,29)\}$.  

 $\Rightarrow$ $C_{96}(19,29,46) \notin Ad_{96}(C_{96}(5,43,46))$.
 
 $\Rightarrow$ $C_{96}(5,43,46)$ and $C_{96}(19,29,46)$ are Type-2 isomorphic w.r.t. $m$ = 2.
 
\item [\rm (j)]  For  $s$ = 2, 10, 14, 22, 26, 34, 38, 46, using the definition of $\theta_{n,m,t}$, we get,  
 \\
 $\theta_{96,2,12}(C_{96}(7,s,41))$ = $\theta_{96,2,12}(C_{96}(7,s,41, 55,96-s,89))$ = $C_{96}(\theta_{96,2,12}(7,s,41, 55,96-s,89))$ 
 
 \hspace{2.7cm} = $C_{96}(31,s,65, 79,96-s,17)$ = $C_{48}(s,17,31)$. 
 
 $\Rightarrow$ $C_{96}(7,s,41)$ $\cong$ $C_{96}(s,17,31)$ for $s$ = 2, 10, 14, 22, 26, 34, 38, 46. Also, for $s$ = 2, 10, 14, 22, 26, 34, 38, 46, we get,
\\ 
$Ad_{96}(C_{96}(7,s,41)$ = $Ad_{96}(C_{96}(7,s,41, 55,96-s,89))$ 
  = $\{C_{96}(x(7,s,41, 55,96-s,89)): x = 1$, 
$5,7,11,13,17,19,23,25,29,31,35,37,41,43,47,49,53,55,59,61,65,67,71,77,79,83,85,89,91\}$. 
 
 From the above relation, we find $Ad_{96}(C_{96}(7,s,41))$ for $s$ = 2, 10, 14, 22, 26, 34, 38, 46 and will show that $C_{96}(s,17,31) \notin Ad_{96}(C_{96}(7,s,41))$.
 
\item [\rm (j1)]  $Ad_{96}(C_{96}(2,7,41))$ = $\{C_{96}(2,7,41)$, $C_{96}(10,13,35)$, $C_{96}(1,14,47)$, $C_{96}(19,22,29)$,  
 
 \hfill  $C_{96}(5,26,43)$, $C_{96}(23,25,34)$, $C_{96}(11,37,38)$, $C_{96}(17,31,46)\}$.  

 $\Rightarrow$ $C_{96}(2,17,31) \notin Ad_{96}(C_{96}(2,7,41))$.
 
 $\Rightarrow$ $C_{96}(2,7,41)$ and $C_{96}(2,17,31)$ are Type-2 isomorphic w.r.t. $m$ = 2.

\item [\rm (j2)]  $Ad_{96}(C_{96}(7,10,41))$ = $\{C_{96}(7,10,41)$, $C_{96}(13,35,46)$, $C_{96}(1,26,47)$, $C_{96}(14,19,29)$,  
 
 \hfill  $C_{96}(5,34,43)$, $C_{96}(22,23,25)$, $C_{96}(2,11,37)$, $C_{96}(17,31,38)\}$.  

 $\Rightarrow$ $C_{96}(10,17,31) \notin Ad_{96}(C_{96}(7,10,41))$.
 
 $\Rightarrow$ $C_{96}(7,10,41)$ and $C_{96}(10,17,31)$ are Type-2 isomorphic w.r.t. $m$ = 2.

\item [\rm (j3)]  $Ad_{96}(C_{96}(7,14,41))$ = $\{C_{96}(7,14,41)$, $C_{96}(13,26,35)$, $C_{96}(1,2,47)$, $C_{96}(19,29,38)$,  
 
 \hfill  $C_{96}(5,10,43)$, $C_{96}(23,25,46)$, $C_{96}(11,22,37)$, $C_{96}(17,31,34)\}$.  

 $\Rightarrow$ $C_{96}(14,17,31) \notin Ad_{96}(C_{96}(7,14,41))$.
 
 $\Rightarrow$ $C_{96}(7,14,41)$ and $C_{96}(14,17,31)$ are Type-2 isomorphic w.r.t. $m$ = 2.

\item [\rm (j4)]  $Ad_{96}(C_{96}(7,22,41))$ = $\{C_{96}(7,22,41)$, $C_{96}(13,14,35)$, $C_{96}(1,38,47)$, $C_{96}(19,29,46)$,  
 
 \hfill  $C_{96}(2,5,43)$, $C_{96}(10,23,25)$, $C_{96}(11,34,37)$, $C_{96}(17,26,31)\}$.  

 $\Rightarrow$ $C_{96}(17,22,31) \notin Ad_{96}(C_{96}(7,22,41))$.
 
 $\Rightarrow$ $C_{96}(7,22,41)$ and $C_{96}(17,22,31)$ are Type-2 isomorphic w.r.t. $m$ = 2.

\item [\rm (j5)]  $Ad_{96}(C_{96}(7,26,41))$ = $\{C_{96}(7,26,41)$, $C_{96}(13,34,35)$, $C_{96}(1,10,47)$, $C_{96}(2,19,29)$,  
 
 \hfill  $C_{96}(5,43,46)$, $C_{96}(23,25,38)$, $C_{96}(11,14,37)$, $C_{96}(17,22,31)\}$.  

 $\Rightarrow$ $C_{96}(17,26,31) \notin Ad_{96}(C_{96}(7,26,41))$.
 
 $\Rightarrow$ $C_{96}(7,26,41)$ and $C_{96}(17,26,31)$ are Type-2 isomorphic w.r.t. $m$ = 2.

\item [\rm (j6)]  $Ad_{96}(C_{96}(7,34,41))$ = $\{C_{96}(7,34,41)$, $C_{96}(13,22,35)$, $C_{96}(1,46,47)$, $C_{96}(10,19,29)$,  
 
 \hfill  $C_{96}(5,38,43)$, $C_{96}(2,23,25)$, $C_{96}(11,26,37)$, $C_{96}(14,17,31)\}$.  

 $\Rightarrow$ $C_{96}(17,31,34) \notin Ad_{96}(C_{96}(7,34,41))$.
 
 $\Rightarrow$ $C_{96}(7,34,41)$ and $C_{96}(17,31,34)$ are Type-2 isomorphic w.r.t. $m$ = 2.

\item [\rm (j7)]  $Ad_{96}(C_{96}(7,38,41))$ = $\{C_{96}(7,38,41)$, $C_{96}(2,13,35)$, $C_{96}(1,22,47)$, $C_{96}(19,29,34)$,  
 
 \hfill  $C_{96}(5,14,43)$, $C_{96}(23,25,26)$, $C_{96}(11,37,46)$, $C_{96}(10,17,31)\}$.  

 $\Rightarrow$ $C_{96}(17,31,38) \notin Ad_{96}(C_{96}(7,38,41))$.
 
 $\Rightarrow$ $C_{96}(7,38,41)$ and $C_{96}(17,31,38)$ are Type-2 isomorphic w.r.t. $m$ = 2.

\item [\rm (j8)]  $Ad_{96}(C_{96}(7,41,46))$ = $\{C_{96}(7,41,46)$, $C_{96}(13,35,38)$, $C_{96}(1,34,47)$, $C_{96}(19,26,29)$,  
 
 \hfill  $C_{96}(5,22,43)$, $C_{96}(14,23,25)$, $C_{96}(10,11,37)$, $C_{96}(2,17,31)\}$.  

 $\Rightarrow$ $C_{96}(17,31,46) \notin Ad_{96}(C_{96}(7,41,46))$.
 
 $\Rightarrow$ $C_{96}(7,41,46)$ and $C_{96}(17,31,46)$ are Type-2 isomorphic w.r.t. $m$ = 2.

\item [\rm (k)]  For  $s$ = 2, 10, 14, 22, 26, 34, 38, 46, using the definition of $\theta_{n,m,t}$, we get,  
 \\
 $\theta_{96,2,12}(C_{96}(9,s,39))$ = $\theta_{96,2,12}(C_{96}(9,s,39, 57,96-s,87))$ = $C_{96}(\theta_{96,2,12}(9,s,39, 57,96-s,87))$ 
 
 \hspace{2.7cm} = $C_{96}(33,s,63, 81,96-s,15)$ = $C_{48}(s,15,33)$. 
 
 $\Rightarrow$ $C_{96}(9,s,39)$ $\cong$ $C_{96}(s,15,33)$ for $s$ = 2, 10, 14, 22, 26, 34, 38, 46. Also, for $s$ = 2, 10, 14, 22, 26, 34, 38, 46, we get,
\\ 
$Ad_{96}(C_{96}(9,s,39)$ = $Ad_{96}(C_{96}(9,s,39, 57,96-s,87))$ 
  = $\{C_{96}(x(9,s,39, 57,96-s,87)): x = 1$, 
$5,7,11,13,17,19,23,25,29,31,35,37,41,43,47,49,53,55,59,61,65,67,71,77,79,83,85,89,91\}$. 
 
 From the above relation, we find $Ad_{96}(C_{96}(9,s,39))$ for $s$ = 2, 10, 14, 22, 26, 34, 38, 46 and will show that $C_{96}(s,15,33) \notin Ad_{96}(C_{96}(9,s,39))$.
 
\item [\rm (k1)]  $Ad_{96}(C_{96}(2,9,39))$ = $\{C_{96}(2,9,39)$, $C_{96}(3,10,45)$, $C_{96}(14,15,33)$, $C_{96}(3,22,45)$,  
 
 \hfill  $C_{96}(21,26,27)$, $C_{96}(9,34,39)$, $C_{96}(21,27,38)$, $C_{96}(15,33,46)\}$.  

 $\Rightarrow$ $C_{96}(2,15,33) \notin Ad_{96}(C_{96}(2,9,39))$.
 
 $\Rightarrow$ $C_{96}(2,9,39)$ and $C_{96}(2,15,33)$ are Type-2 isomorphic w.r.t. $m$ = 2.

\item [\rm (k2)]  $Ad_{96}(C_{96}(9,10,39))$ = $\{C_{96}(9,10,39)$, $C_{96}(3,45,46)$, $C_{96}(15,26,33)$, $C_{96}(3,14,45)$,  
 
 \hfill  $C_{96}(21,27,34)$, $C_{96}(9,22,39)$, $C_{96}(2,21,27)$, $C_{96}(15,33,38)\}$.  

 $\Rightarrow$ $C_{96}(10,15,33) \notin Ad_{96}(C_{96}(9,10,39))$.
 
 $\Rightarrow$ $C_{96}(9,10,39)$ and $C_{96}(10,15,33)$ are Type-2 isomorphic w.r.t. $m$ = 2.

\item [\rm (k3)]  $Ad_{96}(C_{96}(9,14,39))$ = $\{C_{96}(9,14,39)$, $C_{96}(3,26,45)$, $C_{96}(2,15,33)$, $C_{96}(3,38,45)$,  
 
 \hfill  $C_{96}(10,21,27)$, $C_{96}(9,39,46)$, $C_{96}(21,22,27)$, $C_{96}(15,33,34)\}$.  

 $\Rightarrow$ $C_{96}(14,15,33) \notin Ad_{96}(C_{96}(9,14,39))$.
 
 $\Rightarrow$ $C_{96}(9,14,39)$ and $C_{96}(14,15,33)$ are Type-2 isomorphic w.r.t. $m$ = 2.

\item [\rm (k4)]  $Ad_{96}(C_{96}(9,22,39))$ = $\{C_{96}(9,22,39)$, $C_{96}(3,14,45)$, $C_{96}(15,33,38)$, $C_{96}(3,45,46)$,  
 
 \hfill  $C_{96}(2,21,27)$, $C_{96}(9,10,39)$, $C_{96}(21,27,34)$, $C_{96}(15,26,33)\}$.  

 $\Rightarrow$ $C_{96}(15,22,33) \notin Ad_{96}(C_{96}(9,22,39))$.
 
 $\Rightarrow$ $C_{96}(9,22,39)$ and $C_{96}(15,22,33)$ are Type-2 isomorphic w.r.t. $m$ = 2.

\item [\rm (k5)]  $Ad_{96}(C_{96}(9,26,39))$ = $\{C_{96}(9,26,39)$, $C_{96}(3,34,45)$, $C_{96}(10,15,33)$, $C_{96}(2,3,45)$,  
 
 \hfill  $C_{96}(21,27,46)$, $C_{96}(9,38,39)$, $C_{96}(14,21,27)$, $C_{96}(15,22,33)\}$.  

 $\Rightarrow$ $C_{96}(15,26,33) \notin Ad_{96}(C_{96}(9,26,39))$.
 
 $\Rightarrow$ $C_{96}(9,26,39)$ and $C_{96}(15,26,33)$ are Type-2 isomorphic w.r.t. $m$ = 2.

\item [\rm (k6)]  $Ad_{96}(C_{96}(9,34,39))$ = $\{C_{96}(9,34,39)$, $C_{96}(3,22,45)$, $C_{96}(15,33,46)$, $C_{96}(3,10,45)$,  
 
 \hfill  $C_{96}(21,27,38)$, $C_{96}(2,9,39)$, $C_{96}(21,26,27)$, $C_{96}(14,15,33)\}$.  

 $\Rightarrow$ $C_{96}(15,33,34) \notin Ad_{96}(C_{96}(9,34,39))$.
 
 $\Rightarrow$ $C_{96}(9,34,39)$ and $C_{96}(15,33,34)$ are Type-2 isomorphic w.r.t. $m$ = 2.

\item [\rm (k7)]  $Ad_{96}(C_{96}(9,38,39))$ = $\{C_{96}(9,38,39)$, $C_{96}(2,3,45)$, $C_{96}(15,22,33)$, $C_{96}(3,34,45)$,  
 
 \hfill  $C_{96}(14,21,27)$, $C_{96}(9,26,39)$, $C_{96}(21,27,46)$, $C_{96}(10,15,33)\}$.  

 $\Rightarrow$ $C_{96}(15,33,38) \notin Ad_{96}(C_{96}(9,38,39))$.
 
 $\Rightarrow$ $C_{96}(9,38,39)$ and $C_{96}(15,33,38)$ are Type-2 isomorphic w.r.t. $m$ = 2.

\item [\rm (k8)]  $Ad_{96}(C_{96}(9,39,46))$ = $\{C_{96}(9,39,46)$, $C_{96}(3,38,45)$, $C_{96}(15,33,34)$, $C_{96}(3,26,45)$,  
 
 \hfill  $C_{96}(21,22,27)$, $C_{96}(9,14,39)$, $C_{96}(10,21,27)$, $C_{96}(2,15,33)\}$.  

 $\Rightarrow$ $C_{96}(15,33,46) \notin Ad_{96}(C_{96}(9,39,46))$.
 
 $\Rightarrow$ $C_{96}(9,39,46)$ and $C_{96}(15,33,46)$ are Type-2 isomorphic w.r.t. $m$ = 2.

\item [\rm (l)]  For  $s$ = 2, 10, 14, 22, 26, 34, 38, 46, using the definition of $\theta_{n,m,t}$, we get,  
 \\
 $\theta_{96,2,12}(C_{96}(11,s,37))$ = $\theta_{96,2,12}(C_{96}(11,s,37, 59,96-s,85))$ 

 \hspace{3cm} = $C_{96}(\theta_{96,2,12}(11,s,37, 59,96-s,85))$ 
 
 \hspace{3cm} = $C_{96}(35,s,61, 83,96-s,13)$ = $C_{48}(s,13,35)$. 
 
 $\Rightarrow$ $C_{96}(11,s,37)$ $\cong$ $C_{96}(s,13,35)$ for $s$ = 2, 10, 14, 22, 26, 34, 38, 46. Also, for $s$ = 2, 10, 14, 22, 26, 34, 38, 46, we get,
\\ 
$Ad_{96}(C_{96}(11,s,37)$ = $Ad_{96}(C_{96}(11,s,37, 59,96-s,85))$ 
  = $\{C_{96}(x(11,s,37, 59,96-s,85)): x = 1$, 
$5,7,11,13,17,19,23,25,29,31,35,37,41,43,47,49,53,55,59,61,65,67,71,77,79,83,85,89,91\}$. 
 
 From the above relation, we find $Ad_{96}(C_{96}(11,s,37))$ for $s$ = 2, 10, 14, 22, 26, 34, 38, 46 and will show that $C_{96}(s,13,35) \notin Ad_{96}(C_{96}(11,s,37))$.
 
\item [\rm (l1)]  $Ad_{96}(C_{96}(2,11,37))$ = $\{C_{96}(2,11,37)$, $C_{96}(7,10,41)$, $C_{96}(14,19,29)$, $C_{96}(22,23,25)$,  
 
 \hfill  $C_{96}(1,26,47)$, $C_{96}(5,34,43)$, $C_{96}(17,31,38)$, $C_{96}(13,35,46)\}$.  

 $\Rightarrow$ $C_{96}(2,13,35) \notin Ad_{96}(C_{96}(2,11,37))$.
 
 $\Rightarrow$ $C_{96}(2,11,37)$ and $C_{96}(2,13,35)$ are Type-2 isomorphic w.r.t. $m$ = 2.

\item [\rm (l2)]  $Ad_{96}(C_{96}(10,11,37))$ = $\{C_{96}(10,11,37)$, $C_{96}(7,41,46)$, $C_{96}(19,26,29)$, $C_{96}(14,23,25)$,  
 
 \hfill  $C_{96}(1,34,47)$, $C_{96}(5,22,43)$, $C_{96}(2,17,31)$, $C_{96}(13,35,38)\}$.  

 $\Rightarrow$ $C_{96}(10,13,35) \notin Ad_{96}(C_{96}(10,11,37))$.
 
 $\Rightarrow$ $C_{96}(10,11,37)$ and $C_{96}(10,13,35)$ are Type-2 isomorphic w.r.t. $m$ = 2.

\item [\rm (l3)]  $Ad_{96}(C_{96}(11,14,37))$ = $\{C_{96}(11,14,37)$, $C_{96}(7,26,41)$, $C_{96}(2,19,29)$, $C_{96}(23,25,38)$,  
 
 \hfill  $C_{96}(1,10,47)$, $C_{96}(5,43,46)$, $C_{96}(17,22,31)$, $C_{96}(13,34,35)\}$.  

 $\Rightarrow$ $C_{96}(13,14,35) \notin Ad_{96}(C_{96}(11,14,37))$.
 
 $\Rightarrow$ $C_{96}(11,14,37)$ and $C_{96}(13,14,35)$ are Type-2 isomorphic w.r.t. $m$ = 2.

\item [\rm (l4)]  $Ad_{96}(C_{96}(11,22,37))$ = $\{C_{96}(11,22,37)$, $C_{96}(7,14,41)$, $C_{96}(19,29,38)$, $C_{96}(23,25,46)$,  
 
 \hfill  $C_{96}(1,2,47)$, $C_{96}(5,10,43)$, $C_{96}(17,31,34)$, $C_{96}(13,26,35)\}$.  

 $\Rightarrow$ $C_{96}(13,22,35) \notin Ad_{96}(C_{96}(11,22,37))$.
 
 $\Rightarrow$ $C_{96}(11,22,37)$ and $C_{96}(13,22,35)$ are Type-2 isomorphic w.r.t. $m$ = 2.

\item [\rm (l5)]  $Ad_{96}(C_{96}(11,26,37))$ = $\{C_{96}(11,26,37)$, $C_{96}(7,34,41)$, $C_{96}(10,19,29)$, $C_{96}(2,23,25)$,  
 
 \hfill  $C_{96}(1,46,47)$, $C_{96}(5,38,43)$, $C_{96}(14,17,31)$, $C_{96}(13,22,35)\}$.  

 $\Rightarrow$ $C_{96}(13,26,35) \notin Ad_{96}(C_{96}(11,26,37))$.
 
 $\Rightarrow$ $C_{96}(11,26,37)$ and $C_{96}(13,26,35)$ are Type-2 isomorphic w.r.t. $m$ = 2.

\item [\rm (l6)]  $Ad_{96}(C_{96}(11,34,37))$ = $\{C_{96}(11,34,37)$, $C_{96}(7,22,41)$, $C_{96}(19,29,46)$, $C_{96}(10,23,25)$,  
 
 \hfill  $C_{96}(1,38,47)$, $C_{96}(2,5,43)$, $C_{96}(17,26,31)$, $C_{96}(13,14,35)\}$.  

 $\Rightarrow$ $C_{96}(13,34,35) \notin Ad_{96}(C_{96}(11,34,37))$.
 
 $\Rightarrow$ $C_{96}(11,34,37)$ and $C_{96}(13,34,35)$ are Type-2 isomorphic w.r.t. $m$ = 2.

\item [\rm (l7)]  $Ad_{96}(C_{96}(11,37,38))$ = $\{C_{96}(11,37,38)$, $C_{96}(2,7,41)$, $C_{96}(19,22,29)$, $C_{96}(23,25,34)$,  
 
 \hfill  $C_{96}(1,14,47)$, $C_{96}(5,26,43)$, $C_{96}(17,31,46)$, $C_{96}(10,13,35)\}$.  

 $\Rightarrow$ $C_{96}(13,35,38) \notin Ad_{96}(C_{96}(11,37,38))$.
 
 $\Rightarrow$ $C_{96}(11,37,38)$ and $C_{96}(13,35,38)$ are Type-2 isomorphic w.r.t. $m$ = 2.

\item [\rm (l8)]  $Ad_{96}(C_{96}(11,37,46))$ = $\{C_{96}(11,37,46)$, $C_{96}(7,38,41)$, $C_{96}(19,29,34)$, $C_{96}(23,25,26)$,  
 
 \hfill  $C_{96}(1,22,47)$, $C_{96}(5,14,43)$, $C_{96}(10,17,31)$, $C_{96}(2,13,35)\}$.  

 $\Rightarrow$ $C_{96}(13,35,46) \notin Ad_{96}(C_{96}(11,37,46))$.
 
 $\Rightarrow$ $C_{96}(11,37,46)$ and $C_{96}(13,35,46)$ are Type-2 isomorphic w.r.t. $m$ = 2. \hfill $\Box$
\end{enumerate}
 
 \begin{prm} \label{p9} {\rm For $s$ = 6,18,30,42,  show that the following pairs of circulant graphs are Type-2 isomorphic w.r.t. $m$ = 2.  
 \begin{enumerate}
 	\item [\rm (u)]  $C_{96}(1,s,47)$ and $C_{48}(s,23,25)$;  
 	\item [\rm (v)]  $C_{96}(3,s,45)$ and $C_{48}(s,21,27)$;  
 	\item [\rm (w)]  $C_{96}(5,s,43)$ and $C_{48}(s,19,29)$;  
 	\item [\rm (x)]  $C_{96}(7,s,41)$ and $C_{48}(s,17,31)$;  
 	\item [\rm (y)]  $C_{96}(9,s,39)$ and $C_{48}(s,15,33)$; and  
 	\item [\rm (z)]  $C_{96}(11,s,37)$ and $C_{48}(s,13,35)$ 
where $\gcd(96, s)$ = 6.  
 \end{enumerate} }	
 \end{prm}
\noindent
{\bf Solution.}\quad For $s$ = 6, 18, 30, 42, let

$R_1$ = $\{1,s,47\}$, $R_2$ = $\{3,s,45\}$, $R_3$ = $\{5,s,43\}$, $R_4$ = $\{7,s,41\}$, $R_5$ = $\{9,s,39\}$, $R_6$ = $\{11,s,37\}$,  

$S_1$ = $\{s,23,25\}$, $S_2$ = $\{s,21,27\}$, $S_3$ = $\{s,19,29\}$, $S_4$ = $\{s,17,31\}$, $S_5$ = $\{s,15,33\}$, $S_6$ = $\{s,13,35\}$. Here, $n$ = 96 = $2^2\times 3\times 2^3$ and $r$ = $s$ so that $s\in R_i,S_i$ and $\gcd(n, s)$ = $\gcd(96, 6)$ = 6 = $\gcd(96,s )$ for $s$ = 6, 18, 30, 42 and $i$ = 1 to 6. Using the definition of $\theta_{n,m,t}$, we get, $\theta_{96,2,12}(s)$ = $s$ and $\theta_{96,2,12}(96-s)$ = $96-s$. Now, consider different cases.
  \begin{enumerate}
 	\item [\rm (u)]  For  $s$ = 6, 18, 30, 42, using the definition of $\theta_{n,m,t}$, we get,  
 \\
 $\theta_{96,2,12}(C_{96}(1,s,47))$ = $\theta_{96,2,12}(C_{96}(1,s,47, 49,96-s,95))$ = $C_{96}(\theta_{96,2,12}(1,s,47, 49,96-s,95))$ 
 
 \hspace{2.7cm} = $C_{96}(25,s,71, 73,96-s,23)$ = $C_{48}(s,23,25)$. 
 
 $\Rightarrow$ $C_{96}(1,s,47)$ $\cong$ $C_{96}(s,23,25)$ for $s$ = 6, 18, 30, 42. Also, for $s$ = 6, 18, 30, 42, we get,
\\ 
$Ad_{96}(C_{96}(1,s,47)$ = $Ad_{96}(C_{96}(1,s,47, 49,96-s,95))$ 
  = $\{C_{96}(x(1,s,47, 49,96-s,95)): x = 1$, 
$5,7,11,13,17,19,23,25,29,31,35,37,41,43,47,49,53,55,59,61,65,67,71,77,79,83,85,89,91\}$. 
 
 From the above relation, we find $Ad_{96}(C_{96}(1,s,47))$ for $s$ = 6, 18, 30, 42 and will show that $C_{96}(s,23,25) \notin Ad_{96}(C_{96}(1,s,47))$.
 
\item [\rm (u1)]  $Ad_{96}(C_{96}(1,6,47))$ = $\{C_{96}(1,6,47)$, $C_{96}(5,30,43)$, $C_{96}(7,41,42)$, $C_{96}(11,30,37)$, 
 
 \hfill  $C_{96}(13,18,35)$, $C_{96}(6,17,31)$, $C_{96}(18,19,29)$, $C_{96}(23,25,42)\}$.  

 $\Rightarrow$ $C_{96}(6,23,25) \notin Ad_{96}(C_{96}(1,6,47))$.
 
 $\Rightarrow$ $C_{96}(1,6,47)$ and $C_{96}(6,23,25)$ are Type-2 isomorphic w.r.t. $m$ = 2.
 
\item [\rm (u2)]  $Ad_{96}(C_{96}(1,18,47))$ = $\{C_{96}(1,18,47)$, $C_{96}(5,6,43)$, $C_{96}(7,30,41)$, $C_{96}(6,11,37)$, 
 
 \hfill  $C_{96}(13,35,42)$, $C_{96}(17,18,31)$, $C_{96}(19,29,42)$, $C_{96}(23,25,30)\}$.  

 $\Rightarrow$ $C_{96}(18,23,25) \notin Ad_{96}(C_{96}(1,18,47))$.
 
 $\Rightarrow$ $C_{96}(1,18,47)$ and $C_{96}(18,23,25)$ are Type-2 isomorphic w.r.t. $m$ = 2.
 
\item [\rm (u3)]  $Ad_{96}(C_{96}(1,30,47))$ = $\{C_{96}(1,30,47)$, $C_{96}(5,42,43)$, $C_{96}(7,18,41)$, $C_{96}(11,37,42)$, 
 
 \hfill  $C_{96}(6,13,35)$, $C_{96}(17,30,31)$, $C_{96}(6,19,29)$, $C_{96}(18,23,25)\}$.  

 $\Rightarrow$ $C_{96}(23,25,30) \notin Ad_{96}(C_{96}(1,30,47))$.
 
 $\Rightarrow$ $C_{96}(1,30,47)$ and $C_{96}(23,25,30)$ are Type-2 isomorphic w.r.t. $m$ = 2.
 
\item [\rm (u4)]  $Ad_{96}(C_{96}(1,42,47))$ = $\{C_{96}(1,42,47)$, $C_{96}(5,18,43)$, $C_{96}(6,7,41)$, $C_{96}(11,18,37)$, 
 
 \hfill  $C_{96}(13,30,35)$, $C_{96}(17,31,42)$, $C_{96}(19,29,30)$, $C_{96}(6,23,25)\}$.  

 $\Rightarrow$ $C_{96}(23,25,42) \notin Ad_{96}(C_{96}(1,42,47))$.
 
 $\Rightarrow$ $C_{96}(1,42,47)$ and $C_{96}(23,25,42)$ are Type-2 isomorphic w.r.t. $m$ = 2.

\item [\rm (v)] When $s$ = 6,18,30,42, $C_{96}(3,s,45)$ = $3.C_{32}(1,t,15)$ and $C_{96}(s,21,27)$ = $3.C_{32}(t,7,9)$ where $s$ = $3t$. For $t$ = 2,6,10,14, $C_{32}(1,t,15)$ and $C_{32}(t,7,9)$ are Type-2 isomorphic w.r.t. $m$ = 2 by the following. $t$ = 2,6,10,14,

$\theta_{32,2,4}(C_{32}(1,t,15))$ = $\theta_{32,2,4}(C_{32}(1,t,15, 17,32-t,31))$ = $C_{32}(\theta_{32,2,4}(1,t,15, 17,32-t,31))$ 
 
 \hspace{2.7cm} = $C_{32}(9,t,23, 25,32-t,7)$ = $C_{32}(t,7,9)$.

$\Rightarrow$  $C_{32}(1,t,15)$ $\cong$ $C_{32}(t,7,9)$ for $t$ = 2,6,10,14. Also,
 
 $Ad_{32}(C_{32}(1,2,15))$ = $\{C_{32}(1,2,15)$, $C_{32}(3,6,13)$, $C_{32}(5,10,11)$, $C_{32}(7,9,14)\}$;

 $Ad_{32}(C_{32}(1,6,15))$ = $\{C_{32}(1,6,15)$, $C_{32}(3,13,14)$, $C_{32}(2,5,11)$, $C_{32}(7,9,10)\}$;

 $Ad_{32}(C_{32}(1,10,15))$ = $\{C_{32}(1,10,15)$, $C_{32}(2,3,13)$, $C_{32}(5,11,14)$, $C_{32}(6,7,9)\}$; and

 $Ad_{32}(C_{32}(1,14,15))$ = $\{C_{32}(1,14,15)$, $C_{32}(3,10,13)$, $C_{32}(5,6,11)$, $C_{32}(2,7,9)\}$.

$\Rightarrow$  $C_{32}(t,7,9)\notin Ad_{32}(C_{32}(1,t,15))$ for $t$ = 2,6,10,14. 

$\Rightarrow$  $C_{32}(1,t,15)$ and $C_{32}(t,7,9)$ are Type-2 isomorphic w.r.t. $m$ = 2  for $t$ = 2,6,10,14 where $\gcd(32, t)$ = 2 = $m$. 

This implies, $3.C_{32}(1,t,15)$ = $C_{96}(3,3t,45)$ and $3.C_{32}(t,7,9)$ = $C_{96}(3t,21,27)$ are Type-2 isomorphic w.r.t. $m$ = 2 for $t$ = 2,6,10,14 where $\gcd(32, t)$ = 2 = $m$. 

\item [\rm (w)]  For  $s$ = 6, 18, 30, 42, using the definition of $\theta_{n,m,t}$, we get,  
 \\
 $\theta_{96,2,12}(C_{96}(5,s,43))$ = $\theta_{96,2,12}(C_{96}(5,s,43, 53,96-s,91))$ = $C_{96}(\theta_{96,2,12}(5,s,43, 53,96-s,91))$ 
 
 \hspace{2.7cm} = $C_{96}(29,s,67, 77,96-s,19)$ = $C_{48}(s,19,29)$. 
 
 $\Rightarrow$ $C_{96}(5,s,43)$ $\cong$ $C_{96}(s,19,29)$ for $s$ = 6, 18, 30, 42. Also, for $s$ = 6, 18, 30, 42, we get,
\\ 
5$Ad_{96}(C_{96}(5,s,43)$ = $Ad_{96}(C_{96}(5,s,43, 53,96-s,91))$ 
  = $\{C_{96}(x(5,s,43, 53,96-s,91)): x = 1$, 
$5,7,11,13,17,19,23,25,29,31,35,37,41,43,47,49,53,55,59,61,65,67,71,77,79,83,85,89,91\}$. 
 
 From the above relation, we find $Ad_{96}(C_{96}(5,s,43))$ for $s$ = 6, 18, 30, 42 and will show that $C_{96}(s,19,29) \notin Ad_{96}(C_{96}(5,s,43))$.
 
\item [\rm (w1)]  $Ad_{96}(C_{96}(5,6,43))$ = $\{C_{96}(5,6,43)$, $C_{96}(23,25,30)$, $C_{96}(13,35,42)$, $C_{96}(7,30,41)$, 
 
 \hfill  $C_{96}(17,18,31)$, $C_{96}(6,11,37)$, $C_{96}(1,18,47)$, $C_{96}(19,29,42)\}$.  

 $\Rightarrow$ $C_{96}(6,19,29) \notin Ad_{96}(C_{96}(5,6,43))$.
 
 $\Rightarrow$ $C_{96}(5,6,43)$ and $C_{96}(6,19,29)$ are Type-2 isomorphic w.r.t. $m$ = 2.
 
\item [\rm (w2)]  $Ad_{96}(C_{96}(5,18,43))$ = $\{C_{96}(5,18,43)$, $C_{96}(6,23,25)$, $C_{96}(13,30,35)$, $C_{96}(6,7,41)$, 
 
 \hfill  $C_{96}(17,31,42)$, $C_{96}(11,18,37)$, $C_{96}(1,42,47)$, $C_{96}(19,29,30)\}$.  

 $\Rightarrow$ $C_{96}(18,19,29) \notin Ad_{96}(C_{96}(5,18,43))$.
 
 $\Rightarrow$ $C_{96}(5,18,43)$ and $C_{96}(18,19,29)$ are Type-2 isomorphic w.r.t. $m$ = 2.
 
\item [\rm (w3)]  $Ad_{96}(C_{96}(5,30,43))$ = $\{C_{96}(5,30,43)$, $C_{96}(23,25,42)$, $C_{96}(13,18,35)$, $C_{96}(7,41,42)$, 
 
 \hfill  $C_{96}(6,17,31)$, $C_{96}(11,30,37)$, $C_{96}(1,6,47)$, $C_{96}(18,19,29)\}$.  

 $\Rightarrow$ $C_{96}(19,29,30) \notin Ad_{96}(C_{96}(5,30,43))$.
 
 $\Rightarrow$ $C_{96}(5,30,43)$ and $C_{96}(19,29,30)$ are Type-2 isomorphic w.r.t. $m$ = 2.
 
\item [\rm (w4)]  $Ad_{96}(C_{96}(5,42,43))$ = $\{C_{96}(5,42,43)$, $C_{96}(18,23,25)$, $C_{9}(6,13,35)$, $C_{96}(7,18,41)$, 
 
 \hfill  $C_{96}(17,30,31)$, $C_{96}(11,37,42)$, $C_{96}(1,30,47)$, $C_{96}(6,19,29)\}$.  

 $\Rightarrow$ $C_{96}(19,29,42) \notin Ad_{96}(C_{96}(5,42,43))$.
 
 $\Rightarrow$ $C_{96}(5,42,43)$ and $C_{96}(19,29,42)$ are Type-2 isomorphic w.r.t. $m$ = 2.
 
\item [\rm (x)]  For  $s$ = 6, 18, 30, 42, using the definition of $\theta_{n,m,t}$, we get,  
 \\
 $\theta_{96,2,12}(C_{96}(7,s,41))$ = $\theta_{96,2,12}(C_{96}(7,s,41, 55,96-s,89))$ = $C_{96}(\theta_{96,2,12}(7,s,41, 55,96-s,89))$ 
 
 \hspace{2.7cm} = $C_{96}(31,s,65, 79,96-s,17)$ = $C_{48}(s,17,31)$. 
 
 $\Rightarrow$ $C_{96}(7,s,41)$ $\cong$ $C_{96}(s,17,31)$ for $s$ = 6, 18, 30, 42. Also, for $s$ = 6, 18, 30, 42, we get,
\\ 
5$Ad_{96}(C_{96}(7,s,41)$ = $Ad_{96}(C_{96}(7,s,41, 55,96-s,89))$ 
  = $\{C_{96}(x(7,s,41, 55,96-s,89)): x = 1$, 
$5,7,11,13,17,19,23,25,29,31,35,37,41,43,47,49,53,55,59,61,65,67,71,77,79,83,85,89,91\}$. 
 
 From the above relation, we find $Ad_{96}(C_{96}(7,s,41))$ for $s$ = 6, 18, 30, 42 and will show that $C_{96}(s,17,31) \notin Ad_{96}(C_{96}(7,s,41))$.
 
\item [\rm (x1)]  $Ad_{96}(C_{96}(6,7,41))$ = $\{C_{96}(6,7,41)$, $C_{96}(13,30,35)$, $C_{96}(1,42,47)$, $C_{96}(19,29,30)$,  
 
 \hfill $C_{96}(5,18,43)$, $C_{96}(6,23,25)$, $C_{96}(11,18,37)$, $C_{96}(17,31,42)\}$.  

 $\Rightarrow$ $C_{96}(6,17,31) \notin Ad_{96}(C_{96}(6,7,41))$.
 
 $\Rightarrow$ $C_{96}(6,7,41)$ and $C_{96}(6,17,31)$ are Type-2 isomorphic w.r.t. $m$ = 2.
 
\item [\rm (x2)]  $Ad_{96}(C_{96}(7,18,41))$ = $\{C_{96}(7,18,41)$, $C_{96}(6,13,35)$, $C_{96}(1,30,47)$, $C_{96}(6,19,29)$, 
 
 \hfill  $C_{96}(5,42,43)$, $C_{96}(18,23,25)$, $C_{96}(11,37,42)$, $C_{96}(17,30,31)\}$.  

 $\Rightarrow$ $C_{96}(18,17,31) \notin Ad_{96}(C_{96}(7,18,41))$.
 
 $\Rightarrow$ $C_{96}(7,18,41)$ and $C_{96}(18,17,31)$ are Type-2 isomorphic w.r.t. $m$ = 2.
 
\item [\rm (x3)]  $Ad_{96}(C_{96}(7,30,41))$ = $\{C_{96}(7,30,41)$, $C_{96}(13,35,42)$, $C_{96}(1,18,47)$, $C_{96}(19,29,42)$, 
 
 \hfill  $C_{96}(5,6,43)$, $C_{96}(23,25,30)$, $C_{96}(6,11,37)$, $C_{96}(17,18,31)\}$.  

 $\Rightarrow$ $C_{96}(17,30,31) \notin Ad_{96}(C_{96}(7,30,41))$.
 
 $\Rightarrow$ $C_{96}(7,30,41)$ and $C_{96}(17,30,31)$ are Type-2 isomorphic w.r.t. $m$ = 2.
 
\item [\rm (x4)]  $Ad_{96}(C_{96}(7,41,42))$ = $\{C_{96}(7,41,42)$, $C_{96}(13,18,35)$, $C_{9}(1,6,47)$, $C_{96}(18,19,29)$, 
 
 \hfill  $C_{96}(5,30,43)$, $C_{96}(23,25,42)$, $C_{96}(11,30,37)$, $C_{96}(6,17,31)\}$.  

 $\Rightarrow$ $C_{96}(17,31,42) \notin Ad_{96}(C_{96}(7,41,42))$.
 
 $\Rightarrow$ $C_{96}(7,41,42)$ and $C_{96}(17,31,42)$ are Type-2 isomorphic w.r.t. $m$ = 2.
 
\item [\rm (y)] When $s$ = 6,18,30,42, $C_{96}(9,s,39)$ = $3.C_{32}(3,t,13)$ and $C_{96}(s,15,33)$ = $3.C_{32}(t,5,11)$ where $s$ = $3t$. For $t$ = 2,6,10,14, $C_{32}(3,t,13)$ and $C_{32}(t,5,11)$ are Type-2 isomorphic w.r.t. $m$ = 2 by the following. $t$ = 2,6,10,14,

$\theta_{32,2,4}(C_{32}(3,t,13))$ = $\theta_{32,2,4}(C_{32}(3,t,13, 19,32-t,29))$ = $C_{32}(\theta_{32,2,4}(3,t,13, 19,32-t,29))$ 
 
 \hspace{3cm} = $C_{32}(11,t,21, 27,32-t,5)$ = $C_{32}(t,5,11)$.

$\Rightarrow$  $C_{32}(3,t,13)$ $\cong$ $C_{32}(t,5,11)$ for $t$ = 2,6,10,14. Also,
 
 $Ad_{32}(C_{32}(2,3,13))$ = $\{C_{32}(2,3,13)$, $C_{32}(6,7,9)$, $C_{32}(1,10,15)$, $C_{32}(5,11,14)\}$;

 $Ad_{32}(C_{32}(3,6,13))$ = $\{C_{32}(3,6,13)$, $C_{32}(7,9,14)$, $C_{32}(1,2,15)$, $C_{32}(5,10,11)\}$;

 $Ad_{32}(C_{32}(3,10,13))$ = $\{C_{32}(3,10,13)$, $C_{32}(2,7,9)$, $C_{32}(1,14,15)$, $C_{32}(5,6,11)\}$; and

 $Ad_{32}(C_{32}(3,14,13))$ = $\{C_{32}(3,14,13)$, $C_{32}(7,9,10)$, $C_{32}(1,6,15)$, $C_{32}(2,5,11)\}$.

$\Rightarrow$  $C_{32}(t,5,11)\notin Ad_{32}(C_{32}(3,t,13))$ for $t$ = 2,6,10,14. 

$\Rightarrow$  $C_{32}(3,t,13)$ and $C_{32}(t,5,11)$ are Type-2 isomorphic w.r.t. $m$ = 2  for $t$ = 2,6,10,14 where $\gcd(32, t)$ = 2 = $m$. 

This implies, $3.C_{32}(3,t,13)$ = $C_{96}(9,3t,39)$ and $3.C_{32}(t,5,11)$ = $C_{96}(3t,15,33)$ are Type-2 isomorphic w.r.t. $m$ = 2 for $t$ = 2,6,10,14 where $\gcd(32, t)$ = 2 = $m$. 

\item [\rm (z)]  For  $s$ = 6, 18, 30, 42, using the definition of $\theta_{n,m,t}$, we get,  
 \\
 $\theta_{96,2,12}(C_{96}(11,s,37))$ = $\theta_{96,2,12}(C_{96}(11,s,37, 59,96-s,85))$ 

\hfill = $C_{96}(\theta_{96,2,12}(11,s,37, 59,96-s,85))$ = $C_{96}(35,s,61, 83,96-s,13)$ = $C_{48}(s,13,35)$. 
 
 $\Rightarrow$ $C_{96}(11,s,37)$ $\cong$ $C_{96}(s,13,35)$ for $s$ = 6, 18, 30, 42. Also, for $s$ = 6, 18, 30, 42, we get,
\\ 
$Ad_{96}(C_{96}(11,s,37)$ = $Ad_{96}(C_{96}(11,s,37, 59,96-s,85))$ 
  = $\{C_{96}(x(11,s,37, 59,96-s,85)): x = 1$, 
$5,7,11,13,17,19,23,25,29,31,35,37,41,43,47,49,53,55,59,61,65,67,71,77,79,83,85,89,91\}$. 
 
 From the above relation, we find $Ad_{96}(C_{96}(11,s,37))$ for $s$ = 6, 18, 30, 42 and will show that $C_{96}(s,13,35) \notin Ad_{96}(C_{96}(11,s,37))$.
 
\item [\rm (z1)]  $Ad_{96}(C_{96}(6,11,37))$ = $\{C_{96}(6,11,37)$, $C_{96}(7,30,41)$, $C_{96}(19,29,42)$, $C_{96}(23,25,30)$,  
 
 \hfill $C_{96}(1,18,47)$, $C_{96}(5,6,43)$, $C_{96}(17,18,31)$, $C_{96}(13,35,42)\}$.  

 $\Rightarrow$ $C_{96}(6,13,35) \notin Ad_{96}(C_{96}(6,11,37))$.
 
 $\Rightarrow$ $C_{96}(6,11,37)$ and $C_{96}(6,13,35)$ are Type-2 isomorphic w.r.t. $m$ = 2.
 
\item [\rm (z2)]  $Ad_{96}(C_{96}(11,18,37))$ = $\{C_{96}(11,18,37)$, $C_{96}(6,7,41)$, $C_{96}(19,29,30)$, $C_{96}(6,23,25)$, 
 
 \hfill  $C_{96}(1,42,47)$, $C_{96}(5,18,43)$, $C_{96}(17,31,42)$, $C_{96}(13,30,35)\}$.  

 $\Rightarrow$ $C_{96}(13,18,35) \notin Ad_{96}(C_{96}(11,18,37))$.
 
 $\Rightarrow$ $C_{96}(11,18,37)$ and $C_{96}(13,18,35)$ are Type-2 isomorphic w.r.t. $m$ = 2.
 
\item [\rm (z3)]  $Ad_{96}(C_{96}(1,30,37))$ = $\{C_{96}(11,30,37)$, $C_{96}(7,41,42)$, $C_{96}(18,19,29)$, $C_{96}(23,25,42)$, 
 
 \hfill  $C_{96}(1,6,47)$, $C_{96}(5,30,43)$, $C_{96}(6,17,31)$, $C_{96}(13,18,35)\}$.  

 $\Rightarrow$ $C_{96}(13,30,35) \notin Ad_{96}(C_{96}(11,30,37))$.
 
 $\Rightarrow$ $C_{96}(11,30,37)$ and $C_{96}(13,30,35)$ are Type-2 isomorphic w.r.t. $m$ = 2.
 
\item [\rm (z4)]  $Ad_{96}(C_{96}(11,37,42))$ = $\{C_{96}(13,35,42)$, $C_{96}(7,18,41)$, $C_{9}(6,19,29)$, $C_{96}(18,23,25)$, 
 
 \hfill  $C_{96}(1,30,47)$, $C_{96}(5,42,43)$, $C_{96}(17,30,31)$, $C_{96}(6,13,35)\}$.  

 $\Rightarrow$ $C_{96}(13,35,42) \notin Ad_{96}(C_{96}(11,37,42))$.
 
 $\Rightarrow$ $C_{96}(11,37,42)$ and $C_{96}(13,35,42)$ are Type-2 isomorphic w.r.t. $m$ = 2. \hfill $\Box$
\end{enumerate} 

\medskip
\noindent
{\bf Type-1 sets of $C_{96}(r_1,r_2,r_3)$ related to problems \ref{p7} to \ref{p9}}

\medskip
\noindent
We present here the following  Type-1 sets which are related to problems \ref{p7} to \ref{p9}.
\begin{enumerate}
\item [\rm (a1)]  $Ad_{96}(C_{96}(1,4,47))$ = $\{C_{96}(1,4,47)$, $C_{96}(5,20,43)$, $C_{96}(7,28,41)$, $C_{96}(11,37,44)$, 

\hfill $C_{96}(13,35,44)$, $C_{96}(17,28,31)$, $C_{96}(19,20,29)$, $C_{96}((4,23,25))\}$;  

\item [\rm (a2)]  $Ad_{96}(C_{96}(1,8,47))$ = $\{C_{96}(1,8,47)$, $C_{96}(5,40,43)$, $C_{96}(7,40,41)$, $C_{96}(8,11,37)$, 

\hfill $C_{96}(8,13,35)$, $C_{96}(17,31,40)$, $C_{96}(19,29,40)$, $C_{96}((8,23,25))\}$;  

\item [\rm (a3)]  $Ad_{96}(C_{96}(1,12,47))$ = $\{C_{96}(1,12,47)$, $C_{96}(5,36,43)$, $C_{96}(7,12,41)$, $C_{96}(11,36,37)$, 

\hfill $C_{96}(13,35,36)$, $C_{96}(12,17,31)$, $C_{96}(19,29,36)$, $C_{96}((12,23,25))\}$;  

\item [\rm (a4)]  $Ad_{96}(C_{96}(1,16,47))$ = $\{C_{96}(1,16,47)$, $C_{96}(5,16,43)$, $C_{96}(7,16,41)$, $C_{96}(11,16,37)$, 

\hfill $C_{96}(13,16,35)$, $C_{96}(16,17,31)$, $C_{96}(16,19,29)$, $C_{96}((16,23,25))\}$;  

\item [\rm (a5)]  $Ad_{96}(C_{96}(1,20,47))$ = $\{C_{96}(1,20,47)$, $C_{96}(4,5,43)$, $C_{96}(7,41,44)$, $C_{96}(11,28,37)$, 

\hfill $C_{96}(13,28,35)$, $C_{96}(17,31,44)$, $C_{96}(4,19,29)$, $C_{96}((20,23,25))\}$;  

\item [\rm (a6)]  $Ad_{96}(C_{96}(1,24,47))$ = $\{C_{96}(1,24,47)$, $C_{96}(5,24,43)$, $C_{96}(7,24,41)$, $C_{96}(11,24,37)$, 

\hfill $C_{96}(13,24,35)$, $C_{96}(17,24,31)$, $C_{96}(19,24,29)$, $C_{96}((23,24,25))\}$;  

\item [\rm (a7)]  $Ad_{96}(C_{96}(1,28,47))$ = $\{C_{96}(1,28,47)$, $C_{96}(5,43,44)$, $C_{96}(4,7,41)$, $C_{96}(11,20,37)$, 

\hfill $C_{96}(13,20,35)$, $C_{96}(4,17,31)$, $C_{96}(19,29,44)$, $C_{96}((23,25,28))\}$;  

\item [\rm (a8)]  $Ad_{96}(C_{96}(1,32,47))$ = $\{C_{96}(1,32,47)$, $C_{96}(5,32,43)$, $C_{96}(7,32,41)$, $C_{96}(11,32,37)$, 

\hfill $C_{96}(13,32,35)$, $C_{96}(17,32,31)$, $C_{96}(19,32,29)$, $C_{96}((23,25,32))\}$;  

\item [\rm (a9)]  $Ad_{96}(C_{96}(1,36,47))$ = $\{C_{96}(1,36,47)$, $C_{96}(5,12,43)$, $C_{96}(7,36,41)$, $C_{96}(11,12,37)$, 

\hfill $C_{96}(12,13,35)$, $C_{96}(17,31,36)$, $C_{96}(12,19,29)$, $C_{96}((23,25,36))\}$;  

\item [\rm (a10)]  $Ad_{96}(C_{96}(1,40,47))$ = $\{C_{96}(1,40,47)$, $C_{96}(5,8,43)$, $C_{96}(7,8,41)$, $C_{96}(11,37,40)$, 

\hfill $C_{96}(13,35,40)$, $C_{96}(8,17,31)$, $C_{96}(8,19,29)$, $C_{96}((23,25,40))\}$;  

\item [\rm (a11)]  $Ad_{96}(C_{96}(1,44,47))$ = $\{C_{96}(1,44,47)$, $C_{96}(5,28,43)$, $C_{96}(7,20,41)$, $C_{96}(4,11,37)$, 

\hfill $C_{96}(4,13,35)$, $C_{96}(17,20,31)$, $C_{96}(19,28,29)$, $C_{96}((23,25,44))\}$;  

\item [\rm (a12)]  $Ad_{96}(C_{96}(1,47,48))$ = $\{C_{96}(1,47,48)$, $C_{96}(5,43,48)$, $C_{96}(7,41,48)$, $C_{96}(11,37,48)$, 

\hfill $C_{96}(13,35,48)$, $C_{96}(17,31,48)$, $C_{96}(19,29,48)$, $C_{96}((23,25,48))\}$;  

\item [\rm (b1)]  $Ad_{96}(C_{96}(3,4,45))$ = $\{C_{96}(3,4,45)$, $C_{96}(15,20,33)$, $C_{96}(21,27,28)$, $C_{96}(15,33,44)$, 

\hfill $C_{96}(9,39,44)$, $C_{96}(3,28,45)$, $C_{96}(9,20,39)$, $C_{48}(4,21,27)\}$;  

\item [\rm (b2)]  $Ad_{96}(C_{96}(3,8,45))$ = $\{C_{96}(3,8,45)$,  $C_{96}(15,33,40)$, $C_{96}(21,27,40)$, $C_{96}(8,15,33)$, 

\hfill $C_{96}(8,9,39)$, $C_{96}(3,40,45)$, $C_{96}(9,29,40)$, $C_{48}(8,21,27)\}$;  

\item [\rm (b3)]  $Ad_{96}(C_{96}(3,12,45))$ = $\{C_{96}(3,12,45))$,  $C_{96}(15,33,36)$, $C_{96}(12,21,27)$, $C_{96}(15,33,36)$, 

\hfill $C_{96}(9,36,39)$, $C_{96}(3,12,45)$, $C_{96}(9,29,36)$, $C_{48}(12,21,27)\}$;  

\item [\rm (b4)]  $Ad_{96}(C_{96}(3,16,45))$ = $\{C_{96}(3,16,45)$,  $C_{96}(15,16,33)$, $C_{96}(16,21,27)$, $C_{96}(15,16,33)$, 

\hfill $C_{96}(9,16,39)$, $C_{96}(3,16,45)$, $C_{96}(9,16,29)$, $C_{48}(16,21,27)\}$;  

\item [\rm (b5)]  $Ad_{96}(C_{96}(3,20,45))$ = $\{C_{96}(3,20,45)$,  $C_{96}(4,15,33)$, $C_{96}(21,27,44)$, $C_{96}(15,28,33)$, 

\hfill $C_{96}(9,28,39)$, $C_{96}(3,44,45)$, $C_{96}(4,9,29)$, $C_{48}(20,21,27)\}$;  

\item [\rm (b6)]  $Ad_{96}(C_{96}(3,24,45))$ = $\{C_{96}(3,24,45)$,  $C_{96}(15,24,33)$, $C_{96}(21,24,27)$, $C_{96}(15,24,33)$, 

\hfill $C_{96}(9,24,39)$, $C_{96}(3,24,45)$, $C_{96}(9,24,29)$, $C_{48}(21,24,27)\}$;  

\item [\rm (b7)]  $Ad_{96}(C_{96}(3,28,45))$ = $\{C_{96}(3,28,45)$,  $C_{96}(15,33,44)$, $C_{96}(4,21,27)$, $C_{96}(15,20,44)$, 

\hfill $C_{96}(9,20,39)$, $C_{96}(3,4,45)$, $C_{96}(9,29,44)$, $C_{48}(21,28,27)\}$;  

\item [\rm (b8)]  $Ad_{96}(C_{96}(3,32,45))$ = $\{C_{96}(3,32,45)$,  $C_{96}(15,32,33)$, $C_{96}(21,27,32)$, $C_{96}(15,32,32)$, 

\hfill $C_{96}(9,32,39)$, $C_{96}(3,32,45)$, $C_{96}(9,32,29)$, $C_{48}(21,27,32)\}$;  

\item [\rm (b9)]  $Ad_{96}(C_{96}(3,36,45))$ = $\{C_{96}(3,36,45)$,  $C_{96}(12,15,33)$, $C_{96}(21,27,36)$, $C_{96}(12,15,33)$, 

\hfill $C_{96}(9,12,39)$, $C_{96}(3,36,45)$, $C_{96}(9,12,29)$, $C_{48}(21,27,36)\}$;  

\item [\rm (b10)]  $Ad_{96}(C_{96}(3,40,45))$ = $\{C_{96}(3,40,45)$,  $C_{96}(8,15,33)$, $C_{96}(8,21,27)$, $C_{96}(15,33,40)$, 

\hfill $C_{96}(9,39,40)$, $C_{96}(3,8,45)$, $C_{96}(8,9,29)$, $C_{48}(40,21,27)\}$;  

\item [\rm (b11)]  $Ad_{96}(C_{96}(3,44,45))$ = $\{C_{96}(3,44,45)$,  $C_{96}(15,28,33)$, $C_{96}(20,21,27)$, $C_{96}(4,15,33)$, 

\hfill $C_{96}(4,9,39)$, $C_{96}(3,20,45)$, $C_{96}(9,28,29)$, $C_{48}(21,27,44)\}$;  

\item [\rm (b12)]  $Ad_{96}(C_{96}(3,48,45))$ = $\{C_{96}(3,48,45)$,  $C_{96}(15,33,48)$, $C_{96}(21,27,48)$, $C_{96}(15,33,48)$, 

\hfill $C_{96}(9,39,48)$, $C_{96}(3,45,48)$, $C_{96}(9,29,48)$, $C_{48}(21,27,48)\}$;  

item [\rm (c1)]  $Ad_{96}(C_{96}(4,5,43))$ = $\{C_{96}(4,5,43)$, $C_{96}(20,23,25)$, $C_{96}(13,28,35)$, $C_{96}(7,41,44)$, 

\hfill $C_{96}(17,31,44)$, $C_{96}(11,28,37)$, $C_{96}(1,20,47)$, $C_{48}(4,19,29)\}$;  

\item [\rm (c2)]  $Ad_{96}(C_{96}(5,8,43))$ = $\{C_{96}(5,8,43)$,  $C_{96}(23,25,40)$, $C_{96}(13,35,40)$, $C_{96}(7,8,41)$, 

\hfill $C_{96}(8,17,31)$, $C_{96}(11,37,40)$, $C_{96}(1,40,47)$, $C_{48}(8,19,29)\}$;  

\item [\rm (c3)]  $Ad_{96}(C_{96}(5,12,43))$ = $\{C_{96}(5,12,43)$,  $C_{96}(23,25,36)$, $C_{96}(12,13,35)$, $C_{96}(7,36,41)$, 

\hfill $C_{96}(17,31,36)$, $C_{96}(11,12,37)$, $C_{96}(1,36,47)$, $C_{48}(12,19,29)\}$;  

\item [\rm (c4)]  $Ad_{96}(C_{96}(5,16,43))$ = $\{C_{96}(5,16,43)$,  $C_{96}(16,23,25)$, $C_{96}(13,16,35)$, $C_{96}(7,16,41)$, 

\hfill $C_{96}(16,17,31)$, $C_{96}(11,16,37)$, $C_{96}(1,16,47)$, $C_{48}(16,19,29)\}$;  

\item [\rm (c5)]  $Ad_{96}(C_{96}(5,20,43))$ = $\{C_{96}(5,20,43)$,  $C_{96}(4,23,25)$, $C_{96}(13,35,44)$, $C_{96}(7,28,41)$, 

\hfill $C_{96}(17,28,31)$, $C_{96}(11,37,44)$, $C_{96}(1,4,47)$, $C_{48}(19,20,29)\}$;  

\item [\rm (c6)]  $Ad_{96}(C_{96}(5,24,43))$ = $\{C_{96}(5,24,43)$,  $C_{96}(23,24,25)$, $C_{96}(13,24,35)$, $C_{96}(7,24,41)$, 

\hfill $C_{96}(17,24,31)$, $C_{96}(11,24,37)$, $C_{96}(1,24,47)$, $C_{48}(19,24,29)\}$;  

\item [\rm (c7)]  $Ad_{96}(C_{96}(5,28,43))$ = $\{C_{96}(5,28,43)$,  $C_{96}(23,25,44)$, $C_{96}(4,13,35)$, $C_{96}(7,20,41)$, 

\hfill $C_{96}(17,20,31)$, $C_{96}(4,11,37)$, $C_{96}(1,44,47)$, $C_{48}(19,28,29)\}$;  

\item [\rm (c8)]  $Ad_{96}(C_{96}(5,32,43))$ = $\{C_{96}(5,32,43)$,  $C_{96}(23,25,32)$, $C_{96}(13,32,35)$, $C_{96}(7,32,41)$, 

\hfill $C_{96}(17,31,32)$, $C_{96}(11,32,37)$, $C_{96}(1,32,47)$, $C_{48}(19,29,32)\}$;  

\item [\rm (c9)]  $Ad_{96}(C_{96}(5,36,43))$ = $\{C_{96}(5,36,43)$,  $C_{96}(12,23,25)$, $C_{96}(13,35,36)$, $C_{96}(7,12,41)$, 

\hfill $C_{96}(12,17,31)$, $C_{96}(11,36,37)$, $C_{96}(1,12,47)$, $C_{48}(19,29,36)\}$;  

\item [\rm (c10)]  $Ad_{96}(C_{96}(5,40,43))$ = $\{C_{96}(5,40,43)$,  $C_{96}(8,23,25)$, $C_{96}(8,13,35)$, $C_{96}(7,40,41)$, 

\hfill $C_{96}(17,31,40)$, $C_{96}(8,11,37)$, $C_{96}(1,8,47)$, $C_{48}(19,29,40)\}$;  

\item [\rm (c11)]  $Ad_{96}(C_{96}(5,43,44))$ = $\{C_{96}(5,43,44)$,  $C_{96}(23,25,28)$, $C_{96}(13,20,35)$, $C_{96}(4,7,41)$, 

\hfill $C_{96}(4,17,31)$, $C_{96}(11,20,37)$, $C_{96}(1,28,47)$, $C_{48}(19,29,44)\}$;  

\item [\rm (c12)]  $Ad_{96}(C_{96}(5,43,48))$ = $\{C_{96}(5,43,48)$,  $C_{96}(23,25,48)$, $C_{96}(13,35,48)$, $C_{96}(7,41,48)$, 

\hfill $C_{96}(17,31,48)$, $C_{96}(11,37,48)$, $C_{96}(1,47,48)$, $C_{48}(19,29,48)\}$;  

\item [\rm (d1)]  $Ad_{96}(C_{96}(4,7,41))$ = $\{C_{96}(4,7,41)$, $C_{96}(13,20,35)$, $C_{96}(1,28,47)$, $C_{96}(19,29,44)$, 

\hfill $C_{96}(5,43,44)$, $C_{96}(23,25,28)$, $C_{96}(11,20,37)$, $C_{48}(4,17,31)\}$;  

\item [\rm (d2)]  $Ad_{96}(C_{96}(7,8,41))$ = $\{C_{96}(7,8,41)$,  $C_{96}(13,35,40)$, $C_{96}(1,40,47)$, $C_{96}(8,19,29)$, 

\hfill $C_{96}(5,8,43)$, $C_{96}(23,25,40)$, $C_{96}(11,37,40)$, $C_{48}(8,17,31)\}$;  

\item [\rm (d3)]  $Ad_{96}(C_{96}(7,12,41))$ = $\{C_{96}(7,12,41)$,  $C_{96}(13,35,36)$, $C_{96}(1,12,47)$, $C_{96}(19,29,36)$, 

\hfill $C_{96}(5,36,43)$, $C_{96}(12,23,25)$, $C_{96}(11,36,37)$, $C_{48}(12,17,31)\}$;  

\item [\rm (d4)]  $Ad_{96}(C_{96}(7,16,41))$ = $\{C_{96}(7,16,41)$,  $C_{96}(13,16,35)$, $C_{96}(1,16,47)$, $C_{96}(16,19,29)$, 

\hfill $C_{96}(5,16,43)$, $C_{96}(16,23,25)$, $C_{96}(11,16,37)$, $C_{48}(16,17,31)\}$;  

\item [\rm (d5)]  $Ad_{96}(C_{96}(7,20,41))$ = $\{C_{96}(7,20,41)$,  $C_{96}(4,13,35)$, $C_{96}(1,44,47)$, $C_{96}(19,28,29)$, 

\hfill $C_{96}(5,28,43)$, $C_{96}(23,25,44)$, $C_{96}(4,11,37)$, $C_{48}(17,20,31)\}$;  

\item [\rm (d6)]  $Ad_{96}(C_{96}(7,24,41))$ = $\{C_{96}(7,24,41)$,  $C_{96}(13,24,35)$, $C_{96}(1,24,47)$, $C_{96}(19,24,29)$, 

\hfill $C_{96}(5,24,43)$, $C_{96}(23,24,25)$, $C_{96}(11,24,37)$, $C_{48}(24,17,31)\}$;  

\item [\rm (d7)]  $Ad_{96}(C_{96}(7,28,41))$ = $\{C_{96}(7,28,41)$,  $C_{96}(13,35,44)$, $C_{96}(1,4,47)$, $C_{96}(19,20,29)$, 

\hfill $C_{96}(5,20,43)$, $C_{96}(4,23,25)$, $C_{96}(11,37,44)$, $C_{48}(17,28,31)\}$;  

\item [\rm (d8)]  $Ad_{96}(C_{96}(7,32,41))$ = $\{C_{96}(7,32,41)$,  $C_{96}(13,32,35)$, $C_{96}(1,32,47)$, $C_{96}(19,29,32)$, 

\hfill $C_{96}(5,32,43)$, $C_{96}(23,25,32)$, $C_{96}(11,32,37)$, $C_{48}(17,31,32)\}$;  

\item [\rm (d9)]  $Ad_{96}(C_{96}(7,36,41))$ = $\{C_{96}(7,36,41)$,  $C_{96}(13,12,35)$, $C_{96}(1,36,47)$, $C_{96}(12,19,29)$, 

\hfill $C_{96}(5,12,43)$, $C_{96}(23,25,36)$, $C_{96}(11,12,37)$, $C_{48}(17,31,36)\}$;  

\item [\rm (d10)]  $Ad_{96}(C_{96}(7,40,41))$ = $\{C_{96}(7,40,41)$,  $C_{96}(8,13,35)$, $C_{96}(1,8,47)$, $C_{96}(19,40,44)$, 

\hfill $C_{96}(5,40,43)$, $C_{96}(8,23,25)$, $C_{96}(8,11,37)$, $C_{48}(17,31,40)\}$;  

\item [\rm (d11)]  $Ad_{96}(C_{96}(7,41,44))$ = $\{C_{96}(7,41,44)$,  $C_{96}(13,28,35)$, $C_{96}(1,20,47)$, $C_{96}(4,19,29)$, 

\hfill $C_{96}(4,5,43)$, $C_{96}(20,23,25)$, $C_{96}(11,28,37)$, $C_{48}(17,31,44)\}$;  

\item [\rm (d12)]  $Ad_{96}(C_{96}(7,41,48))$ = $\{C_{96}(7,41,48)$,  $C_{96}(13,35,48)$, $C_{96}(1,47,48)$, $C_{96}(19,29,48)$, 

\hfill $C_{96}(5,43,48)$, $C_{96}(23,25,48)$, $C_{96}(11,37,48)$, $C_{48}(17,31,48)\}$;  

\item [\rm (e1)]  $Ad_{96}(C_{96}(4,9,39))$ = $\{C_{96}(4,9,39)$, $C_{96}(3,20,45)$, $C_{96}(15,28,33)$, $C_{96}(3,44,45)$, 

\hfill  $C_{96}(21,27,44)$, $C_{96}(9,28,39)$, $C_{96}(20,21,27)$, $C_{48}(4,15,33)\}$;  

\item [\rm (e2)]  $Ad_{96}(C_{96}(8,9,39))$ = $\{C_{96}(8,9,39)$, $C_{96}(3,40,45)$, $C_{96}(15,33,40)$, $C_{96}(3,8,45)$, 

\hfill  $C_{96}(8,21,27)$, $C_{96}(9,39,40)$, $C_{96}(21,27,40)$, $C_{48}(8,15,33)\}$;  

\item [\rm (e3)]  $Ad_{96}(C_{96}(9,12,39))$ = $\{C_{96}(9,12,39)$, $C_{96}(3,36,45)$, $C_{96}(12,15,33)$, $C_{96}(3,36,45)$, 

\hfill  $C_{96}(21,27,36)$, $C_{96}(9,12,39)$, $C_{96}(21,27,36)$, $C_{48}(12,15,33)\}$;  

\item [\rm (e4)]  $Ad_{96}(C_{96}(9,16,39))$ = $\{C_{96}(9,16,39)$, $C_{96}(3,16,45)$, $C_{96}(15,16,33)$, $C_{96}(3,16,45)$, 

\hfill  $C_{96}(16,21,27)$, $C_{96}(9,16,39)$, $C_{96}(16,21,27)$, $C_{48}(15,16,33)\}$;  

\item [\rm (e5)]  $Ad_{96}(C_{96}(9,20,39))$ = $\{C_{96}(9,20,39)$, $C_{96}(3,4,45)$, $C_{96}(15,33,44)$, $C_{96}(3,28,45)$, 

\hfill  $C_{96}(21,27,28)$, $C_{96}(9,39,44)$, $C_{96}(4,21,27)$, $C_{48}(15,20,33)\}$;  

\item [\rm (e6)]  $Ad_{96}(C_{96}(9,24,39))$ = $\{C_{96}(9,24,39)$, $C_{96}(3,24,45)$, $C_{96}(15,24,33)$, $C_{96}(3,24,45)$, 

\hfill  $C_{96}(21,24,27)$, $C_{96}(9,24,39)$, $C_{96}(21,24,27)$, $C_{48}(15,24,33)\}$;  

\item [\rm (e7)]  $Ad_{96}(C_{96}(9,28,39))$ = $\{C_{96}(9,28,39)$, $C_{96}(3,44,45)$, $C_{96}(4,15,33)$, $C_{96}(3,20,45)$, 

\hfill  $C_{96}(20,21,27)$, $C_{96}(4,9,39)$, $C_{96}(21,27,44)$, $C_{48}(15,28,33)\}$;  

\item [\rm (e8)]  $Ad_{96}(C_{96}(9,32,39))$ = $\{C_{96}(9,32,39)$, $C_{96}(3,32,45)$, $C_{96}(15,32,33)$, $C_{96}(3,32,45)$, 

\hfill  $C_{96}(21,27,32)$, $C_{96}(9,32,39)$, $C_{96}(21,27,32)$, $C_{48}(15,32,33)\}$;  

\item [\rm (e9)]  $Ad_{96}(C_{96}(9,36,39))$ = $\{C_{96}(9,36,39)$, $C_{96}(3,12,45)$, $C_{96}(15,36,33)$, $C_{96}(3,12,45)$, 

\hfill  $C_{96}(12,21,27)$, $C_{96}(9,36,39)$, $C_{96}(12,21,27)$, $C_{48}(15,33,36)\}$;  

\item [\rm (e10)]  $Ad_{96}(C_{96}(9,39,40))$ = $\{C_{96}(9,39,40)$, $C_{96}(3,8,45)$, $C_{96}(8,15,33)$, $C_{96}(3,40,45)$, 

\hfill  $C_{96}(21,27,40)$, $C_{96}(8,9,39)$, $C_{96}(8,21,27)$, $C_{48}(15,33,40)\}$;  

\item [\rm (e11)]  $Ad_{96}(C_{96}(9,39,44))$ = $\{C_{96}(9,39,44)$, $C_{96}(3,28,45)$, $C_{96}(15,20,33)$, $C_{96}(3,4,45)$, 

\hfill  $C_{96}(4,21,27)$, $C_{96}(9,20,39)$, $C_{96}(21,28,27)$,  $C_{48}(15,33,44)\}$;  

\item [\rm (e12)]  $Ad_{96}(C_{96}(9,39,48))$ = $\{C_{96}(9,39,48)$, $C_{96}(3,45,48)$, $C_{96}(15,33,48)$, $C_{96}(3,45,48)$, 

\hfill  $C_{96}(21,27,48)$, $C_{96}(9,39,48)$, $C_{96}(21,27,48)$, $C_{48}(15,33,48)\}$;  

\item [\rm (f1)]  $Ad_{96}(C_{96}(4,11,37))$ = $\{C_{96}(4,11,37)$, $C_{96}(7,20,41)$,  $C_{96}(19,28,29)$,  $C_{96}(23,25,44)$, 

\hfill  $C_{96}(1,44,47)$,  $C_{96}(5,28,43)$,  $C_{96}(17,20,31)$, $C_{48}(4,13,35)\}$;  

\item [\rm (f2)]  $Ad_{96}(C_{96}(8,11,37))$ = $\{C_{96}(8,11,37)$,  $C_{96}(7,40,41)$,  $C_{96}(19,29,40)$,  $C_{96}(8,23,25)$, 

\hfill  $C_{96}(1,8,47)$,  $C_{96}(5,40,43)$,  $C_{96}(17,31,40)$, $C_{48}(8,13,35)\}$;  

\item [\rm (f3)]  $Ad_{96}(C_{96}(11,12,37))$ = $\{C_{96}(11,12,37)$,  $C_{96}(7,36,41)$,  $C_{96}(12,19,29)$,  $C_{96}(23,25,36)$, 

\hfill  $C_{96}(1,36,47)$,  $C_{96}(5,12,43)$,  $C_{96}(17,31,36)$, $C_{48}(12,13,35)\}$;  

\item [\rm (f4)]  $Ad_{96}(C_{96}(11,16,37))$ = $\{C_{96}(11,16,37)$,  $C_{96}(7,16,41)$,  $C_{96}(16,19,29)$,  $C_{96}(16,23,25)$, 

\hfill  $C_{96}(1,16,47)$,  $C_{96}(5,16,43)$,  $C_{96}(16,17,31)$, $C_{48}(13,16,35)\}$;  

\item [\rm (f5)]  $Ad_{96}(C_{96}(11,20,37))$ = $\{C_{96}(11,20,37)$,  $C_{96}(4,7,41)$,  $C_{96}(19,29,44)$,  $C_{96}(23,25,28)$, 

\hfill  $C_{96}(1,28,47)$,  $C_{96}(5,43,44)$,  $C_{96}(4,17,31)$, $C_{48}(13,20,35)\}$;  

\item [\rm (f6)]  $Ad_{96}(C_{96}(11,24,37))$ = $\{C_{96}(11,24,37)$,  $C_{96}(7,24,41)$,  $C_{96}(19,24,29)$,  $C_{96}(23,24,25)$, 

\hfill  $C_{96}(1,24,47)$,  $C_{96}(5,24,43)$,  $C_{96}(17,24,31)$, $C_{48}(13,24,35)\}$;  

\item [\rm (f7)]  $Ad_{96}(C_{96}(11,28,37))$ = $\{C_{96}(11,28,37)$,  $C_{96}(7,41,44)$,  $C_{96}(4,19,29)$,  $C_{96}(20,23,25)$, 

\hfill  $C_{96}(1,20,47)$,  $C_{96}(4,5,43)$,  $C_{96}(17,31,44)$, $C_{48}(13,28,35)\}$;  

\item [\rm (f8)]  $Ad_{96}(C_{96}(11,32,37))$ = $\{C_{96}(11,32,37)$,  $C_{96}(7,32,41)$,  $C_{96}(19,29,32)$,  $C_{96}(23,25,32)$, 

\hfill  $C_{96}(1,32,47)$,  $C_{96}(5,32,43)$,  $C_{96}(17,31,32)$, $C_{48}(13,32,35)\}$;  

\item [\rm (f9)]  $Ad_{96}(C_{96}(11,36,37))$ = $\{C_{96}(11,36,37)$,  $C_{96}(7,12,41)$,  $C_{96}(19,29,36)$,  $C_{96}(12,23,25)$, 

\hfill  $C_{96}(1,12,47)$,  $C_{96}(5,36,43)$,  $C_{96}(12,17,31)$, $C_{48}(13,35,36)\}$;  

 \item [\rm (f10)]  $Ad_{96}(C_{96}(11,37,40))$ = $\{C_{96}(11,37,40)$,  $C_{96}(7,8,41)$,  $C_{96}(8,19,29)$,  $C_{96}(23,25,40)$, 

\hfill  $C_{96}(1,40,47)$,  $C_{96}(5,8,43)$,  $C_{96}(8,17,31)$, $C_{48}(13,35,40)\}$;  

 \item [\rm (f11)]  $Ad_{96}(C_{96}(11,37,44))$ = $\{C_{96}(11,37,44)$,  $C_{96}(7,28,41)$,  $C_{96}(19,20,29)$,  $C_{96}(4,23,25)$, 

\hfill  $C_{96}(1,4,47)$,  $C_{96}(5,20,43)$,  $C_{96}(17,28,31)$, $C_{48}(13,35,44)\}$;  

\item [\rm (f12)]  $Ad_{96}(C_{96}(11,37,48))$ = $\{C_{96}(11,37,48)$,  $C_{96}(7,41,48)$,  $C_{96}(19,29,48)$,  $C_{96}(23,25,48)$, 

\hfill  $C_{96}(1,47,48)$,  $C_{96}(5,43,48)$,  $C_{96}(17,31,48)$, $C_{48}(13,35,48)\}$.  

\item [\rm (g1)]  $Ad_{96}(C_{96}(1,2,47))$ = $\{C_{96}(1,2,47)$, $C_{96}(5,10,43)$, $C_{96}(7,14,41)$, $C_{96}(11,22,37)$, 
 
 \hfill  $C_{96}(13,26,35)$, $C_{96}(17,31,34)$, $C_{96}(19,29,38)$, $C_{96}(23,25,46)\}$.  

\item [\rm (g2)]  $Ad_{96}(C_{96}(1,10,47))$ = $\{C_{96}(1,10,47)$, $C_{96}(5,43,46)$, $C_{96}(7,26,41)$, $C_{96}(11,14,37)$, 
 
 \hfill  $C_{96}(13,34,35)$, $C_{96}(17,22,31)$, $C_{96}(2,19,29)$, $C_{96}(23,25,38)\}$.  
 
\item [\rm (g3)]  $Ad_{96}(C_{96}(1,14,47))$ = $\{C_{96}(1,14,47)$, $C_{96}(5,26,43,)$, $C_{96}(2,7,41)$, $C_{96}(11,37,38)$, 
 
 \hfill  $C_{96}(10,13,35)$, $C_{96}(17,31,46)$, $C_{96}(19,22,29)$, $C_{96}(23,25,34)\}$.  
 
\item [\rm (g4)]  $Ad_{96}(C_{96}(1,22,47))$ = $\{C_{96}(1,22,47)$, $C_{96}(5,14,43,)$, $C_{96}(7,38,41)$, $C_{96}(11,37,46)$, 
 
 \hfill  $C_{96}(2,13,35)$, $C_{96}(10,17,31)$, $C_{96}(19,29,34)$, $C_{96}(23,25,26)\}$.  
 
\item [\rm (g5)]  $Ad_{96}(C_{96}(1,26,47))$ = $\{C_{96}(1,26,47)$, $C_{96}(5,34,43)$, $C_{96}(7,10,41)$, $C_{96}(2,11,37)$, 
 
 \hfill  $C_{96}(13,35,46)$, $C_{96}(17,31,38)$, $C_{96}(14,19,29)$, $C_{96}(22,23,25)\}$.  
 
\item [\rm (g6)]  $Ad_{96}(C_{96}(1,34,47))$ = $\{C_{96}(1,34,47)$, $C_{96}(5,22,43)$, $C_{96}(7,41,46)$, $C_{96}(10,11,37)$, 
 
 \hfill  $C_{96}(13,35,38)$, $C_{96}(2,17,31)$, $C_{96}(19,26,29)$, $C_{96}(14,23,25)\}$.  
 
\item [\rm (g7)]  $Ad_{96}(C_{96}(1,38,47))$ = $\{C_{96}(1,38,47)$, $C_{96}(2,5,43)$, $C_{96}(7,22,41)$, $C_{96}(11,34,37)$, 
 
 \hfill  $C_{96}(13,14,35)$, $C_{96}(17,26,31)$, $C_{96}(19,29,46)$, $C_{96}(10,23,25)\}$.  
 
\item [\rm (g8)]  $Ad_{96}(C_{96}(1,46,47))$ = $\{C_{96}(1,46,47)$, $C_{96}(5,38,43)$, $C_{96}(7,34,41)$, $C_{96}(11,26,37)$, 
 
 \hfill  $C_{96}(13,22,35)$, $C_{96}(14,17,31)$, $C_{96}(10,19,29)$, $C_{96}(2,23,25)\}$.  
 
\item [\rm (h1)]  $Ad_{96}(C_{96}(2,3,45))$ = $\{C_{96}(2,3,45)$, $C_{96}(10,15,33)$, $C_{96}(14,21,27)$, $C_{96}(15,22,33)$,  
 
 \hfill  $C_{96}(9,26,39)$, $C_{96}(3,34,45)$, $C_{96}(9,38,39)$, $C_{96}(21,27,46)\}$.  

 \item [\rm (h2)]  $Ad_{96}(C_{96}(3,10,45))$ = $\{C_{96}(3,10,45)$, $C_{96}(15,33,46)$, $C_{96}(21,26,27)$, $C_{96}(14,15,33)$,  
 
 \hfill  $C_{96}(9,34,39)$, $C_{96}(3,22,45)$, $C_{96}(2,9,39)$, $C_{96}(21,27,38)\}$.  

 \item [\rm (h3)]  $Ad_{96}(C_{96}(3,14,45))$ = $\{C_{96}(3,14,45)$, $C_{96}(15,26,33,)$, $C_{96}(2,21,27)$, $C_{96}(15,33,38)$,  
 
 \hfill  $C_{96}(9,10,39)$, $C_{96}(3,45,46)$, $C_{96}(9,22,39)$, $C_{96}(21,27,34)\}$.  

 \item [\rm (h4)]  $Ad_{96}(C_{96}(3,22,45))$ = $\{C_{96}(3,22,45)$, $C_{96}(14,15,33)$, $C_{96}(21,27,38)$, $C_{96}(15,33,46)$,  
 
 \hfill  $C_{96}(2,9,39)$, $C_{96}(3,10,45)$, $C_{96}(9,34,39)$, $C_{96}(21,26,27)\}$.  

\item [\rm (h5)]  $Ad_{96}(C_{96}(3,26,45))$ = $\{C_{96}(3,26,45)$, $C_{96}(15,33,34)$, $C_{96}(10,21,27)$, $C_{96}(2,15,33)$,  
 
 \hfill  $C_{96}(9,39,46)$, $C_{96}(3,38,45)$, $C_{96}(9,14,39)$, $C_{96}(21,22,27)\}$.  

\item [\rm (h6)]  $Ad_{96}(C_{96}(3,34,45))$ = $\{C_{96}(3,34,45)$, $C_{96}(15,22,33)$, $C_{96}(21,27,46)$, $C_{96}(10,15,33)$,  
 
 \hfill  $C_{96}(9,38,39,)$, $C_{96}(2,3,45)$, $C_{96}(9,26,39)$, $C_{96}(14,21,27)\}$.  

\item [\rm (h7)]  $Ad_{96}(C_{96}(3,38,45))$ = $\{C_{96}(3,38,45)$, $C_{96}(2,15,33)$, $C_{96}(21,22,27)$, $C_{96}(15,33,34)$,  
 
 \hfill  $C_{96}(9,14,39,)$, $C_{96}(3,26,45)$, $C_{96}(9,39,46)$, $C_{96}(10,21,27)\}$.  

\item [\rm (h8)]  $Ad_{96}(C_{96}(3,45,46))$ = $\{C_{96}(3,45,46)$, $C_{96}(15,33,38)$, $C_{96}(21,27,34)$, $C_{96}(15,26,33)$,  
 
 \hfill  $C_{96}(9,22,39,)$, $C_{96}(3,14,45)$, $C_{96}(9,10,39)$, $C_{96}(2,21,27)\}$.  

\item [\rm (i1)]  $Ad_{96}(C_{96}(2,5,43))$ = $\{C_{96}(2,5,43)$, $C_{96}(10,23,25)$, $C_{96}(13,14,35)$, $C_{96}(7,22,41)$,  
 
 \hfill  $C_{96}(17,26,31)$, $C_{96}(11,34,37)$, $C_{96}(1,38,47)$, $C_{96}(19,29,46)\}$.  

\item [\rm (i2)]  $Ad_{96}(C_{96}(5,10,43))$ = $\{C_{96}(5,10,43)$, $C_{96}(23,25,46)$, $C_{96}(13,26,35)$, $C_{96}(7,14,41)$,  
 
 \hfill  $C_{96}(17,31,34)$, $C_{96}(11,22,37)$, $C_{96}(1,2,47)$, $C_{96}(19,29,38)\}$.  

\item [\rm (i3)]  $Ad_{96}(C_{96}(5,14,43))$ = $\{C_{96}(5,14,43)$, $C_{96}(23,25,26)$, $C_{96}(2,13,35)$, $C_{96}(7,38,41)$,  
 
 \hfill  $C_{96}(10,17,31)$, $C_{96}(11,37,46)$, $C_{96}(1,22,47)$, $C_{96}(19,29,34)\}$.  
 
\item [\rm (i4)]  $Ad_{96}(C_{96}(5,22,43))$ = $\{C_{96}(5,22,43)$, $C_{96}(14,23,25)$, $C_{96}(13,35,38)$, $C_{96}(7,41,46)$,  
 
 \hfill  $C_{96}(2,17,31)$, $C_{96}(10,11,37)$, $C_{96}(1,34,47)$, $C_{96}(19,26,29)\}$.  

\item [\rm (i5)]  $Ad_{96}(C_{96}(5,26,43))$ = $\{C_{96}(5,26,43)$, $C_{96}(23,25,34)$, $C_{96}(10,13,35)$, $C_{96}(2,7,41)$,  
 
 \hfill  $C_{96}(17,31,46)$, $C_{96}(11,37,38)$, $C_{96}(1,14,47)$, $C_{96}(19,22,29)\}$.  

\item [\rm (i6)]  $Ad_{96}(C_{96}(5,34,43))$ = $\{C_{96}(5,34,43)$, $C_{96}(22,23,25)$, $C_{96}(13,35,46)$, $C_{96}(7,10,41)$,  
 
 \hfill  $C_{96}(17,31,38)$, $C_{96}(2,11,37)$, $C_{96}(1,26,47)$, $C_{96}(14,19,29)\}$.  

\item [\rm (i7)]  $Ad_{96}(C_{96}(5,38,43))$ = $\{C_{96}(5,38,43)$, $C_{96}(2,23,25)$, $C_{96}(13,22,35)$, $C_{96}(7,34,41)$,  
 
 \hfill  $C_{96}(14,17,31)$, $C_{96}(11,26,37)$, $C_{96}(1,46,47)$, $C_{96}(10,19,29)\}$.  

\item [\rm (i8)]  $Ad_{96}(C_{96}(5,43,46))$ = $\{C_{96}(5,43,46)$, $C_{96}(,23,25)$, $C_{96}(13,,35)$, $C_{96}(7,,41)$,  
 
 \hfill  $C_{96}(,17,31)$, $C_{96}(11,,37)$, $C_{96}(1,,47)$, $C_{96}(,19,29)\}$.  

\item [\rm (j1)]  $Ad_{96}(C_{96}(2,7,41))$ = $\{C_{96}(2,7,41)$, $C_{96}(10,13,35)$, $C_{96}(1,14,47)$, $C_{96}(19,22,29)$,  
 
 \hfill  $C_{96}(5,26,43)$, $C_{96}(23,25,34)$, $C_{96}(11,37,38)$, $C_{96}(17,31,46)\}$.  

\item [\rm (j2)]  $Ad_{96}(C_{96}(7,10,41))$ = $\{C_{96}(7,10,41)$, $C_{96}(13,35,46)$, $C_{96}(1,26,47)$, $C_{96}(14,19,29)$,  
 
 \hfill  $C_{96}(5,34,43)$, $C_{96}(22,23,25)$, $C_{96}(2,11,37)$, $C_{96}(17,31,38)\}$.  

\item [\rm (j3)]  $Ad_{96}(C_{96}(7,14,41))$ = $\{C_{96}(7,14,41)$, $C_{96}(13,26,35)$, $C_{96}(1,2,47)$, $C_{96}(19,29,38)$,  
 
 \hfill  $C_{96}(5,10,43)$, $C_{96}(23,25,46)$, $C_{96}(11,22,37)$, $C_{96}(17,31,34)\}$.  

\item [\rm (j4)]  $Ad_{96}(C_{96}(7,22,41))$ = $\{C_{96}(7,22,41)$, $C_{96}(13,14,35)$, $C_{96}(1,38,47)$, $C_{96}(19,29,46)$,  
 
 \hfill  $C_{96}(2,5,43)$, $C_{96}(10,23,25)$, $C_{96}(11,34,37)$, $C_{96}(17,26,31)\}$.  

\item [\rm (j5)]  $Ad_{96}(C_{96}(7,26,41))$ = $\{C_{96}(7,26,41)$, $C_{96}(13,34,35)$, $C_{96}(1,10,47)$, $C_{96}(2,19,29)$,  
 
 \hfill  $C_{96}(5,43,46)$, $C_{96}(23,25,38)$, $C_{96}(11,14,37)$, $C_{96}(17,22,31)\}$.  

\item [\rm (j6)]  $Ad_{96}(C_{96}(7,34,41))$ = $\{C_{96}(7,34,41)$, $C_{96}(13,22,35)$, $C_{96}(1,46,47)$, $C_{96}(10,19,29)$,  
 
 \hfill  $C_{96}(5,38,43)$, $C_{96}(2,23,25)$, $C_{96}(11,26,37)$, $C_{96}(14,17,31)\}$.  

\item [\rm (j7)]  $Ad_{96}(C_{96}(7,38,41))$ = $\{C_{96}(7,38,41)$, $C_{96}(2,13,35)$, $C_{96}(1,22,47)$, $C_{96}(19,29,34)$,  
 
 \hfill  $C_{96}(5,14,43)$, $C_{96}(23,25,26)$, $C_{96}(11,37,46)$, $C_{96}(10,17,31)\}$.  

\item [\rm (j8)]  $Ad_{96}(C_{96}(7,41,46))$ = $\{C_{96}(7,41,46)$, $C_{96}(13,35,38)$, $C_{96}(1,34,47)$, $C_{96}(19,26,29)$,  
 
 \hfill  $C_{96}(5,22,43)$, $C_{96}(14,23,25)$, $C_{96}(10,11,37)$, $C_{96}(2,17,31)\}$.  

\item [\rm (k1)]  $Ad_{96}(C_{96}(2,9,39))$ = $\{C_{96}(2,9,39)$, $C_{96}(3,10,45)$, $C_{96}(14,15,33)$, $C_{96}(3,22,45)$,  
 
 \hfill  $C_{96}(21,26,27)$, $C_{96}(9,34,39)$, $C_{96}(21,27,38)$, $C_{96}(15,33,46)\}$.  

\item [\rm (k2)]  $Ad_{96}(C_{96}(9,10,39))$ = $\{C_{96}(9,10,39)$, $C_{96}(3,45,46)$, $C_{96}(15,26,33)$, $C_{96}(3,14,45)$,  
 
 \hfill  $C_{96}(21,27,34)$, $C_{96}(9,22,39)$, $C_{96}(2,21,27)$, $C_{96}(15,33,38)\}$.  

\item [\rm (k3)]  $Ad_{96}(C_{96}(9,14,39))$ = $\{C_{96}(9,14,39)$, $C_{96}(3,26,45)$, $C_{96}(2,15,33)$, $C_{96}(3,38,45)$,  
 
 \hfill  $C_{96}(10,21,27)$, $C_{96}(9,39,46)$, $C_{96}(21,22,27)$, $C_{96}(15,33,34)\}$.  

\item [\rm (k4)]  $Ad_{96}(C_{96}(9,22,39))$ = $\{C_{96}(9,22,39)$, $C_{96}(3,14,45)$, $C_{96}(15,33,38)$, $C_{96}(3,45,46)$,  
 
 \hfill  $C_{96}(2,21,27)$, $C_{96}(9,10,39)$, $C_{96}(21,27,34)$, $C_{96}(15,26,33)\}$.  

\item [\rm (k5)]  $Ad_{96}(C_{96}(9,26,39))$ = $\{C_{96}(9,26,39)$, $C_{96}(3,34,45)$, $C_{96}(10,15,33)$, $C_{96}(2,3,45)$,  
 
 \hfill  $C_{96}(21,27,46)$, $C_{96}(9,38,39)$, $C_{96}(14,21,27)$, $C_{96}(15,22,33)\}$.  

\item [\rm (k6)]  $Ad_{96}(C_{96}(9,34,39))$ = $\{C_{96}(9,34,39)$, $C_{96}(3,22,45)$, $C_{96}(15,33,46)$, $C_{96}(3,10,45)$,  
 
 \hfill  $C_{96}(21,27,38)$, $C_{96}(2,9,39)$, $C_{96}(21,26,27)$, $C_{96}(14,15,33)\}$.  

\item [\rm (k7)]  $Ad_{96}(C_{96}(9,38,39))$ = $\{C_{96}(9,38,39)$, $C_{96}(2,3,45)$, $C_{96}(15,22,33)$, $C_{96}(3,34,45)$,  
 
 \hfill  $C_{96}(14,21,27)$, $C_{96}(9,26,39)$, $C_{96}(21,27,46)$, $C_{96}(10,15,33)\}$.  

\item [\rm (k8)]  $Ad_{96}(C_{96}(9,39,46))$ = $\{C_{96}(9,39,46)$, $C_{96}(3,38,45)$, $C_{96}(15,33,34)$, $C_{96}(3,26,45)$,  
 
 \hfill  $C_{96}(21,22,27)$, $C_{96}(9,14,39)$, $C_{96}(10,21,27)$, $C_{96}(2,15,33)\}$.  

\item [\rm (l1)]  $Ad_{96}(C_{96}(2,11,37))$ = $\{C_{96}(2,11,37)$, $C_{96}(7,10,41)$, $C_{96}(14,19,29)$, $C_{96}(22,23,25)$,  
 
 \hfill  $C_{96}(1,26,47)$, $C_{96}(5,34,43)$, $C_{96}(17,31,38)$, $C_{96}(13,35,46)\}$.  

\item [\rm (l2)]  $Ad_{96}(C_{96}(10,11,37))$ = $\{C_{96}(10,11,37)$, $C_{96}(7,41,46)$, $C_{96}(19,26,29)$, $C_{96}(14,23,25)$,  
 
 \hfill  $C_{96}(1,34,47)$, $C_{96}(5,22,43)$, $C_{96}(2,17,31)$, $C_{96}(13,35,38)\}$.  

\item [\rm (l3)]  $Ad_{96}(C_{96}(11,14,37))$ = $\{C_{96}(11,14,37)$, $C_{96}(7,26,41)$, $C_{96}(2,19,29)$, $C_{96}(23,25,38)$,  
 
 \hfill  $C_{96}(1,10,47)$, $C_{96}(5,43,46)$, $C_{96}(17,22,31)$, $C_{96}(13,34,35)\}$.  

\item [\rm (l4)]  $Ad_{96}(C_{96}(11,22,37))$ = $\{C_{96}(11,22,37)$, $C_{96}(7,14,41)$, $C_{96}(19,29,38)$, $C_{96}(23,25,46)$,  
 
 \hfill  $C_{96}(1,2,47)$, $C_{96}(5,10,43)$, $C_{96}(17,31,34)$, $C_{96}(13,26,35)\}$.  

\item [\rm (l5)]  $Ad_{96}(C_{96}(11,26,37))$ = $\{C_{96}(11,26,37)$, $C_{96}(7,34,41)$, $C_{96}(10,19,29)$, $C_{96}(2,23,25)$,  
 
 \hfill  $C_{96}(1,46,47)$, $C_{96}(5,38,43)$, $C_{96}(14,17,31)$, $C_{96}(13,22,35)\}$.  

\item [\rm (l6)]  $Ad_{96}(C_{96}(11,34,37))$ = $\{C_{96}(11,34,37)$, $C_{96}(7,22,41)$, $C_{96}(19,29,46)$, $C_{96}(10,23,25)$,  
 
 \hfill  $C_{96}(1,38,47)$, $C_{96}(2,5,43)$, $C_{96}(17,26,31)$, $C_{96}(13,14,35)\}$.  

\item [\rm (l7)]  $Ad_{96}(C_{96}(11,37,38))$ = $\{C_{96}(11,37,38)$, $C_{96}(2,7,41)$, $C_{96}(19,22,29)$, $C_{96}(23,25,34)$,  
 
 \hfill  $C_{96}(1,14,47)$, $C_{96}(5,26,43)$, $C_{96}(17,31,46)$, $C_{96}(10,13,35)\}$.  

\item [\rm (l8)]  $Ad_{96}(C_{96}(11,37,46))$ = $\{C_{96}(11,37,46)$, $C_{96}(7,38,41)$, $C_{96}(19,29,34)$, $C_{96}(23,25,26)$,  
 
 \hfill  $C_{96}(1,22,47)$, $C_{96}(5,14,43)$, $C_{96}(10,17,31)$, $C_{96}(2,13,35)\}$.  

\item [\rm (u1)]  $Ad_{96}(C_{96}(1,6,47))$ = $\{C_{96}(1,6,47)$, $C_{96}(5,30,43)$, $C_{96}(7,41,42)$, $C_{96}(11,30,37)$, 
 
 \hfill  $C_{96}(13,18,35)$, $C_{96}(6,17,31)$, $C_{96}(18,19,29)$, $C_{96}(23,25,42)\}$.  

\item [\rm (u2)]  $Ad_{96}(C_{96}(1,18,47))$ = $\{C_{96}(1,18,47)$, $C_{96}(5,6,43)$, $C_{96}(7,30,41)$, $C_{96}(6,11,37)$, 
 
 \hfill  $C_{96}(13,35,42)$, $C_{96}(17,18,31)$, $C_{96}(19,29,42)$, $C_{96}(23,25,30)\}$.  

\item [\rm (u3)]  $Ad_{96}(C_{96}(1,30,47))$ = $\{C_{96}(1,30,47)$, $C_{96}(5,42,43)$, $C_{96}(7,18,41)$, $C_{96}(11,37,42)$, 
 
 \hfill  $C_{96}(6,13,35)$, $C_{96}(17,30,31)$, $C_{96}(6,19,29)$, $C_{96}(18,23,25)\}$.  

\item [\rm (u4)]  $Ad_{96}(C_{96}(1,42,47))$ = $\{C_{96}(1,42,47)$, $C_{96}(5,18,43)$, $C_{96}(6,7,41)$, $C_{96}(11,18,37)$, 
 
 \hfill  $C_{96}(13,30,35)$, $C_{96}(17,31,42)$, $C_{96}(19,29,30)$, $C_{96}(6,23,25)\}$.  

\item [\rm (v1)]   $Ad_{32}(C_{32}(1,2,15))$ = $\{C_{32}(1,2,15)$, $C_{32}(3,6,13)$, $C_{32}(5,10,11)$, $C_{32}(7,9,14)\}$;

\item [\rm (v2)]   $Ad_{32}(C_{32}(1,6,15))$ = $\{C_{32}(1,6,15)$, $C_{32}(3,13,14)$, $C_{32}(2,5,11)$, $C_{32}(7,9,10)\}$;

\item [\rm (v3)]   $Ad_{32}(C_{32}(1,10,15))$ = $\{C_{32}(1,10,15)$, $C_{32}(2,3,13)$, $C_{32}(5,11,14)$, $C_{32}(6,7,9)\}$; 

\item [\rm (v4)]   $Ad_{32}(C_{32}(1,14,15))$ = $\{C_{32}(1,14,15)$, $C_{32}(3,10,13)$, $C_{32}(5,6,11)$, $C_{32}(2,7,9)\}$.

\item [\rm (w1)]  $Ad_{96}(C_{96}(5,6,43))$ = $\{C_{96}(5,6,43)$, $C_{96}(23,25,30)$, $C_{96}(13,35,42)$, $C_{96}(7,30,41)$, 
 
 \hfill  $C_{96}(17,18,31)$, $C_{96}(6,11,37)$, $C_{96}(1,18,47)$, $C_{96}(19,29,42)\}$.  

\item [\rm (w2)]  $Ad_{96}(C_{96}(5,18,43))$ = $\{C_{96}(5,18,43)$, $C_{96}(6,23,25)$, $C_{96}(13,30,35)$, $C_{96}(6,7,41)$, 
 
 \hfill  $C_{96}(17,31,42)$, $C_{96}(11,18,37)$, $C_{96}(1,42,47)$, $C_{96}(19,29,30)\}$.  

\item [\rm (w3)]  $Ad_{96}(C_{96}(5,30,43))$ = $\{C_{96}(5,30,43)$, $C_{96}(23,25,42)$, $C_{96}(13,18,35)$, $C_{96}(7,41,42)$, 
 
 \hfill  $C_{96}(6,17,31)$, $C_{96}(11,30,37)$, $C_{96}(1,6,47)$, $C_{96}(18,19,29)\}$.  

\item [\rm (w4)]  $Ad_{96}(C_{96}(5,42,43))$ = $\{C_{96}(5,42,43)$, $C_{96}(18,23,25)$, $C_{9}(6,13,35)$, $C_{96}(7,18,41)$, 
 
 \hfill  $C_{96}(17,30,31)$, $C_{96}(11,37,42)$, $C_{96}(1,30,47)$, $C_{96}(6,19,29)\}$.  

\item [\rm (x1)]  $Ad_{96}(C_{96}(6,7,41))$ = $\{C_{96}(6,7,41)$, $C_{96}(13,30,35)$, $C_{96}(1,42,47)$, $C_{96}(19,29,30)$,  
 
 \hfill $C_{96}(5,18,43)$, $C_{96}(6,23,25)$, $C_{96}(11,18,37)$, $C_{96}(17,31,42)\}$.  

\item [\rm (x2)]  $Ad_{96}(C_{96}(7,18,41))$ = $\{C_{96}(7,18,41)$, $C_{96}(6,13,35)$, $C_{96}(1,30,47)$, $C_{96}(6,19,29)$, 
 
 \hfill  $C_{96}(5,42,43)$, $C_{96}(18,23,25)$, $C_{96}(11,37,42)$, $C_{96}(17,30,31)\}$.  

\item [\rm (x3)]  $Ad_{96}(C_{96}(7,30,41))$ = $\{C_{96}(7,30,41)$, $C_{96}(13,35,42)$, $C_{96}(1,18,47)$, $C_{96}(19,29,42)$, 
 
 \hfill  $C_{96}(5,6,43)$, $C_{96}(23,25,30)$, $C_{96}(6,11,37)$, $C_{96}(17,18,31)\}$.  

\item [\rm (x4)]  $Ad_{96}(C_{96}(7,41,42))$ = $\{C_{96}(7,41,42)$, $C_{96}(13,18,35)$, $C_{9}(1,6,47)$, $C_{96}(18,19,29)$, 
 
 \hfill  $C_{96}(5,30,43)$, $C_{96}(23,25,42)$, $C_{96}(11,30,37)$, $C_{96}(6,17,31)\}$.  

\item [\rm (y1)]   $Ad_{32}(C_{32}(2,3,13))$ = $\{C_{32}(2,3,13)$, $C_{32}(6,7,9)$, $C_{32}(1,10,15)$, $C_{32}(5,11,14)\}$;

\item [\rm (y2)]   $Ad_{32}(C_{32}(3,6,13))$ = $\{C_{32}(3,6,13)$, $C_{32}(7,9,14)$, $C_{32}(1,2,15)$, $C_{32}(5,10,11)\}$;

\item [\rm (y3)]   $Ad_{32}(C_{32}(3,10,13))$ = $\{C_{32}(3,10,13)$, $C_{32}(2,7,9)$, $C_{32}(1,14,15)$, $C_{32}(5,6,11)\}$; 

\item [\rm (y4)]   $Ad_{32}(C_{32}(3,14,13))$ = $\{C_{32}(3,14,13)$, $C_{32}(7,9,10)$, $C_{32}(1,6,15)$, $C_{32}(2,5,11)\}$.

\item [\rm (z1)]  $Ad_{96}(C_{96}(6,11,37))$ = $\{C_{96}(6,11,37)$, $C_{96}(7,30,41)$, $C_{96}(19,29,42)$, $C_{96}(23,25,30)$,  
 
 \hfill $C_{96}(1,18,47)$, $C_{96}(5,6,43)$, $C_{96}(17,18,31)$, $C_{96}(13,35,42)\}$.  

\item [\rm (z2)]  $Ad_{96}(C_{96}(11,18,37))$ = $\{C_{96}(11,18,37)$, $C_{96}(6,7,41)$, $C_{96}(19,29,30)$, $C_{96}(6,23,25)$, 
 
 \hfill  $C_{96}(1,42,47)$, $C_{96}(5,18,43)$, $C_{96}(17,31,42)$, $C_{96}(13,30,35)\}$.  

\item [\rm (z3)]  $Ad_{96}(C_{96}(1,30,37))$ = $\{C_{96}(11,30,37)$, $C_{96}(7,41,42)$, $C_{96}(18,19,29)$, $C_{96}(23,25,42)$, 
 
 \hfill  $C_{96}(1,6,47)$, $C_{96}(5,30,43)$, $C_{96}(6,17,31)$, $C_{96}(13,18,35)\}$.  

\item [\rm (z4)]  $Ad_{96}(C_{96}(11,37,42))$ = $\{C_{96}(13,35,42)$, $C_{96}(7,18,41)$, $C_{9}(6,19,29)$, $C_{96}(18,23,25)$, 
 
 \hfill  $C_{96}(1,30,47)$, $C_{96}(5,42,43)$, $C_{96}(17,30,31)$, $C_{96}(6,13,35)\}$.  
\end{enumerate} 

\section{Conclusion}

 In \cite{v2-2} - \cite{v2-4}, the author established that the number of pairs of Type-2 isomorphic circulant graphs of orders 16, 24 and 32 are 8, 32 and 384; and the number of triples of Type-2 isomorphic circulant graphs of order 27 is 12 and presented all. The author feels that a lot of scope is there for further research and proposes the following open problems on this topic.

Here, we propose the following open problems on circulant graphs of order 48.

\begin{oprm} \label{op1} {\rm The following pairs of circulant graphs are non-isomorphic for $s$ = 3, 9, 15, 21.
\begin{enumerate}	\item [\rm (a)]  $C_{48}(1,s,23)$ and $C_{48}(s,11,13)$;  and 
			\item [\rm (b)]  $C_{48}(5,s,19)$ and $C_{48}(s,7,17)$. \hfill $\Box$ 
\end{enumerate} }	
\end{oprm}

\begin{oprm} \label{op2} {\rm Find all pairs of isomorphic circulant graphs of order 48. \hfill $\Box$}	
\end{oprm}

\begin{oprm} \label{op3} {\rm Find all pairs of Type-1 isomorphic circulant graphs of order 48. \hfill $\Box$}	
\end{oprm}

\begin{oprm} \label{op4} {\rm Find all pairs of Type-2 isomorphic circulant graphs of order 48. \hfill $\Box$}	
\end{oprm}
 
\begin{oprm} \label{op5} {\rm Find all pairs of isomorphic circulant graphs of order 81. \hfill $\Box$}	
\end{oprm}

\begin{oprm} \label{op6} {\rm Find all pairs of Type-1 isomorphic circulant graphs of order 81. \hfill $\Box$}	
\end{oprm}

\begin{oprm} \label{op7} {\rm Find all pairs of Type-2 isomorphic circulant graphs of order 81. \hfill $\Box$}	
\end{oprm}

\begin{oprm} \label{op8} {\rm Find all pairs of isomorphic circulant graphs of order 96. \hfill $\Box$}	
\end{oprm}

\begin{oprm} \label{op9} {\rm Find all pairs of Type-1 isomorphic circulant graphs of order 96. \hfill $\Box$}	
\end{oprm}

\begin{oprm} \label{op10} {\rm Find all pairs of Type-2 isomorphic circulant graphs of order 96. \hfill $\Box$}	
\end{oprm}

\vspace{.1cm}
\noindent
\textbf{Declaration of competing interest}\quad The author declares that he has no conflict of interest.


\end{document}